\documentclass[11pt, reqno]{amsart}

\usepackage[OT2, T1]{fontenc}

\usepackage[usenames,dvipsnames]{color}

\usepackage{titletoc}

\usepackage{url}
\usepackage{amsmath}
\usepackage{array}
\usepackage{graphicx}
\usepackage{amsfonts}
\usepackage{amssymb}
\usepackage{amsthm}
\definecolor{link}{RGB}{11,0,128}
\usepackage[colorlinks=true,citecolor=NavyBlue,linkcolor=Bittersweet,urlcolor=link]{hyperref}
\usepackage{paralist}
\usepackage{colonequals}		
\usepackage{mathabx}		
\usepackage{amscd}
\usepackage[all,cmtip]{xy}	
\usepackage{sseq}			
\usepackage{verbatim}
\usepackage{parskip}
\usepackage{microtype}
\usepackage[in]{fullpage}
\usepackage{mathrsfs} 			
\usepackage{cleveref}
\usepackage{enumitem}
\usepackage{moreenum}			
\usepackage[alphabetic,lite]{amsrefs} 	
\usepackage{stmaryrd}	
\usepackage{subfiles}
\usepackage{ragged2e}

\DeclareSymbolFont{cyrletters}{OT2}{wncyr}{m}{n}
\DeclareMathSymbol{\Sha}{\mathalpha}{cyrletters}{"58}

\newcommand{\gA}{\alpha}
\newcommand{\gB}{\beta}
\newcommand{\gG}{\gamma}

\newcommand{\bA}{\mathbb{A}}
\newcommand{\bC}{\mathbb{C}}

\newcommand{\bF}{\mathbb{F}}
\newcommand{\bG}{\mathbb{G}}

\newcommand{\bP}{\mathbb{P}}
\newcommand{\bQ}{\mathbb{Q}}

\newcommand{\bZ}{\mathbb{Z}}

\newcommand{\cE}{\mathcal{E}}
\newcommand{\cF}{\mathcal{F}}
\newcommand{\cG}{\mathcal{G}}

\newcommand{\cI}{\mathcal{I}}

\newcommand{\cM}{\mathcal{M}}

\newcommand{\fm}{\mathfrak{m}}

\newcommand{\sA}{\mathscr{A}}
\newcommand{\sB}{\mathscr{B}}
\newcommand{\sC}{\mathscr{C}}

\newcommand{\sG}{\mathscr{G}}

\newcommand{\sI}{\mathscr{I}}

\newcommand{\sO}{\mathscr{O}}

\newcommand{\sU}{\mathscr{U}}
\newcommand{\sV}{\mathscr{V}}

\newcommand{\sX}{\mathscr{X}}
\newcommand{\sY}{\mathscr{Y}}

\newcommand{\ra}{\rightarrow}

\newcommand{\xra}{\xrightarrow}

\newcommand{\hra}{\hookrightarrow}

\newcommand{\wt}{\widetilde}
\newcommand{\wh}{\widehat}

\newcommand{\pr}{^{\prime}}

\newcommand{\ce}{\colonequals}
\newcommand{\ov}{\overline}

\newcommand{\sm}{\mathrm{sm}}
\newcommand{\sing}{\mathrm{sing}}
\renewcommand{\b}{\textbf}
\newcommand{\surjects}{\twoheadrightarrow}

\newcommand{\tensor}{\otimes} 		
\newcommand{\isomto}{\overset{\sim}{\longrightarrow}}

\newcommand{\st}{{\mathrm{st}}}		
\newcommand{\llb}{\llbracket}		
\newcommand{\rrb}{\rrbracket}		

\renewcommand{\i}{^{-1}}
\renewcommand{\th}{^{\mathrm{th}}}

\providecommand{\abs}[1]{\left\lvert#1\right\rvert}

\providecommand{\In}[1]{\left\langle#1\right\rangle}
\providecommand{\p}[1]{\left(#1\right)}

\providecommand{\f}[2]{\frac{#1}{#2}}
\DeclareMathOperator{\Ker}{Ker}			
\DeclareMathOperator{\Coker}{Coker}		
\DeclareMathOperator{\im}{Im}			
\DeclareMathOperator{\Spec}{Spec}		
\DeclareMathOperator{\Fitt}{Fitt}			
\DeclareMathOperator{\Frac}{Frac}		
\DeclareMathOperator{\id}{id}			
\DeclareMathOperator{\ord}{ord}	
\DeclareMathOperator{\GL}{GL}		
\DeclareMathOperator{\Aut}{Aut}		
\DeclareMathOperator{\Pic}{Pic}		
\DeclareMathOperator{\lcm}{lcm}		
\DeclareMathOperator{\Isom}{Isom}		
\newcommand{\ba}{\begin{aligned}}
\newcommand{\ea}{\end{aligned}}
\newcommand{\be}{\begin{equation}}
\newcommand{\ee}{\end{equation}}
\newcommand{\pf}{\begin{proof}}
\newcommand{\bpf}{\begin{proof}}
\newcommand{\epf}{\end{proof}}
\newcommand{\bthm}{\begin{thm}}
\newcommand{\ethm}{\end{thm}}
\newcommand{\bprop}{\begin{prop}}
\newcommand{\eprop}{\end{prop}}
\newcommand{\bcor}{\begin{cor}}
\newcommand{\ecor}{\end{cor}}
\newcommand{\brem}{\begin{rem}}
\newcommand{\erem}{\end{rem}}
\newcommand{\brems}{\begin{rems} \hfill \begin{enumerate}[label=\b{\thesubsubsection.},ref=\thesubsubsection]}
\newcommand{\bremstweak}{\begin{rems-tweak} \hfill \begin{enumerate}[label=\b{\thesubsection.},ref=\thesubsection]}
\newcommand{\remi}{\addtocounter{subsubsection}{1} \item}
\newcommand{\remitweak}{\addtocounter{subsection}{1} \item}
\newcommand{\erems}{\end{enumerate} \end{rems}}
\newcommand{\eremstweak}{\end{enumerate} \end{rems-tweak}}
\newcommand{\blem}{\begin{lemma}}
\newcommand{\elem}{\end{lemma}}
\newcommand{\bconj}{\begin{conj}}
\newcommand{\econj}{\end{conj}}
\newcommand{\bprob}{\begin{Problem}}
\newcommand{\eprob}{\end{Problem}}
\newcommand{\bq}{\begin{Q}}
\newcommand{\eq}{\end{Q}}
\newcommand{\benum}{\begin{enumerate}[label={{\upshape(\alph*)}}]}
\newcommand{\benuma}{\begin{enumerate}[label={{\upshape(\arabic*)}}]}
\newcommand{\benumr}{\begin{enumerate}[label={{\upshape(\roman*)}}]}
\newcommand{\eenum}{\end{enumerate}}
\newcommand{\bc}{}
\newcommand{\bd}{\begin{defn}}
\newcommand{\ed}{\end{defn}}
\newcommand{\beg}{\begin{eg}}
\newcommand{\eeg}{\end{eg}}
\newcommand{\bcl}{\begin{claim}}
\newcommand{\ecl}{\end{claim}}
\newcommand{\lab}{\label}

\newcommand{\q}{\quad}
\newcommand{\qq}{\quad\quad}
\newcommand{\qqq}{\quad\quad\quad}
\newcommand{\qqqq}{\quad\quad\quad\quad}
\newcommand{\qqqqq}{\quad\quad\quad\quad\quad}

\newcommand{\qqqqqqqq}{\quad\quad\quad\quad\quad\quad\quad\quad}

\newcommand{\Ell}{\cE\ell\ell}
\newcommand{\EEl}{\overline{\cE\ell\ell}}

\newcommand{\naive}{\mathrm{naive}}
\newcommand{\Tate}{\underline{\mathrm{Tate}}}

\newcommand{\tst}{\textstyle}

\newcommand*{\QED}{\hfill\ensuremath{\qed}}

\theoremstyle{plain}
\newtheorem{thm}[subsubsection]{Theorem}
\Crefname{thm}{Theorem}{Theorems}

\newtheorem{thm-tweak}[subsection]{Theorem}
\Crefname{thm-tweak}{Theorem}{Theorems}

\Crefname{rethm}{Theorem}{Theorem}
\newtheorem{prop}[subsubsection]{Proposition}
\Crefname{prop}{Proposition}{Propositions} 

\newtheorem{prop-tweak}[subsection]{Proposition}
\Crefname{prop-tweak}{Proposition}{Propositions} 

\newtheorem{Q}[subsubsection]{Question}
\Crefname{Q}{Question}{Questions}
\Crefname{eg}{Example}{Examples}
\newtheorem{Problem}[subsection]{Problem}
\Crefname{Problem}{Problem}{Problems}
\newtheorem{conj}[subsubsection]{Conjecture}
\Crefname{conj}{Conjecture}{Conjectures}
\newtheorem{cor}[subsubsection]{Corollary}
\Crefname{cor}{Corollary}{Corollaries}

\newtheorem{lemma}[subsubsection]{Lemma}

\newtheorem{lemma-tweak}[subsection]{Lemma}
\Crefname{lemma-tweak}{Lemma}{Lemmas}

\Crefname{subprop}{Proposition}{Propositions}
\Crefname{subcor}{Corollary}{Corollaries}
\Crefname{sublem}{Lemma}{Lemmas}

\theoremstyle{remark}
\newtheorem{claim}[equation]{Claim}
\Crefname{claim}{Claim}{Claims}

\Crefname{subrem}{Remark}{Remarks}

\theoremstyle{definition}
\newtheorem{defn}[subsubsection]{Definition}
\Crefname{defn}{Definition}{Definitions}
\newtheorem{conv}[subsubsection]{Convention}
\Crefname{conv}{Convention}{Conventions}
\newtheorem{eg}[subsubsection]{Example}
\newtheorem{rem}[subsubsection]{Remark}
\Crefname{rem}{Remark}{Remarks}

\newtheorem{rem-tweak}[subsection]{Remark}
\Crefname{rem-tweak}{Remark}{Remarks}

\newtheorem*{rems}{Remarks}

\newtheorem*{rems-tweak}{Remarks}

\newtheoremstyle{subsection-tweak}
   {11pt}
   {3pt}%
   {}
   {}%
   {\bfseries}
   {}%
   {.5em}
   {\thmnumber{\@{#1}{}\@{#2}.}%
    \thmnote{~{\bfseries#3.}}}

\Crefname{innercustomconj}{Conjecture}{Conjecture}

\theoremstyle{subsection-tweak}
\newtheorem{pp}[subsubsection]{}
\newcommand{\bpp}{\begin{pp}}
\newcommand{\epp}{\end{pp}}

\theoremstyle{subsection-tweak}
\newtheorem{pp-tweak}[subsection]{}

\numberwithin{equation}{subsubsection}

\makeatletter
\def\@tocline#1#2#3#4#5#6#7{
    \begingroup 
    \@ifempty{#4}{%
    }{%
    }%

    \parindent\z@ \leftskip#3\relax \advance\leftskip\@tempdima\relax
    #5\hskip-\@tempdima
      \ifcase #1
       \or\or \hskip 2em \or \hskip 1em \else \hskip 3em \fi%
      #6\nobreak\relax
    \dotfill\hbox to\@pnumwidth{\@tocpagenum{#7}}\par
    \nobreak
    \endgroup
  }
 \def\l@section{\@tocline{1}{0pt}{1pc}{}{}}

\renewcommand{\tocsection}[3]{%
  \indentlabel{\@ifnotempty{#2}{\makebox[1.3em][l]{%
    \ignorespaces#1 \bfseries{#2}.\hfill}}}\bfseries{#3}
    \vspace{1.5pt}}

\renewcommand{\tocsubsection}[3]{%
  \indentlabel{\@ifnotempty{#2}{\hspace*{-0.5em}\makebox[2.1em][l]{%
    \ignorespaces#1#2.\hfill}}}#3
    \vspace{1.5pt}}

\makeatother 

\usepackage[foot]{amsaddr}

\begin{document}

\author{K\k{e}stutis \v{C}esnavi\v{c}ius}
\title{A modular description of $\sX_0(n)$}
\date{\today}
\subjclass[2010]{Primary 11G18; Secondary 14D22, 14D23, 14G35.}
\keywords{Elliptic curve, generalized elliptic curve, level structure, modular curve, moduli stack}
\address{Mathematisches Institut, Universit\"{a}t Bonn, Endenicher Allee 60, D-53115, Bonn, Germany}
\email{kestutis@math.uni-bonn.de}

\begin{abstract} 
As we explain, when a positive integer $n$ is not squarefree, even over $\bC$ the moduli stack that parametrizes generalized elliptic curves equipped with an ample cyclic subgroup of order $n$ does not agree at the cusps with the $\Gamma_0(n)$-level modular stack $\sX_0(n)$ defined by Deligne and Rapoport via normalization. Following a suggestion of Deligne, we present a refined moduli stack of ample cyclic subgroups of order $n$ that does recover $\sX_0(n)$ over $\bZ$ for all $n$. The resulting modular description enables us to extend the regularity theorem of Katz and Mazur: $\sX_0(n)$ is also regular at the cusps. We also prove such regularity for $\sX_1(n)$ and several other modular stacks, some of which have been treated by Conrad by a different method. For the proofs we introduce a tower of compactifications $\EEl_m$ of the stack $\Ell$ that parametrizes elliptic curves---the ability to vary $m$ in the tower permits robust reductions of the analysis of Drinfeld level structures on generalized elliptic curves to elliptic curve cases via congruences. 
 \end{abstract}

\maketitle

\vspace{-1em}

\hypersetup{
    linktoc=page,     
}

\renewcommand*\contentsname{}
\q\\
\tableofcontents

\newpage


\section{Introduction}

\begin{pp-tweak}[Algebraic stacks that refine $X_0(n)$]
The study of the compactification $X_0(n)$ of the coarse moduli space of the algebraic stack $\sY_0(n)$ that parametrizes elliptic curves equipped with a cyclic subgroup of order $n$ is key for many arithmetic problems, so one seeks to understand the arithmetic properties of $X_0(n)$, especially over $\bZ$. For this, it is desirable to conceptualize the construction of $X_0(n)$ by realizing it as a coarse moduli space of an algebraic stack that compactifies $\sY_0(n)$.
 
The sought compactifying stack $\sX_0(n)$ was defined by Deligne and Rapoport in \cite{DR73}*{IV.3.3} via a normalization procedure. However, $\sX_0(n)$ lacks an \emph{a priori} moduli interpretation, so instead one often considers the stack $\sX_0(n)^\naive$ that parametrizes generalized elliptic curves whose smooth locus is equipped with a cyclic subgroup of order $n$ that is ample, i.e.,~meets every irreducible component of every geometric fiber. Even though $\sX_0(n)^\naive$ is algebraic, has $X_0(n)$ as its coarse moduli space, and agrees with $\sX_0(n)$ on the elliptic curve locus, it seems to have been overlooked that
\[
\text{If $n$ is not squarefree, then $\sX_0(n)$ and $\sX_0(n)^\naive$ are genuinely different, even over $\bC$.}
\]
\end{pp-tweak}

\begin{pp-tweak}[Pathologies of $\sX_0(p^2)^\naive$] \lab{non-rep}
To explain the difference, we set $n \ce p^2$ for some prime $p$, let $\sX(1)$ denote the stack that parametrizes those generalized elliptic curves whose geometric fibers are integral, and consider the structure morphism
\[
c\colon \sX_0(p^2)^\naive \ra \sX(1)
\]
which in terms of the moduli interpretation forgets the subgroup and contracts the generalized elliptic curve with respect to the identity section. We claim that the morphism $c$ is not representable.

To see this, let $E$ be the standard $p$-gon over $\bC$ and let $\zeta_{p^2} \in \bC^\times$ be a primitive root of unity of order $p^2$. Then $E^\sm = \bG_m \times \bZ/p\bZ$ and each of the $\mu_p$ worth of automorphisms of $E$ fixing $\bG_m \times \{ 0\}$ stabilizes the cyclic subgroup $\In{(\zeta_{p^2}, 1)}$ of order $p^2$. Each such automorphism contracts to the identity, so $c$ is not representable.

In contrast, the morphism 
\[
\sX_0(p^2) \ra \sX(1)
\]
is representable by construction, so the $\sX(1)$-stacks $\sX_0(p^2)^\naive$ and $\sX_0(p^2)$ are not isomorphic. The same $p$-gon example carried out over $\ov{\bF}_p$ shows that $\sX_0(p^2)^\naive$ is not even Deligne--Mumford (whereas $\sX_0(p^2)$ is), a pathology that has already been pointed out in \cite{Edi90}*{1.1.1.1} and \cite{Con07a}.
\end{pp-tweak}

\begin{pp-tweak}[A modular description of $\sX_0(n)$] \lab{mod-desc}
One of the main goals of this paper is to refine the definition of $\sX_0(n)^\naive$ to obtain a moduli interpretation of $\sX_0(n)$ even when $n$ is not squarefree. The elliptic curve locus needs no refinement, so the issue is to incorporate the cusps in a way that avoids the nonrepresentability of $c$ phenomenon. For this, we follow a suggestion of Deligne explained in \cite{Del15}. To present Deligne's idea, we assume that $n = p^2$ for a prime $p$ and work over $\bZ[\f{1}{p}]$. 

In vague terms, the idea is to subsume the automorphisms causing the nonrepresentability of $c$ into the moduli problem. To make this possible, the data being parametrized will involve algebraic stacks and not merely schemes. In precise terms, the moduli problem that in Chapter \ref{Gamma-0-case} will be proved to recover $\sX_0(p^2)_{\bZ[\f{1}{p}]}$ assigns to every $\bZ[\f{1}{p}]$-scheme $S$ the groupoid of tuples
\[
(E \ra S,\, G,\, S_{(1)},\, S_{(p)},\, S_{(p^2)},\, \cG_{(1)},\, \cG_{(p)},\, \cG_{(p^2)})
\]
consisting of
\begin{itemize}
\item
A generalized elliptic curve $E \ra S$;

\item
A cyclic subgroup $G \subset E_{S - S^{\infty}}$ of order $p^2$ over the elliptic curve locus $S - S^{\infty}$;

\item
Open subschemes $S_{(1)}$, $S_{(p)}$, and $S_{(p^2)}$ of $S$ that cover $S$, have $S - S^{\infty}$ as their pairwise intersections, and such that the degenerate geometric fibers of $E_{S_{(1)}}$ and $E_{S_{(p)}}$ are $1$-gons and those of $E_{S_{(p^2)}}$ are $p^2$-gons;

\item
Ample cyclic subgroups $\cG_{(1)} \subset E^\sm_{S_{(1)}}$ and $\cG_{(p^2)} \subset E^\sm_{S_{(p^2)}}$ of order $p^2$ that recover $G$ over~$S - S^\infty$;

\item
An ample cyclic subgroup $\cG_{(p)} \subset \cE_{(p)}^\sm$ of order $p^2$ of the universal generalized elliptic curve $\cE_{(p)}$ whose degenerate geometric fibers are $p$-gons and whose contraction is $E_{S_{(p)}}$, subject to the requirement that $\cG_{(p)}$ recovers $G$ over $S - S^\infty$ (over which $\cE_{(p)}$ is identified with $E$). 
\end{itemize}

In essence, the moduli problem parametrizes generalized elliptic curves equipped with an ample cyclic subgroup of order $p^2$ with the caveat that over the part of the degeneracy locus prone to the nonrepresentability of $c$  the subgroup has been upgraded to live inside a suitable universal ``decontraction'' $\cE_{(p)}$ (which is an algebraic stack and not a scheme). The role of the $S_{(p^i)}$ is to remember the subdivision of the degeneracy locus $S^\infty$---without $S_{(1)}$ and $S_{(p)}$ we cannot single out those $1$-gon degenerate geometric fibers of $E$ that were ``meant'' to be $p$-gons but had to be ``upgraded'' in order to avoid the nonrepresentability of $c$.
\end{pp-tweak}

\begin{pp-tweak}[Incorporating bad characteristics]
After the work of Drinfeld and of Katz and Mazur, the extension of the above modular description of $\sX_0(p^2)_{\bZ[\f{1}{p}]}$ to $\sX_0(p^2)$ is a matter of technique. However, new difficulties at the cusps in characteristic $p$ force us to impose an additional coherence requirement on $\cG_{(p)}$, a requirement that holds automatically away from $p$ and also on the elliptic curve locus (see \S\ref{coherence} and \Cref{compatible-gens}) and that seems well suited for the analysis of $\cG_{(p)}$ even over $\bZ[\f{1}{p}]$. With this proviso, we prove that for any $n$ the analogue of the moduli problem described in \S\ref{mod-desc} gives a moduli interpretation for $\sX_0(n)$. We then use this moduli interpretation to prove the following extension of a regularity theorem of Katz and Mazur:
\end{pp-tweak}

\begin{thm-tweak}[\Cref{X0n-main}~\ref{X0nM-a}] \lab{X0n-reg}
The Deligne--Mumford stack $\sX_0(n)$ is regular.
\end{thm-tweak}

In fact, $\sX_0(n)_{\bZ[\f{1}{n}]}$ is even $\bZ[\f{1}{n}]$-smooth by \cite{DR73}*{IV.6.7}, whereas the elliptic curve locus $\sY_0(n)$ is regular by \cite{KM85}*{5.1.1}, so \Cref{X0n-reg} was known away from the closed substack of the cusps that lies in characteristics dividing $n$. 

In the proof of \Cref{X0n-reg}, the eventual source of regularity is the combination of \cite{DR73}*{V.4.13} and \cite{KM85}*{5.1.1} that proves the regularity of another modular stack $\sX(n)$. The reduction to $\sX(n)$ rests on the moduli interpretation of $\sX_0(n)$ and on the regularity of $\sY_0(n)$. In particular, no stage of the argument requires any computations with universal deformation rings, other than what comes in from \cite{KM85}*{Ch.~5--6} through our reliance on the regularity of $\sY(n)$ and $\sY_0(n)$.

We use \Cref{X0n-reg} and the moduli interpretation of $\sX_0(n)$ to prove that the coarse moduli space $X_0(n)$ is regular in a neighborhood of the cusps (see \Cref{coarse-reg}). This regularity is not new (see the introduction of Chapter \ref{coarse-spaces}) but our proof seems more conceptual.

\begin{pp-tweak}[The compactifications $\EEl_m$]
We have been vague about the base of the universal ``decontraction'' $\cE_{(p)}$. For the construction of this base in general (beyond $n = p^2$), it is natural to fix an $m \in \bZ_{\ge 1}$ and to consider the $\bZ$-stack $\EEl_m$ that parametrizes those generalized elliptic curves whose degenerate geometric fibers are $m$-gons. We prove in \Cref{elln-props} that $\EEl_m$ is algebraic, as well as proper and smooth over $\bZ$, albeit is not Deligne--Mumford unless $m = 1$. Thus, each $\EEl_m$ compactifies the stack $\Ell$ that parametrizes elliptic curves, and $\EEl_1$ is the compactification that is sometimes called $\ov{\cM}_{1, 1}$.

As we describe in section \ref{reconstr}, the compactifications $\EEl_m$ form an infinite tower, with transition maps given by contractions of generalized elliptic curves. This tower is the backbone of our study of $\sX_0(n)$ and of several other ``classical'' modular curves. For these curves, the most important moduli-theoretic phenomenon that is not seen on the elliptic curve locus is the fact that ``forgetful'' contractions change generalized elliptic curves that underlie level structures. The ability to vary $m$ in the tower $\{\EEl_m\}_{m\mid m'}$ allows us to isolate the part of this phenomenon that has nothing to do with level structures. The remaining part that is specific to the level structure at hand may then be studied via ``congruences'' that reduce to the elliptic curve case.
\end{pp-tweak}

\begin{pp-tweak}[Other modular curves] \lab{Con-Con}
To illustrate the utility of $\EEl_m$, let us consider the stack $\sX(n)^\naive$ that parametrizes pairs consisting of a generalized elliptic curve $E \ra S$ with $n$-gon degenerate geometric fibers and a Drinfeld $(\bZ/n\bZ)^2$-structure on $E^\sm[n]$. (In the end, $\sX(n)^\naive$ agrees with $\sX(n)$ mentioned earlier and gives $\sX(n)$ a moduli interpretation.) Using the work of Katz and Mazur, we prove via ``mod $n$ congruences with elliptic curves'' that the forgetful map 
\[
\sX(n)^\naive \ra \EEl_n
\]
is representable and finite locally free of rank $\#\GL_2(\bZ/n\bZ)$. It follows that $\sX(n)^\naive$ is algebraic, proper and flat over $\bZ$, and even Cohen--Macaulay. Other proofs of these properties of $\sX(n)^\naive$ have been given by Conrad in \cite{Con07a}: the proof of the  algebraicity used Hilbert schemes via tricanonical embeddings, whereas the Cohen--Macaulay property required a detailed analysis of the universal deformation rings at the cusps (in addition to the work of Katz and Mazur on the elliptic curve~locus).

The relations with $\EEl_m$ together with the ``congruence method'' that crucially uses the work of Katz and Mazur allow us to reprove the main results of \cite{Con07a} in Chapter \ref{moduli-problems}. These include the moduli interpretations and the regularity of the modular stacks $\sX(n)$ and $\sX_1(n)$ (as well as some variants) and the construction of Hecke correspondences for $\sX_1(n)$. The latter takes advantage of the theory of isogenies of generalized elliptic curves developed in Chapter \ref{quotients}. Away from the level, the moduli interpretations and the regularity have been proved by Deligne and Rapoport in \cite{DR73}*{IV.3.5 and IV.4.14}; away from the cusps, they have been proved by Katz and Mazur in \cite{KM85}*{5.1.1}. Prior to the work of Conrad, \cite{Con07a}, the moduli interpretations and the regularity of $\sX(n)$ and $\sX_1(n)$ (among others) have been considered in an unfinished manuscript of Edixhoven \cite{Edi01}*{esp.~2.1.2}.
\end{pp-tweak}

\begin{pp-tweak}[Reliance on the literature]
For what concerns generalized elliptic curves and Drinfeld level structures on them, we wish to explicate the logical dependence of our work on the three main references that we use: \cite{DR73}, \cite{KM85}, and \cite{Con07a}. 
\begin{itemize}
\item 
We rely on \cite{DR73} almost in its entirety; the sections of op.~cit.~that are logically independent from the work of this paper are \cite{DR73}*{II.\S3, V.\S2--3, VI.\S2--6, and VII.\S3--4}. 

\item 
We make essential use of the results of \cite{KM85}*{Ch.~1--6} and extend some of them to generalized elliptic curves (see, in particular, section \ref{KM-fest}), but have no need for the results of \cite{KM85}*{Ch.~7--14} (other than for comparison in \Cref{KM-comp} and Remarks~\ref{KM-version} and \ref{rem-more-reg}). 

\item
We use some auxiliary general results from the introductory sections 2.1 and 2.2 of \cite{Con07a} but the rest of op.~cit.~is logically independent from our work (as mentioned in \S\ref{Con-Con}, we give different proofs to the main results of \cite{Con07a}).
\end{itemize}
\end{pp-tweak}

\begin{pp-tweak}[Notation and conventions] \lab{conv}
We let $\Ell$ denote the $\bZ$-stack that, for variable schemes $S$, parametrizes elliptic curves $E \ra S$. More precisely, for a scheme $S$, the objects (resp.,~the morphisms) of the groupoid $\Ell(S)$ are the elliptic curves $E \ra S$ (resp.,~the isomorphisms between elliptic curves over $S$) and, for a scheme morphism $S' \ra S$, the induced functor $\Ell(S) \ra \Ell(S')$ is $E \mapsto E \times_S S'$. We use the analogous meaning of `parametrizes' when defining other stacks.  Other than in the introduction, we use the notation $\sX_{\Gamma_0(n)}$ (resp.,~$\sX_{\Gamma_1(n)}$, etc.) introduced in \S\ref{genl-level} for stacky modular curves defined via normalization and the notation $\sX_0(n)$ (resp.,~$\sX_1(n)$, etc.) for stacks defined in terms of a moduli problem; once we prove that $\sX_{\Gamma_0(n)} = \sX_0(n)$ (and similarly in the other cases), we use the two notations interchangeably.

We use the definition of an fpqc cover for which all Zariski covers are fpqc; explicitly, $S' \ra S$ is an fpqc cover if it is flat and every affine open $U \subset S$ is the union of images of finitely many affine opens of $S'$. An $S$-scheme $S'$ is an fppf cover (or simply fppf) if $S' \ra S$ is faithfully flat and locally of finite presentation. For a scheme $S$, we let $S^{\mathrm{red}}$ denote its associated reduced scheme. For an $S$-group algebraic space $G$, we let $G^0$ denote the subsheaf of sections that fiberwise factor through the identity component. We let $\sX^\sm$ and $\Delta_{\sX/S}$ denote the smooth locus and the diagonal of a morphism $\sX \ra S$. For a field $k$, we let $\ov{k}$ denote a choice of its algebraic closure. A geometric point is the spectrum of an algebraically closed field. For an $n \in \bZ_{\ge 1}$, we set $\phi(n)\ce \#(\bZ/n\bZ)^\times$.

For what concerns algebraic stack and algebraic space conventions, we follow \cite{SP}, except that `representable' stands for `representable by algebraic spaces.' In particular, quasi-compactness or separatedness of the diagonal are not part of the definition, but in practice end up being present (along with even stronger properties). An algebraic stack is Deligne--Mumford if its diagonal is  unramified---for the equivalence with the \'{e}tale atlas definition in the presence of quasi-compactness and separatedness of the diagonal, see \cite{LMB00}*{8.1}. The relative dimension (at a point) of a smooth morphism of algebraic stacks is the difference of the relative dimensions (at a lift of the point) of the morphisms from a smooth atlas of the source, cf.~\cite{LMB00}*{bottom of p.~98}.
\end{pp-tweak}

\subsection*{Acknowledgements}
 I thank Pierre Deligne for correspondence about the moduli interpretation in the $\Gamma_0(n)$ case and for permitting me to make his letter \cite{Del15} available. The modular description of $\sX_0(n)$ presented in Chapter \ref{Gamma-0-case} is inspired by the ideas explained there. I thank the referee for a very careful reading of the manuscript and for numerous helpful suggestions. I thank the MathOverflow community---the reading of several anonymous discussions has been useful while working on some aspects of this paper. I thank Rebecca Bellovin, George Boxer, Brian Conrad, Bas Edixhoven, Benedict Gross, Dino Lorenzini, Martin Olsson, Ken Ribet, and Sug Woo Shin for helpful conversations or correspondence.  I thank the Miller Institute for Basic Research in Science at the University of California Berkeley for its support.


\section{Isogenies of generalized elliptic curves} \lab{quotients}

The main goal of this chapter is to expose a robust theory of isogenies of generalized elliptic curves. This theory is the subject of section \ref{quotient-section} and will be useful on several occasions, particularly, for algebraizing homomorphisms of formal generalized elliptic curves in section \ref{GAGAmama} and for constructing Hecke correspondences for $\sX_1(n)$ in section \ref{Hecke}. In order to prepare for the study of isogenies, in section \ref{homomomo} we review several basic concepts, such as that of a homomorphism of generalized elliptic curves, and record some general results that will be useful throughout the paper.


\centering \subsection{Homomorphisms between generalized elliptic curves\nopunct} \lab{homomomo} \hfill
\justify

In this section, we review basic definitions and properties of generalized elliptic curves, building up to the notion of a homomorphism, which will be studied in section \ref{quotient-section}. We assume that the reviewed concepts are familiar, so we concentrate on those aspects that will be used later. We begin with the notion of an $n$-gon, which is needed in order to define generalized elliptic curves.  Informally, an $n$-gon is the curve obtained by glueing $n$-copies of $\bP^1$ in a cyclic manner: the point $0$ of the $i\th$ copy gets identified with the point $\infty$ of the $(i+1)^\st$ copy.

\bd \lab{n-gon-def}
For an $n \in \bZ_{\ge 1}$ and an scheme $S$, the \emph{standard $n$-gon over $S$} is the coequalizer  of 
\[
\tst \xymatrix{ \bigsqcup_{\bZ/n\bZ} S \ \ar@{_(->}@<-.5ex>[r] \ar@{^(->}@<.5ex>[r] & \ \bigsqcup_{\bZ/n\bZ} \bP^1_S, }
\]
where the top (resp.,~the bottom) closed immersion includes the $i\th$ copy of $S$ as the $0$ (resp.,~the $\infty$) section of the $i\th$ (resp.,~$(i + 1)^\st$) copy of $\bP^1_S$. A \emph{N\'{e}ron $n$-gon over $S$} (or an \emph{$n$-gon over $S$}) is an $S$-scheme isomorphic to the standard $n$-gon over $S$. (We often omit `over $S$' if the base is implicit.)
\ed

\brem \lab{n-gon-concrete}
Even though colimits usually do not exist in the category of schemes, the ones used in \Cref{n-gon-def} do exist and their formation commutes with base change in $S$. To see this, one checks directly (or with the help of \cite{Fer03}*{4.3}) that for $n \ge 2$ the sought coequalizer is the base change to $S$ of the glueing of 
\[
\tst \bigsqcup_{i \in \bZ/n\bZ} \Spec \bZ[X_i, Y_i]/(X_iY_i)
\]
obtained by identifying the opens 
\[
\tst \Spec \bZ[Y_i, \f{1}{Y_i}] \qq \text{and} \qq \Spec \bZ[X_{i + 1}, \f{1}{X_{i + 1}}]
\]
via $Y_i = \f{1}{X_{i + 1}}$ for every $i \in \bZ/n\bZ$, and one treats the $n = 1$ case by realizing the standard $1$-gon as the $\bZ/n\bZ$-quotient of the standard $n$-gon, cf.~\cite{Con07a}*{top of p.~215}.
\erem

We recall the definition of a generalized elliptic curve, which is a central notion for this paper.

\bd \lab{gen-ell-def}
A \emph{generalized elliptic curve} over a scheme $S$ is the data of
\begin{itemize}
\item
A proper, flat, finitely presented morphism $E \ra S$ each of whose geometric fibers is either a smooth connected curve of genus $1$ or a N\'{e}ron $n$-gon for some $n \ge 1$, and 

\item
An $S$-morphism $E^\sm \times_S E \xra{+} E$ that restricts to a commutative $S$-group scheme structure on $E^\sm$ for which $+$ becomes an $S$-group action,
\end{itemize}
such that via pullback of line bundles the action $+$ induces the trivial action of $E^\sm$ on $\Pic^0_{E/S}$.
\ed

\brems
\remi \lab{comp-gen-def}
Our definition of a generalized elliptic curve is equivalent to the one given in \cite{DR73}*{II.1.12}: the difference is that we have imposed the requirement that $E^\sm$ acts trivially on $\Pic^0_{E/S}$ at the outset. \emph{Loc.~cit.}~replaces this with the \emph{a priori} milder requirement that on degenerate geometric fibers every translation by a smooth point induces a rotation on the underlying $n$-gon, which ends up being equivalent due to \cite{DR73}*{II.1.7 (ii) and II.1.13}. 

The requirement about the triviality of the induced action on $\Pic^0_{C/S}$ holds automatically on a large part of $E^{\sm}$, namely, it always holds on the relative identity component $(E^{\sm})^0$---to see this, we apply \cite{DR73}*{II.1.14}\footnote{We could also apply \cite{Con07a}*{2.2.1} to avoid using the representability of $\Pic^0_{E/S}$ by a scheme. On the other hand, such representability may be proved as follows: by \cite{Art69a}*{7.3}, the functor $\Pic^0_{E/S}$ is an algebraic space, so \cite{DR73}*{II.2.6~(i)} proves that the map 
\[
(E^\sm)^0 \ra \Pic^0_{E/S} \qq \text{defined by} \qq t \mapsto \sO_E(t) \tensor \sO_E(e)\i
\]
is an open immersion (where $e \in E(S)$ denotes the identity section), and the representability of $\Pic^0_{E/S}$ by a scheme follows from \cite{BLR90}*{6.6/2~(b)} applied to $\Pic^0_{E/S}$ acting on itself by translation (see also \Cref{gen-ell-limit}).} to $\Pic^0_{E/S} \times_S E^\sm$ to get the openness of the locus of $E^\sm$ where the induced action on $\Pic^0_{E/S}$ is trivial, note that this locus is closed under the group law of $E^\sm$, and conclude by noting that it contains the zero section. In particular, every elliptic curve is a generalized elliptic curve, and a generalized elliptic curve $E \ra S$ is an elliptic curve over the open of $S$ over which $E$ is smooth.

\remi \lab{n-gon-gen-ell}
The standard $n$-gon is canonically a generalized elliptic curve: due to its description recalled in \Cref{n-gon-concrete}, its smooth locus is $\bG_m \times \bZ/n\bZ$ and the translation action of this group scheme on itself extends to an action on the $n$-gon. By the previous remark, the triviality of the induced action on $\Pic^0$ may be checked on the geometric fibers using \cite{DR73}*{II.1.7~(ii)}. For later use, we now describe the automorphism functor of this generalized elliptic curve.
\erems

\blem \lab{auto-n-gon}
For a fixed $n \in \bZ_{\ge 1}$, let $E \ra \Spec \bZ$ be the standard $n$-gon generalized elliptic curve. There is the following identification of the automorphism functor of $E$:
\[
\Aut(E) \cong \mu_n \times \bZ/2\bZ,
\]
where the generator of $\bZ/2\bZ$ acts as inversion on $E^\sm$ and, for a scheme $S$ and an index $i \in \bZ/n\bZ$, a section $\zeta \in \mu_n(S)$ acts on the $i\th$ component of $E^\sm_S \cong (\bG_m)_S \times \bZ/n\bZ$ as scaling by $\zeta^i$.
\elem

\bpf
By \cite{DR73}*{II.1.10}, we have 
\[
\Aut(E) \cong \mu_n \rtimes \bZ/2\bZ
\]
with $\mu_n$ and $\bZ/2\bZ$ acting as described above, so we need to argue that $\bZ/2\bZ$ is central in $\Aut(E)$. For this, due to the $\bZ$-universal schematic density of $E^\sm$ in $E$ supplied by \cite{EGAIV3}*{11.10.10}, it suffices to note that every generalized elliptic curve automorphism of a base change of $E$ must commute with inversion on $E^\sm$. 
\epf

We turn to the closed subschemes $E^{\sing} \subset E$ and $S^{\infty, \pi} \subset S$ that measure the degeneration of $E$.

\bd \lab{def-sing}
The \emph{subscheme of nonsmoothness} of a generalized elliptic curve $E \xra{\pi} S$ is the closed subscheme $E^{\sing} \subset E$ defined by the first Fitting ideal sheaf $\Fitt_1(\Omega^1_{E/S}) \subset \sO_E$. The \emph{degeneracy locus} of $E \xra{\pi} S$ is the schematic image $S^{\infty, \pi} \subset S$ of  $E^{\sing}$.
\ed

\brems
\remi \lab{sing-bc}
The closed subscheme $E^{\sing}$ is supported at those points of $E$ at which $\pi$ is not smooth and its formation commutes with arbitrary base change in $S$, see \cite{SGA7I}*{VI, 5.3 and 5.4}. Even though the formation of schematic images often does not commute with nonflat base change, the formation of $S^{\infty, \pi}$ does commute with arbitrary base change, see \cite{Con07a}*{2.1.12}. 

\remi \lab{gen-loc-proj}
By \cite{DR73}*{II.1.15}, we have 
\[
\qq\tst  S^{\infty, \pi} = \bigsqcup_{n \ge 1} S^{\infty, \pi, n}
\] 
for closed subschemes $S^{\infty, \pi, n} \subset S$ such that only finitely many of the $S^{\infty, \pi, n}$ meet a given affine open of $S$ and such that $E_{S^{\infty, \pi, n}}$ is fppf locally on $S^{\infty, \pi, n}$ isomorphic to the standard $n$-gon (which was discussed in Remark \ref{n-gon-gen-ell}). In particular, every generalized elliptic curve $E \xra{\pi} S$ is, Zariski locally on $S$, projective because, by \cite{DR73}*{II.1.20} and \cite{KM85}*{1.2.3}, over the open 
\[
\qq \tst S - \bigsqcup_{n \neq n'} S^{\infty, \pi, n}
\]
the $n'$-torsion subscheme $E^\sm[n'] \subset E$ is a $\pi$-ample relative effective Cartier divisor.
\erems

We record a basic relationship between $E^\sing$ and its schematic image $S^{\infty, \pi}$ in the following lemma:

\blem \lab{sing-map}
For a generalized elliptic curve $E \ra S$, the map 
\[
E^\sing \ra S^{\infty, \pi}
\]
is finite \'{e}tale; it has degree $n$ over $S^{\infty, \pi, n}$.
\elem

\bpf
The map in question exists by the definition of $S^{\infty, \pi}$ and its formation commutes with base change in $S$ by Remark \ref{sing-bc}. We may therefore assume that $S = S^{\infty, \pi, n}$ and that $E$ is the standard $n$-gon. But in this case $E^\sing$ is a disjoint union of $n$ copies of $S$ and there is nothing to prove.
\epf

Degenerate generalized curves possess canonical finite subgroups of multiplicative type and their torsion subgroups are amenable to scrutiny. We make this precise in the following lemma:

\blem \lab{can-submult}
For every generalized elliptic curve $E \xra{\pi} S$ with $S^{\mathrm{red}} = (S^{\infty, \pi})^{\mathrm{red}}$ and every $d \in \bZ_{\ge 1}$, the $d$-torsion $(E^\sm)^0[d]$ is a finite locally free $S$-group scheme  of order $d$ that is \'{e}tale locally on $S$ isomorphic to $\mu_d$. The $S$-group scheme 
\[
E^\sm[d]/(E^\sm)^0[d]
\]
is \'{e}tale and if $m \in \bZ_{\ge 1}$ divides both $d$ and the number of irreducible components of each geometric fiber of $E$, then $(E^\sm[d]/(E^\sm)^0[d])[m]$ is \'{e}tale locally on $S$ isomorphic to $\bZ/m\bZ$.
\elem

\bpf
Due to the fibral criterion for flatness \cite{EGAIV3}*{11.3.11}, the quasi-finite, finitely presented, separated $S$-groups $(E^\sm)^0[d]$ and $E^\sm[d]$ are flat. The fibers of $(E^\sm)^0[d] \ra S$ have degree $d$, so, due to \cite{DR73}*{II.1.19}, the $S$-group $(E^\sm)^0[d]$ is finite locally free of rank $d$. Due to \cite{Con14}*{B.4.1 and B.3.4}, the claim about the \'{e}tale local structure of $(E^\sm)^0[d]$ reduces to case of geometric fibers.

Thanks to the settled claims about $(E^\sm)^0[d]$, \cite{EGAIV3}*{8.11.2} and \cite{SGA3Inew}*{V, 4.1} imply that $E^\sm[d]/(E^\sm)^0[d]$ is a separated, quasi-finite, finitely presented, flat $S$-scheme. By inspecting geometric fibers we see that $E^\sm[d]/(E^\sm)^0[d]$ is \'{e}tale. The \'{e}tale local structure of 
\[
(E^\sm[d]/(E^\sm)^0[d])[m]
\]
may be seen over the strict Henselizations of $S$, and hence even on geometric fibers.
\epf

The focus of Chapter \ref{quotients} is generalized elliptic curve homomorphisms. We recall their definition.

\bd\lab{homo-def}
A \emph{homomorphism} between generalized elliptic curves $E \ra S$ and $E' \ra S$ is an $S$-morphism 
\[
f\colon E \ra E' \qq \text{with} \qq f(E^\sm) \subset E'^{\sm}
\]
that intertwines the group laws of $E^\sm$ and $E'^{\sm}$. Its \emph{kernel} is the $S$-subscheme $\Ker f \ce E \times_{f,\, E',\, e'} S$ of $E$, where $\times_{f,\, E',\, e'}$ denotes the base change along $f$ of the identity section $e'\colon S \ra E'$.
\ed

\brems
\remi \lab{homo-equiv}
Due to the $S$-universal schematic density of $E^\sm$ in $E$ supplied by \cite{EGAIV3}*{11.10.10} and the separatedness of $E' \ra S$, a homomorphism $f$ necessarily also intertwines the group actions $E^\sm \times E \ra E$ and $E'^{\sm} \times E' \ra E'$.

\remi \lab{homo-surj}
If a homomorphism $f$ is surjective, then $f|_{E^\sm}$ is flat and $\Ker f \subset E^\sm$, as may be checked on geometric fibers using the fibral criterion for flatness \cite{EGAIV3}*{11.3.11}. In this case, $\Ker f$ is a finite locally free $S$-subgroup scheme of $E^\sm$.
\erems

\beg \lab{eg-homo}
The constant morphism that factors through $e'$ is a homomorphism, the ``zero homomorphism.'' Any elliptic curve isogeny is also a homomorphism. For a $d \in \bZ_{\ge 1}$, the map 
\[
\bP^1_S \ra \bP^1_S \qq \text{given on homogeneous coordinates by} \qq [x : y] \mapsto [x^d : y^d]
\]
respects $0$ and $\infty$, so it induces an $S$-morphism from the standard $1$-gon over $S$ to itself. This morphism restricts to the $d\th$ power map on the $(\bG_m)_S$ of the smooth locus of the $1$-gon, so it is a homomorphism with kernel $(\mu_d)_S$.
\eeg

\brem \lab{gen-ell-limit}
Generalized elliptic curves are susceptible to limit arguments that reduce to a Noetherian base. More precisely, by \cite{EGAIV2}*{8.8.2~(ii), 8.10.5~(xii), 11.2.6~(ii)}, Zariski locally on $S$, the underlying relative curve $E \ra S$ is the base change of a proper and flat relative curve $E_0 \ra S_0$ for which $S_0$ is of finite type over $\bZ$. Thus, since the formation of $E_0^\sm$ commutes with base change, $E^\sm$ is necessarily of finite presentation. Moreover, by \cite{EGAIV2}*{8.8.2~(i)}, after enlarging $S_0$, the commutative $S$-group action 
\[
E^\sm \times_S E \xra{+} E \qq \text{descends to a commutative $S_0$-group action} \qq E^\sm_0 \times_{S_0} E_0 \xra{+} E_0.
\]
The degenerate geometric fibers of $E_0 \ra S_0$ are N\'{e}ron $n$-gons: indeed, \cite{DR73}*{II.1.3} applies because the condition of having only ordinary double points as singularities is equivalent to the unramifiedness of $E^\sing_0$, whose formation commutes with base change (see~Remark~\ref{sing-bc}), whereas the triviality of the relative dualizing sheaf may be descended from an overfield using specialization techniques. Using Remark~\ref{comp-gen-def} to infer the triviality of the induced action of $E_0^\sm$ on $\Pic^0_{E_0/S_0}$, we conclude that $E_0 \ra S_0$ is a generalized elliptic curve that descends $E \ra S$ to a Noetherian base. Similarly, Zariski locally on $S$, elliptic curve homomorphisms are defined over a base that is of finite type over $\bZ$.

By the limit arguments above, the open immersion $S - S^{\infty, \pi} \hra S$ is always quasi-compact. 
\erem


\centering \subsection{Quotients of generalized elliptic curves by finite locally free subgroups \nopunct} \lab{quotient-section} \hfill

\justify

Even though homomorphisms between generalized elliptic curves are useful in practice, their structural properties are not immediately apparent. Moreover, guided by the theory of isogenies of elliptic curves, one suspects that for any finite locally free $S$-subgroup scheme $G \subset E^\sm$ with $E \ra S$ a generalized elliptic curve, there should be an essentially unique homomorphism $E \ra E'$ with kernel $G$. If $G$ intersects the identity components of the degenerate geometric fibers of $E \ra S$ trivially, then the translation action of $G$ on $E$ is free, the fppf sheaf quotient $E/G$ is a generalized elliptic curve, and 
\[
E \ra E/G
\] 
is the sought ``isogeny.'' This special case is already useful---for instance, such isogenies are discussed in \cite{Con07a}*{2.1.6} and exploited in several key proofs of \emph{op.~cit.}

The goal of this section is to explain how to make sense of isogenies of generalized elliptic curves in general. \Cref{naive-quot} and its proof explain how to build the desired ``quotient by $G$'' homomorphism $E \ra E/G$, and we arrive at the concept of an isogeny in \Cref{def-isog}. With \Cref{naive-quot} in hand, structural properties of arbitrary homomorphisms are susceptible to scrutiny and are detailed in \Cref{Zar-loc-isog,isog-sing}.

We begin with an example that illustrates what $E/G$ should be in a certain degenerate situation.

\beg \lab{n-gon-quot}
Let $E$ be the standard $n$-gon over $\bZ$, and consider the subgroup $\mu_d \subset (E^\sm)^0$ for some $d \in \bZ_{\ge 1}$. We would like to build a generalized elliptic curve homomorphism 
\[
f_d\colon E \ra E' \qq \text{with kernel $\mu_d$.}
\]
By Remark \ref{homo-equiv}, any such $f_d$ is $\mu_d$-equivariant, so it factors through the categorical quotient $E/\mu_d$, which exists because $E$ is projective and $\mu_d$ is finite. We claim that 
\[
E \ra E/\mu_d
\]
is already the desired $f_d\colon E \ra E'$.

This claim follows from the description of $E$ recalled in \Cref{n-gon-concrete}. More precisely, if $n \ge 2$, then on $\Spec\p{ \f{\bZ[X_i, Y_i]}{(X_iY_i)}}$ the action of $\mu_d = \Spec \p{\f{\bZ[T]}{(T^d - 1)}}$ is determined by 
\[
X_i \mapsto X_i \tensor T \qq \text{and} \qq Y_i \mapsto Y_i \tensor T,
\]
so the ring of invariants is the $\bZ$-subalgebra of $\f{\bZ[X_i, Y_i]}{(X_iY_i)}$ generated by $X_i^d$ and $Y_i^d$, and hence $E/\mu_d$ is the standard $n$-gon with the quotient map $E \ra E/\mu_d$ induced by the $d\th$ power map on each $\bP^1_{\bZ}$. The same description holds if $n = 1$, as the same computation performed $\bZ/m\bZ$-equivariantly on the $m$-gon cover for some $m \ge 2$ proves. Thus, the map $E \ra E/\mu_d$ is a homomorphism whose kernel is $\mu_d$, and it is initial among such homomorphisms, so it is the desired $f_d$.
\eeg

\brems
\remi \lab{quot-bc}
\Cref{n-gon-quot} may be carried out over any base scheme $S$, which shows that the formation of $f_d$ commutes with arbitrary base change. In particular, the formation of the categorical quotient $E/\mu_d$ commutes with arbitrary (possibly nonflat) base change.

\remi \lab{isog-non-flat}
For $d > 1$, the ``isogeny'' $E \ra E/\mu_d$ constructed in \Cref{n-gon-quot} is not flat at the singular points, as the formal criterion for flatness \cite{BouAC}*{III, \S5, n$^{\circ}$~2, Thm.~1} reveals. In contrast, every isogeny between elliptic curves is flat.
\erems

\Cref{n-gon-quot} suggests that over an arbitrary base $S$, the desired quotient of a generalized elliptic curve $E \ra S$ by a finite locally free $S$-subgroup $G \subset E^\sm$ may simply be the categorical quotient $E/G$. In \Cref{naive-quot} we prove that this indeed the case. The main issue that needs to be addressed is that the formation of categorical quotients does not in general commute with nonflat base change (as in the special case of forming the ring of invariants under the action of a finite group). Such phenomena do not occur for generalized elliptic curves because the analysis of $E/G$ may be reduced to the cases when $G$ is either diagonalizable or acts freely on $E$.

\bthm \lab{naive-quot}
Let $S$ be a scheme, $E \xra{\pi} S$ a generalized elliptic curve, and $G \subset E^\sm$ an $S$-subgroup scheme that is finite locally free over $S$. There is an $S$-scheme morphism 
\[
q\colon E \ra E/G
\]
that is initial among $G$-equivariant $S$-morphisms from $E$ to an $S$-scheme equipped with the trivial $G$-action {\upshape(}$E$ is equipped with the translation action of $G${\upshape)}. Moreover, $q$ has the following properties.
\benumr
\item \lab{NQ-1}
The formation of $q$ commutes with arbitrary base change in $S$, and $E/G$ is $S$-flat.

\item \lab{NQ-2}
The map $q\colon E \ra E/G$ is surjective, finite, and universally open.

\item \lab{NQ-4}
There is a unique structure of a generalized elliptic curve on 
\[
\qq E/G \ra S
\]
for which $q$ is a homomorphism. For this structure, $q$ induces an $S$-group isomorphism 
\[
\qq E^\sm/G \cong (E/G)^\sm,
\]
where $E^\sm/G$ is the fppf sheaf quotient; in particular, $E^\sm \xra{q} (E/G)^\sm$ is finite locally free.

\item \lab{NQ-5}
If $E$ is an elliptic curve, then $q\colon E \ra E/G$ is an isogeny with kernel $G$.
\eenum
\ethm

\bpf
Zariski locally on $S$ the map $\pi$ is projective (see Remark \ref{gen-loc-proj}), so every finite set of points of any $\pi$-fiber is contained in an affine open of $E$ (see \cite{EGAII}*{4.5.4}). Therefore, by \cite{SGA3Inew}*{V,~4.1~(i)} and its proof, $E$ is covered by $G$-invariant affine opens and the initial $q$ is nothing else than the categorical quotient that is glued together from the rings of invariants of such $G$-invariant affines; moreover, this $q$ is automatically a quotient map on the underlying topological spaces.

Since $G$ acts freely on $E^\sm$, by \cite{SGA3Inew}*{V,~4.1~(iv)}, the open $S$-subscheme 
\[
E^\sm/G \subset E/G
\]
that results from the $G$-invariance of $E^\sm$ is identified with the fppf sheaf quotient of $E^\sm$ by $G$, the map $E^\sm \xra{q} E^\sm/G$ is finite locally free, and the formation of $E^\sm/G$ commutes with base change. 

\benumr
\item
The formation of $E/G$ commutes with flat base change, so we may first assume that $S$ is affine and then use Remark \ref{gen-ell-limit} to assume that $S = \Spec R$ for some Noetherian $R$. Moreover, by the previous paragraph, the claim is clear on the elliptic curve locus, so we may replace $R$ by its completion along the ideal $I \subset R$ that cuts out the degeneracy locus $S^{\infty, \pi} \subset S$ to assume that $R$ is $I$-adically complete and separated.

For such $R$, the intersections 
\[
\qq G_{R/I^j} \cap (E^\sm_{R/I^j})^0 \qq \text{for} \qq j \ge 1
\]
are finite locally free $R/I^j$-subgroup schemes of $G$. By Grothendieck's existence theorem \cite{Ill05}*{8.4.5, 8.4.7}, these subgroups algebraize to a finite locally free $R$-subgroup 
\[
\qq H \subset G \qq \text{with} \qq H \subset (E^\sm)^0.
\]
The $R/I$-fibers of $H$ are of multiplicative type, so $H$ itself is of multiplicative type. At the cost of replacing $R$ by a finite locally free cover we may assume that $H$ is diagonalizable.

By \cite{SGA3Inew}*{I, 4.7.3}, any $R$-module $M$ equipped with an action of a diagonalizable $H$ is a direct sum of $\chi$-isotypic submodules for characters $\chi$ of $H$, so the submodule $M^H$ of $H$-invariants is of formation compatible with arbitrary base change and is $R$-flat if $M$ is. In particular, the categorical quotient $E/H$ is $R$-flat and of formation compatible with base change. As may be checked on geometric $R$-fibers, $G/H$ acts freely on $E/H$, so the further quotient $E/G = (E/H)/(G/H)$ is also $R$-flat and of formation compatible with base change.

\item
The surjectivity of $q$ follows from the first paragraph of the proof. By \cite{SGA3Inew}*{V,~4.1~(ii)}, the morphism $q$ is integral, and hence even finite because it inherits the property of being of finite type from $E \ra S$. In particular, $q$ is universally closed,  so it is also universally open by \cite{Ryd13}*{2.4} (which applies due to  the bottom of p.~636 there and \cite{SGA3Inew}*{V,~4.1~(iii)}).

\item
We begin by arguing that $E/G$ possesses the $S$-scheme properties required in \Cref{gen-ell-def}.

Due to \cite{AM69}*{7.8}, the morphism $E/G \ra S$ inherits finite presentation from $E \ra S$ thanks to the finiteness of $E \ra E/G$ (and an initial reduction to Noetherian $S$ based on \ref{NQ-1}). By \ref{NQ-2}, 
\[
\qq E \ra E/G, \qq \text{and hence also} \qq  E \times_S E \ra E/G \times_S E/G,
\]
is a finite surjection, so the image of $\Delta_{E/S}(E)$ in $E/G \times_S E/G$, i.e., $\Delta_{(E/G)/S}(E/G)$, is closed. In other words,  
the finite type morphism $E/G \ra S$ inherits separatedness from $E \ra S$, so it also inherits properness by \cite{EGAII}*{5.4.3 (ii)}. Finally, $E/G \ra S$ is flat by \ref{NQ-1}. For the fibral properties, due to \ref{NQ-1}, we may assume that $S$ is a geometric point.

If $S$ is a geometric point and $E$ is an elliptic curve, then $E/G$ is its isogenous quotient. If $S$ is a geometric point and $E$ is the standard $N$-gon, then we set 
\[
\qq H \ce G \cap (E^\sm)^0, \qq \text{so} \qq H \cong \mu_d \q \text{for some} \q d \ge 1.
\]
By \Cref{n-gon-quot}, $E \ra E/H$ is a ``self-isogeny'' of the standard $N$-gon, and, by construction, $G/H$ acts freely on $E/H$. Therefore, $E/G$, which is identified with $(E/H)/(G/H)$, is the standard $n$-gon with $n = \f{N}{\#(G/H)}$. This analysis also shows that $q(E^\sm) = (E/G)^\sm$.

Due to the paragraph preceding the proof of \ref{NQ-1}, all that remains to be shown is that the $S$-group scheme structure of $(E/G)^\sm \cong E^\sm/G$ extends to a unique action of $(E/G)^\sm$ on $E/G$; indeed, the induced action on $\Pic^0_{(E/G)/S}$ will automatically be trivial due to the fibral analysis of the previous paragraph and Remark \ref{comp-gen-def}. The uniqueness follows from the separatedness of $E/G$ and the universal schematic density of $(E/G)^\sm$ in $E/G$ supplied by \cite{EGAIV3}*{11.10.10}. For the same reason, for the existence we only need to produce a morphism 
\[
\qq (E/G)^\sm \times_S E/G \ra E/G
\]
that extends the group law of $(E/G)^\sm$---the relevant diagrams that encode the property of being a group scheme will automatically commute. To build this morphism from the one for $E$, it suffices to prove that
\[
\qq E^\sm/G \times_S E/G \cong (E^\sm \times_S E)/(G\times_S G)
\]
where the quotients are categorical. For this isomorphism, it suffices to form the quotient on the right in stages and to note that the formation of $E^\sm/G$ commutes with base change along $E \ra S$ whereas the formation of $E/G$ commutes with base change along $E^\sm/G \ra S$.

\item
By \ref{NQ-4}, $q\colon E \ra E/G$ is a finite locally free homomorphism between elliptic curves over $S$ and its kernel is $G$, i.e.,~$q$ is an isogeny with kernel $G$.
\qedhere
\eenum
\epf

\brem
The categorical quotient $E/G$ may also be analyzed with the tame stack formalism  of Abramovich--Olsson--Vistoli, \cite{AOV08}. For this, the key point is that the quotient stack $[E/G]$ is tame by \cite{AOV08}*{Thm.~3.2} because the automorphism functors of its geometric points are of multiplicative type. Then, since $E/G$ is the coarse moduli space of $[E/G]$ (see~\cite{Con05}*{Thm.~3.1}), $E/G$ is $S$-flat and of formation compatible with arbitrary base change by \cite{AOV08}*{Cor.~3.3}.
\erem

\bpp[The quotient notation] \lab{quot-not}
In the sequel, whenever $E \ra S$ is a generalized elliptic curve and $G \subset E^\sm$ is a finite locally free $S$-subgroup, we write $E/G$ for the generalized elliptic curve constructed in \Cref{naive-quot}. In the following corollary, we record some further properties of this quotient construction that follow from \Cref{naive-quot} and its proof.
\epp

\bcor \lab{isog-prop}
Let $E \ra S$ {\upshape(}resp.,~$E' \ra S${\upshape)} be a fixed {\upshape(}resp.,~variable{\upshape)} generalized elliptic curve over a scheme $S$.
\benum
\item \lab{IP-a}
If $G \subset E^\sm$ is finite locally free $S$-subgroup, then the homomorphism $E \ra E/G$ is initial among homomorphisms $f\colon E \ra E'$ with $G \subset \Ker f$.

\item \lab{IP-b}
If $f\colon E \ra E'$ is a surjective homomorphism, then $\Ker f$ is a finite locally free $S$-subgroup of $E^\sm$, and $\Ker f$ determines $f$ up to an isomorphism in the sense that $f$ induces an isomorphism 
\[
\qq E/(\Ker f) \cong E'.
\]

\item \lab{IP-c}
If $G_1 \subset G_2 \subset E^\sm$ are finite locally free $S$-subgroups, then 
\[
\qq (E/G_1)/(G_2/G_1) \cong E/G_2.
\]
\eenum
\ecor

\bpf \hfill
\benum
\item
The map $f$ is $G$-equivariant for the trivial $G$-action on $E'$, so it uniquely factors through the categorical quotient $E \ra E/G$. It remains to note that the induced map $(E/G)^\sm \ra (E')^\sm$ intertwines the group laws, as may be checked on the fppf cover $E^\sm \ra (E/G)^\sm$.

\item
The first claim was proved in Remark \ref{homo-surj}. Due to \ref{IP-a}, $f$ induces a homomorphism $E/(\Ker f) \ra E'$ that is an isomorphism on the smooth loci. Due to \cite{EGAIV4}*{17.9.5} and the $S$-flatness of $E/(\Ker f)$, checking that $E/(\Ker f) \ra E'$ is an isomorphism may be done on geometric fibers, where it follows from the fact that an endomorphism of the standard $n$-gon that is an automorphism on the smooth locus must be an automorphism.

\item
The claim follows from the universal property of $E \ra E/G_2$ recorded in \ref{IP-a}.
\qedhere
\eenum
\epf

\Cref{isog-prop}~\ref{IP-b} and the analogy with elliptic curves justify the following definition:

\bd \lab{def-isog}
An \emph{isogeny} between generalized elliptic curves $E \ra S$ and $E' \ra S$ is a surjective homomorphism $f \colon E \ra E'$ (so, by \Cref{isog-prop}~\ref{IP-b}, it induces an isomorphism $E' \cong E/(\Ker f)$). The \emph{degree} of an isogeny $f$ is the locally constant function on $S$ given by the order of $\Ker f$.
\ed

The principal difference with the elliptic curve case is that an isogeny between generalized elliptic curves is not necessarily flat (see Remark \ref{isog-non-flat}). As we explain in the following \Cref{Zar-loc-isog} (whose elliptic curve case is \cite{KM85}*{2.4.2}), the structure of an arbitrary homomorphism may be completely understood in terms of isogenies (in turn, by \Cref{isog-prop}~\ref{IP-b}, the structure of an isogeny is completely determined by its kernel).

\bprop \lab{Zar-loc-isog}
Every homomorphism $f \colon E \ra E'$ between generalized elliptic curves $E \ra S$ and $E' \ra S$ is Zariski locally on $S$ either an isogeny or the zero homomorphism.
\eprop

\bpf
Limit arguments described in Remark \ref{gen-ell-limit} allow us to reduce to the case when $S$ is Noetherian, so the claim follows from \cite{MFK94}*{Prop.~6.1}, which proves that on each connected component of $S$ the map $f$ is either surjective (i.e.,~an isogeny) or the zero homomorphism.
\epf

Due to \Cref{Zar-loc-isog}, the following result describes how homomorphisms interact with the subschemes of nonsmoothness and the degeneracy loci of \Cref{def-sing}:

\bprop \lab{isog-sing}
If $f\colon E \ra E'$ is an isogeny between generalized elliptic curves $E \xra{\pi} S$ and $E' \xra{\pi'} S$, then $f|_{E^\sing}$ factors through $E'^{\sing}$ and $S^{\infty, \pi} \subset S^{\infty, \pi'}$.
\eprop

\bpf
The second claim follows from the first because $S^{\infty, \pi}$ (resp.,~$S^{\infty, \pi'}$) is the schematic image of $E^\sing \ra S$ (resp.,~of $E'^{\sing} \ra S$). Moreover, since the formation of all the subschemes in question commutes with base change in $S$ (see Remark \ref{sing-bc}), we may use Remark \ref{gen-loc-proj} to assume that $S = S^{\infty, \pi, n}$ and that $E$ is the standard $n$-gon. 

The intersection $G$ of $\Ker f$ with the relative identity component $(E^\sm)^0 = \bG_m$ is a finite locally free $S$-subgroup scheme of both $\Ker f$ and $\bG_m$. By \Cref{isog-prop}~\ref{IP-b}--\ref{IP-c}, $f$ is identified with the composite 
\[
E \ra E/G \ra (E/G)/((\Ker f) /G)
\]
of isogenies. Therefore, since the assertion about $f|_{E^\sing}$ is compatible with composition, it suffices to treat the cases $G = \Ker f$ and $G = 0$ separately.

Since $\bG_m$ has a unique finite locally free $S$-subgroup of a given order, Zariski locally on $S$ we have $G = \mu_d$ for some $d \in \bZ_{\ge 1}$. Thus, if $G = \Ker f$, then we may assume that $f$ is the $f_d$ described in \Cref{n-gon-quot} (see also Remark \ref{quot-bc}). For this $f_d$, the claim is clear: 
\[
\tst E^\sing \qq \text{is identified with} \qq \bigsqcup_{\bZ/n\bZ} S \q \text{used in \Cref{n-gon-def}}
\]
and $f_d$ is induced by the $d\th$ power map on every $\bP^1_S$ so maps $E^\sing$ to itself.

If $G = 0$, then $f$ is \'{e}tale, so that $\Omega^1_{E/S} \cong f^*\Omega^1_{E'/S}$. By \cite{SGA7I}*{VI, 5.1 (a)}, the formation of the closed subscheme cut out by a Fitting ideal of a finite type quasi-coherent module commutes with pullback to another scheme, so this relation between the sheaves of differentials gives $E^\sing = f\i(E'^{\sing})$.
\epf

The inclusion $S^{\infty, \pi} \subset S^{\infty, \pi'}$ of \Cref{isog-sing} may be sharpened to a precise relation between the corresponding ideal sheaves. We record this in \Cref{isog-raise-it} and \Cref{explain-raise-it}.

\bprop \lab{isog-raise-it}
If $f \colon E \ra E'$ is an isogeny between generalized elliptic curves and if there is a $d \in \bZ_{\ge 1}$ such that for every degenerate geometric fiber $E_{\ov{s}}$ the intersection $(\Ker f)_{\ov{s}} \cap (E^\sm_{\ov{s}})^0$ has rank $d$, then the ideal sheaves in $\sO_S$ of the degeneracy loci $S^{\infty, \pi}$ and $S^{\infty, \pi'}$ of $E$ and $E'$ are related by
\[
\sI_{S^{\infty, \pi'}} = \sI_{S^{\infty, \pi}}^d.
\]
\eprop

\brem \lab{explain-raise-it}
For any $f$, Zariski locally on $S$ there exists a required $d$. In order to prove this, we may assume that $S = S^{\infty, \pi}$ and may work fppf locally on $S$, so Remark \ref{gen-loc-proj} reduces to the case when $E$ is the standard $n$-gon. In this case $\Ker f \cap (E^\sm)^0$ is an open and closed $S$-subgroup of $\Ker f$, and the claim follows from the local constancy of its rank over $S$.
\erem

\bpf[Proof of Proposition~{\upshape\ref{isog-raise-it}}]
It suffices to treat the case when $S = \Spec R$ for some Artinian local ring $(R, \fm)$ that has a separably closed residue field $R/\fm$. The elliptic curve case is clear, so we assume that  $E_{R/\fm}$ is degenerate. Moreover, by \Cref{isog-prop}~\ref{IP-c}, quotients may be taken in stages, so we assume that either 
\[
\Ker f \subset (E^\sm)^0 \qq \text{or}  \qq  \Ker f \cap (E^\sm)^0 = 0.
\]

We begin with the case $\Ker f \cap (E^\sm)^0 = 0$, when $f$ is finite \'{e}tale of rank $\#(\Ker f)$, so that $E^\sing = f\i(E'^\sing)$ by \cite{SGA7I}*{VI,~5.1~(a)}. \Cref{sing-map} then gives the desired $S^{\infty, \pi} = S^{\infty, \pi'}$.

In the remaining case when $\Ker f \subset (E^\sm)^0$, we first replace $S$ by a flat cover to be able to assume that there is a finite \'{e}tale $S$-subgroup $G \subset E^\sm$ such that $G_{R/\fm}$ maps isomorphically to the component group of $E_{R/\fm}^\sm$. Due to the settled $\Ker f \cap (E^\sm)^0 = 0$ case, passage to $E/G$ and $E'/f(G)$ does not affect the degeneracy loci. Therefore, we may replace 
\[
E \q \text{by}\q E/G \qq \text{and} \qq  E' \q \text{by} \q E'/f(G)
\]
to reduce to the case when $E$ is irreducible.

In this situation, since $S$ is Artinian local and strictly Henselian, \cite{DR73}*{VII.2.1} ensures that $E$ is a base change of the Tate curve 
\[
\Tate_1 \ra \Spec \bZ\llb q\rrb
\]
(\emph{loc.~cit.}~proves that $\Tate_1$ realizes $\Spec \bZ\llb q\rrb$ as an \'{e}tale double cover of the formal completion of $\EEl_1$ along $\EEl_1^\infty$; in the notation of \emph{loc.~cit.},~$\Tate_1 = \ov{\sG}_m^q/q^\bZ$). If, moreover, $\Ker f \subset (E^\sm)^0$, then $\Ker f = \mu_{\#(\Ker f)}$ inside $(E^\sm)^0$ (see \Cref{can-submult}), so that we are reduced to the case when 
\[
E \ra S \q  \text{is} \q \Tate_1 \ra \Spec \bZ\llb q\rrb \qq \text{and}\qq \Ker f = \mu_d.
\]
However, in this case the quotient map\footnote{In the notation of \cite{DR73}*{VII.1.10}, we have $\Tate_1(q^d) = \ov{\sG}_m^{q^d}/(q^d)^\bZ$ over $A = \bZ\llb q\rrb$.} $\Tate_1 \ra \Tate_1/\mu_d$ is identified with the map 
\[
\Tate_1 \ra \Tate_1(q^d) \qq \text{induced by ``raising the coordinates to the $d\th$ power,''}
\]
as in \Cref{n-gon-quot} (compare with \cite{Con07a}*{2.5.1}). It remains to recall from \cite{DR73}*{VII.1.11} that the degeneracy locus of $\Tate_1$ (resp.,~of $\Tate_1(q^d)$) is cut out by the principal ideal $(q) \subset \bZ\llb q\rrb$ (resp.,~$(q^d) \subset \bZ\llb q\rrb$).
\epf


\section{Compactifications of the stack of elliptic curves} \lab{compactify-ell}

Our approach to the study of level structures on generalized elliptic curves makes essential use of the tower $\{ \EEl_n\}_{n \mid n'}$ of compactifications of the stack $\Ell$ that parametrizes elliptic curves. The purpose of this chapter is to construct this tower and to detail its properties. We begin with the construction of the individual compactifications $\EEl_n$ in section \ref{Elln-section}, and proceed to expose the transition morphisms $\EEl_{nm} \ra \EEl_n$ in section \ref{reconstr}. Section \ref{coarse-space-section} proves that the coarse moduli space of $(\EEl_n)_S$ is the ``$j$-line'' $\bP^1_S$ for every $n$ and every scheme $S$, whereas section \ref{GAGAmama} uses the global structure of the stacks $\EEl_n$ to algebraize formal generalized elliptic curves and their homomorphisms.

\centering \subsection{The compactification $\EEl_n$ obtained by allowing $n$-gons for a fixed $n$ \nopunct} \lab{Elln-section}\hfill

\justify

The goal of this section is to detail algebro-geometric properties of the $\bZ$-stack $\EEl_n$ obtained from the stack of elliptic curves $\Ell$ by ``adjoining N\'{e}ron $n$-gons'' (see \Cref{def-ell-n}). We prove in \Cref{elln-props} that $\EEl_n$ is a proper and smooth compactification of $\Ell$. This result has already been proved over $\bZ[\f{1}{n}]$ in \cite{DR73}*{IV.2.2}, which uses deformation-theoretic methods through its reliance on \cite{DR73}*{III.1.2}. These methods require the number of the irreducible components of each geometric fiber of the generalized elliptic curve in question to be prime to the characteristic, so they do not seem to work without inverting $n$. A related difficulty is that even though the stack $\EEl_n$ is algebraic, outside the elliptic curve locus it is not Deligne--Mumford in characteristics dividing $n$ (see \Cref{elln-props}~\ref{EN-b}), so $\EEl_n$ may not possess universal deformation rings at some of its geometric points. To overcome these difficulties, we proceed indirectly by exploiting a convenient auxiliary algebraic stack $\sB_n$ whose relationship to $\EEl_n$ is described in \Cref{Bn-input}. 

We begin by defining the stack $\EEl_n$ that we are going to study and later use.

\begin{defn} \lab{def-ell-n}
For an $n \in \bZ_{\ge 1}$, let $\EEl_n$ denote the $\bZ$-stack parametrizing those generalized elliptic curves $E \xra{\pi} S$ whose degenerate geometric fibers are $n$-gons. Let $\EEl_n^\infty$ denote the closed substack of $\EEl_n$ cut out by the degeneracy loci $S^{\infty, \pi}$ (defined in \Cref{def-sing}).
\end{defn}

\brems
\remi
The effectivity of descent data that is needed for $\EEl_n$ to be a $\bZ$-stack (for the fpqc topology) results from the $S$-ampleness of the relative effective Cartier divisor $E^\sm[n] \subset E$.

\remi
The well-definedness of the closed substack $\EEl_n^\infty$ rests on the compatibility (recalled in Remark \ref{sing-bc}) of the formation of the degeneracy locus $S^{\infty, \pi}$ with base change. 
\erems

We turn to the auxiliary stack $\sB_n$ and to its relation to $\EEl_n$.

\bpp[The stack $\sB_n$] \lab{def-Bn}
Following \cite{DR73}*{V.1.3}, for an $n \in \bZ_{\ge 1}$ we let $\sB_n$ be the $\bZ$-stack that, for variable schemes $S$, parametrizes the pairs $(E, G)$ consisting of a generalized elliptic curve $E \ra S$ whose degenerate geometric fibers are $n$-gons and a finite \'{e}tale subgroup $G \subset E^\sm$ that is \'{e}tale locally on $S$ isomorphic to $\bZ/n\bZ$ and meets every irreducible component of every geometric fiber of $E \ra S$. If $n = 1$, then $G$ is the zero subgroup, so $\sB_1 = \EEl_1$.
\epp

\bprop \lab{Bn-input} 
Fix an $n \in \bZ_{\ge 1}$.
\benum
\item \lab{Bn-a}
The $\bZ$-stack $\sB_n$ is Deligne--Mumford and $\bZ$-smooth of relative dimension $1$.

\item \lab{Bn-b}
The morphism 
\[
\qq \sB_n \ra \EEl_n
\]
that forgets $G$ factors through the open substack $\EEl_n^{\text{$n$-}\ord} \subset \EEl_n$ obtained by removing the supersingular elliptic curves in characteristics dividing $n$. The induced morphism 
\[
\qq \sB_n \ra \EEl_n^{\text{$n$-}\ord}
\]
is representable by schemes, separated, quasi-finite, faithfully flat, and of finite presentation. 

\item \lab{Bn-c}
The stack $\EEl_n^{\text{$n$-}\ord}$ is algebraic and $\bZ$-smooth of relative dimension $1$.
\eenum
\eprop

\bpf \hfill
\benum
\item
Both claims follow from \cite{DR73}*{V.1.4}.

\item
The morphism 
\[
\qq q\colon \EEl_n \ra \EEl_1\qq \text{is well defined by}\qq q(E) = E/E^\sm[n]
\]
(see \S\ref{quot-not}), and, as in \cite{DR73}*{VI.1.1}, the $j$-invariant gives the morphism $j\colon \EEl_1 \ra \bP^1_\bZ$. Since $\EEl_n^{\text{$n$-}\ord}$ is the preimage under $j \circ q$ of the open subscheme of $\bP^1_\bZ$ obtained by removing the supersingular $j$-invariants in characteristics dividing $n$, it is indeed an open substack of $\EEl_n$.

The morphism $\sB_n \ra \EEl_n$ factors through $\EEl_n^{\text{$n$-}\ord}$ because a supersingular elliptic curve over an algebraically closed field of positive characteristic $p$ cannot have $\bZ/p\bZ$ as a subgroup. Therefore, our task is to prove that for any generalized elliptic curve $E \ra S$ whose geometric fibers are $n$-gons, ordinary elliptic curves in characteristic dividing $n$, or arbitrary elliptic curves in characteristic not dividing $n$, the functor
\[
\qqq F_0\colon S' \mapsto \{ \text{$S'$-ample subgroups $G \subset E^\sm_{S'}$ that are \'{e}tale locally on $S'$ isomorphic to $\bZ/n\bZ$}\}
\]
on the category of $S$-schemes is representable by a separated, quasi-finite, faithfully flat $S$-scheme $B$ of finite presentation (the $S'$-ampleness of $G$ as a relative effective Cartier divisor on $E_{S'}$ is equivalent to the condition that $G$ meets every irreducible component of every geometric fiber of $E_{S'} \ra S'$). In fact, it suffices to prove the same statement with `faithfully flat' replaced by `flat' and for the functor $F_0'$ obtained by dropping the $S'$-ampleness requirement from the definition of $F_0$: indeed, the surjectivity of $B \ra S$ will follow from the imposed fibral assumptions on $E \ra S$, whereas \cite{EGAIV3}*{9.6.4} together with limit arguments ensures that the inclusion $F_0 \subset F_0'$ is representable by quasi-compact open immersions.

We analyze $F_0'$ by studying the related functor
\[
\qqq F_1\colon S' \mapsto \{ \text{$P \in E^\sm[n](S')$ that define a closed immersion $\bZ/n\bZ \hra E^\sm_{S'}[n]$ by $1 \mapsto P$}\}.
\]
The map $F_1 \ra F_0'$ that sends $P$ to the copy of $\bZ/n\bZ$ that $P$ generates is representable by schemes and finite \'{e}tale of rank $\phi(n)$. Therefore, once we prove that $F_1$ is representable by a finitely presented, separated, quasi-finite (and hence also quasi-affine, see \cite{EGAIV3}*{8.11.2}), flat $S$-scheme, the desired claim about $F_0'$ will follow from \cite{SGA3Inew}*{V,~4.1} (combined with \cite{EGAIV2}*{2.2.11~(iii)} and \cite{EGAIV4}*{17.7.5}).

The $S$-scheme $E^\sm[n]$ represents the functor of $S'$-homomorphisms $\bZ/n\bZ \ra E^\sm_{S'}[n]$. Such a homomorphism is a closed immersion if and only if its corresponding map $f$ of finite locally free $\sO_{S'}$-algebras is surjective, which is an open condition on $S'$ because $\Coker(f)$ is a finitely generated $\sO_{S'}$-module. Therefore, the inclusion $F_1 \subset E^\sm[n]$ is representable by open immersions, and is quasi-compact by limit arguments, so the claims about $F_1$ follow.

\item
Both claims follow from \ref{Bn-b}. More precisely, if $X \ra \sB_n$ is a smooth atlas, then the composed morphism 
\[
\qq X \ra \EEl_n^{\text{$n$-}\ord}
\]
is representable by algebraic spaces, faithfully flat, and locally of finite presentation, so $\EEl_n^{\text{$n$-}\ord}$ is algebraic by \cite{SP}*{\href{http://stacks.math.columbia.edu/tag/06DC}{06DC}} (see also \cite{LMB00}*{10.6} for a related result), whereas, due to \cite{EGAIV4}*{17.7.7}, the $\bZ$-smoothness of $\EEl_n^{\text{$n$-}\ord}$ follows from that of $\sB_n$ (for the relative dimension aspect, one may use \cite{EGAIV2}*{6.1.2}).
\qedhere
\eenum
\epf

With \Cref{Bn-input} in hand, we are ready to address algebro-geometric properties of $\EEl_n$ (see \Cref{Elln-coarse} for some further properties).

\bthm \lab{elln-props} 
Fix an $n \in \bZ_{\ge 1}$.
\benum
\item \lab{EN-a}
The $\bZ$-stack $\EEl_n$ is algebraic with finite diagonal, proper, and smooth of relative dimension~$1$. 

\item \lab{EN-b}
The largest open substack of $\EEl_n$ that is Deligne--Mumford is 
\[
\qq \EEl_n - (\EEl_n^\infty)_{\bZ/n\bZ}.
\]
\item \lab{EN-c}
The morphism $\Spec \bZ \ra \EEl_n^\infty$ that corresponds to the standard $n$-gon is surjective, representable, and finite locally free of rank $2n$. In particular, the proper $\bZ$-algebraic stack $\EEl_n^\infty$ is irreducible, has geometrically irreducible $\bZ$-fibers, and is $\bZ$-smooth of relative dimension $0$.

\item \lab{EN-d}
The closed substack $\EEl_n^\infty \subset \EEl_n$ is a reduced relative effective Cartier divisor over $\Spec \bZ$.
\eenum
\ethm

\brem \lab{DM-make-sense}
In \ref{EN-b}, the largest Deligne--Mumford open substack of the separated $\bZ$-algebraic stack $\EEl_n$ does make sense \emph{a priori}. Indeed, the proof of \cite{Con07a}*{2.2.5 (2)} shows that if $S$ is a scheme and $\sX$ is an $S$-algebraic stack that is covered by $S$-separated open substacks, then there is a unique open substack 
\[
\sU \subset \sX
\]
containing exactly those geometric points of $\sX$ that have an unramified automorphism functor. (Equivalently, $\sU$ contains those $S$-scheme valued points of $\sX$ whose automorphism functors are unramified.) By Nakayama's lemma (or simply by \cite{SP}*{\href{http://stacks.math.columbia.edu/tag/02GF}{02GF~(1)$\Leftrightarrow$(2)}}), the diagonal $\Delta_{\sU/S}$ is unramified, so $\sU$ is Deligne--Mumford, and, by construction, $\sU$ contains every Deligne--Mumford open substack of $\sX$. Even though we take the unramifiedness of the diagonal as our definition of being Deligne--Mumford (see \S\ref{conv}), in the case in hand $\sU$ inherits separatedness from $\EEl_n$, so, by \cite{LMB00}*{8.1}, it also satisfies the \'{e}tale atlas definition of a Deligne--Mumford stack.
\erem

\bpf[Proof of Theorem~{\upshape\ref{elln-props}}]
\hfill
\benum
\item
The stack $\EEl_n$ is a union of open substacks $\Ell$ and $\EEl_n^{\text{$n$-}\ord}$, both of which are algebraic and $\bZ$-smooth of relative dimension $1$ by \Cref{Bn-input}. Therefore, $\EEl_n$ is also algebraic and $\bZ$-smooth of relative dimension $1$. 

By \cite{Con07a}*{3.2.4}, the isomorphism functor of two generalized elliptic curves $E \ra S$ and $E' \ra S$ whose degenerate geometric fibers are $n$-gons is representable by a finite $S$-scheme,\footnote{\lab{no-blowup} Here is a sketch for a proof of this representability that bypasses blowups used in \cite{Con07a}*{3.2.2 and 3.2.4}: as in the proof of \cite{DR73}*{III.2.5}, one uses Hilbert schemes to get representability by a quasi-finite, separated $S$-scheme; then, due to the valuative criterion, the key point is to check that if $S$ is the spectrum of a strictly Henselian discrete valuation ring and $E$ and $E'$ are degenerating elliptic curves with identified generic fibers: $E_\eta = E'_\eta$, then $E = E'$; for this, the theory of N\'{e}ron models (especially, \cite{BLR90}*{7.4/3}) identifies $(E^\sm)^0$ with $(E'^{\sm})^0$ and, since the reductions of $\eta$-rational points are dense in the special fibers, also $E^\sm$ with $E'^{\sm}$; then Zariski's main theorem \cite{BLR90}*{2.3/2$\pr$} produces the graph of the sought identification $E = E'$.}
so $\Delta_{\EEl_n/\bZ}$ is finite and, in particular, $\EEl_n$ is separated. The morphism 
\[
\qq \Ell \sqcup \Spec \bZ \ra \EEl_n
\]
whose restriction to $\Spec \bZ$ corresponds to the standard $n$-gon is surjective on underlying topological spaces, so $\EEl_n$ is quasi-compact, and hence is of finite type over $\bZ$. Its properness therefore results from the valuative criterion \cite{LMB00}*{7.10}, which is satisfied due to the semistable reduction theorem for elliptic curves (and the availability of contractions, which are reviewed in \S\ref{contr}).

\item
In the view of \Cref{DM-make-sense}, we only need to show that 
\[
\qq \EEl_n - (\EEl_n^\infty)_{\bZ/n\bZ}
\]
contains those geometric points $x$ of $\EEl_n$ whose automorphism functor is unramified. If $x$ lies in $\Ell = \EEl_n - \EEl_n^\infty$, then $\Aut(x)$ is unramified by \cite{Del75}*{5.3 (I)} (or by \cite{MFK94}*{Cor.~6.2}). If $x$ lies in $\EEl_n^\infty$, then, by \Cref{auto-n-gon}, $\Aut(x)$ is unramified if and only if $x$ lies in 
\[
\qq \EEl_n^\infty - (\EEl_n^\infty)_{\bZ/n\bZ}.
\]

\item
For the asserted properties of the morphism, it suffices to note that for a generalized elliptic curve $E \xra{\pi} S$ with $S^{\infty, \pi, n} = S$, the functor of isomorphisms between $E$ and the standard $n$-gon is representable by a finite locally free $S$-scheme of rank $2n$, as may be checked fppf locally on $S$ with the help of Remark \ref{gen-loc-proj} and \Cref{auto-n-gon}. The asserted properties of $\ov{\Ell}_n^\infty$ then follow by using \cite{EGAIV4}*{17.7.7} and \cite{EGAIV2}*{6.1.2} for the smoothness aspect.

\item
By \ref{EN-c}, the stack $\ov{\Ell}_n^\infty$ is $\bZ$-smooth, so it is also reduced. For the Cartier divisor claim, we may work over a smooth finite type scheme cover 
\[
\qq X \ra \EEl_n, \qq \text{with} \q X^\infty \subset X \q \text{being the preimage of}\q \EEl_n^\infty.
\]
By \cite{KM85}*{1.1.5.2}, we may also base change from $\bZ$ to an algebraically closed field. Then, for a point $x \in X^\infty$, by \ref{EN-a} and \ref{EN-c}, both $X$ and $X^\infty$ are smooth at $x$ and 
\[
\qq \dim_x X^\infty = \dim_x X - 1.
\]
Thus, $X^\infty \subset X$ is a Weil divisor and, since $X$ is regular, also a desired Cartier divisor.
\qedhere
\eenum
\epf

For later use we record the following proposition taken from \cite{Con07a}*{3.2.4}. 

\bprop\lab{extend-iso}
Let $E \xra{\pi} S$ and $E' \xra{\pi'} S$ be generalized elliptic curves such that 
\[
\qqq S^{\infty, \pi, n} \cap S^{\infty, \pi', m} = \emptyset \qq \text{whenever $n \neq m$.}
\]
\benum
\item \lab{EI-a}
The fppf sheaf $\Isom(E, E')$ that parametrizes generalized elliptic curve isomorphisms is representable by a finite $S$-scheme of finite presentation.

\item \lab{EI-b}
If $S$ is integral and normal and $\eta$ is its generic point, then any $\eta$-isomorphism 
\[
\qq E_\eta \simeq E'_\eta  \qq \text{extends to a unique $S$-isomorphism} \qq E \simeq E'.
\]
\eenum
\eprop

\bpf
Part \ref{EI-a} has essentially been proved in \cref{no-blowup}. Alternatively, 
Zariski locally on $S$ there is an $n \in \bZ_{\ge 1}$ such that $E$ and $E'$ correspond to objects of $\EEl_n$, so \ref{EI-a} is a reformulation of the finiteness of the diagonal of $\EEl_n$ proved in \Cref{elln-props}~\ref{EN-a}. To obtain \ref{EI-b} one combines \ref{EI-a} with the following useful lemma.
\epf

\blem
If $S$ is an integral normal scheme, $\eta$ is its generic point, and $F$ is a finite $S$-scheme, then the pullback map $F(S) \ra F(\eta)$ is bijective.
\elem

\bpf
The injectivity follows from the schematic dominance of $\eta \ra S$ and the separatedness of $F \ra S$. For the surjectivity, we may work Zariski locally on $S$ to assume that $S = \Spec A$. Then the schematic image in $F$ of an $x \in F(\eta)$ is $\Spec B$ for some finite $A$-subalgebra $B \subset \Frac A$. Since $A$ is normal, $A = B$, so the schematic image is the sought extension of $x$ to an element of $F(S)$.
\epf


\centering \subsection{The tower of compactifications \nopunct} \lab{reconstr} \hfill

\justify

The compactifications $\EEl_n$ introduced in the previous section are related to each other: they form an infinite tower in which the transition morphisms 
\[
\EEl_{nm} \ra \EEl_n
\]
encode contractions of generalized elliptic curves. The goal of this section is to use the already established results about $\EEl_n$ to prove several basic properties, such as flatness, of these transition morphisms (see \Cref{Bnm-input}) and to deduce some concrete results about the generalized elliptic curves themselves (see \Cref{raise-it,pump-it-up}). We begin with a brief review of contractions.

\bpp[Contraction with respect to a finite locally free subgroup] \lab{contr}
As is justified in \cite{Con07a}*{top of p.~218} (which is based on \cite{DR73}*{IV.1.2}), if $E \ra S$ is a generalized elliptic curve and $G \subset E^\sm$ is a finite locally free $S$-subgroup, then there is a generalized elliptic curve 
\be\lab{contr-eq}
c_G(E) \ra S \qq \text{equipped with a surjective $S$-scheme morphism} \qq E \ra c_G(E)
\ee
such that 
\begin{itemize}
\item 
the image under $E \ra c_G(E)$ of each disjoint from $G$ irreducible component of a geometric fiber of $E \ra S$ is a single point, and

\item
the map $E \ra c_G(E)$ restricts to a group isomorphism between the open complement of the union of such components and $(c_G(E))^\sm$.
\end{itemize}
In particular, if $E$ is an elliptic curve, then $E = c_G(E)$.

These conditions ensure that $G$ is identified with an $S$-subgroup of $c_G(E)^\sm$ that meets every irreducible component of every geometric fiber of $c_G(E) \ra S$. Due to \cite{DR73}*{IV.1.2}, they also determine the data \eqref{contr-eq} uniquely up to a unique isomorphism. In particular, whenever $G' \subset E^\sm$ is another finite locally free $S$-subgroup that meets the same irreducible components of the geometric fibers of $E \ra S$ as $G$, one gets a canonical identification 
\be \lab{G-Gpr-id}
c_G(E) = c_{G'}(E).
\ee
For the same reason, the formation of $E \ra c_G(E)$ commutes with arbitrary base change in $S$.

We call this $c_G(E)$ \emph{the contraction of $E$ with respect to $G$}. The compatibility of the formation of $c_G(E)$ with base change shows that for every $n, m\in \bZ_{\ge 1}$, the identity map on $\Ell$ extends to the ``contraction'' $\bZ$-morphism
\[
\EEl_{nm} \ra \EEl_n \qq \text{defined by} \qq E \mapsto c_{E^\sm[n]}(E). \qq
\]
Likewise, if $(E, G)$ is classified by the stack $\sB_{nm}$ of \S\ref{def-Bn}, then $(c_{G[n]}(E), G[n])$ is classified by the stack $\sB_n$, so there is the ``contraction'' $\bZ$-morphism
\[
\qqqqq\ \,\sB_{nm} \ra \sB_n \qq \text{defined by} \qq (E, G) \mapsto (c_{G[n]}(E), G[n]). \qq
\]
These and similar morphisms will be called \emph{contractions} or \emph{contraction morphisms} in the sequel (a slight abuse of terminology because it is not substacks of $\EEl_{nm}$ or $\sB_{nm}$ that are getting contracted).
\epp

In many situations, we will need a robust criterion for recognizing algebraic spaces and morphisms that are representable by algebraic spaces. The following lemma, which paraphrases \cite{Con07a}*{2.2.5~(1) and 2.2.7} and may be traced back to \cite{DR73}*{IV.2.6}, is well suited for this task.

\blem \lab{rep-crit}
Let $S$ be a scheme and let $\sX$ and $\sY$ be $S$-algebraic stacks whose diagonals $\Delta_{\sX/S}$ and $\Delta_{\sY/S}$ are quasi-compact and separated.
\benum
\item \lab{RC-a}
The stack $\sX$ is an algebraic space if and only if for every algebraically closed field $\ov{k}$ whose spectrum is equipped with a morphism to $S$, every object $\xi$ of $\sX(\ov{k})$, and every Artinian local $\ov{k}$-algebra $A$, the pullback of $\xi$ to the groupoid $\sX(A)$ has no nonidentity automorphism{\upshape;} if $\sX$ is Deligne--Mumford, then $A = \ov{k}$ suffices.

\item \lab{RC-b}
An $S$-morphism 
\[
\qq f\colon \sX \ra \sY
\]
is representable by algebraic spaces if and only if for every algebraically closed field $\ov{k}$ whose spectrum is equipped with a morphism to $S$, every object $\xi$ of $\sX(\ov{k})$, and every Artinian local $\ov{k}$-algebra $A$, no nonidentity automorphism of the pullback of $\xi$ to $\sX(A)$ is sent to an identity automorphism in $\sY(A)${\upshape;} if $\sX$ is Deligne--Mumford, then $A = \ov{k}$ suffices.
\eenum
\elem

\bpf \hfill
\benum
\item
The necessity is clear. For the sufficiency, due to \cite{Con07a}*{2.2.5~(1)}, it is enough to argue that the assumed condition implies the triviality of the automorphism functor of every $\xi$. This functor is a separated $\ov{k}$-group algebraic space $G$ of finite type, so is   necessarily a scheme due to \cite{Art69a}*{4.2}, and is even $\ov{k}$-\'{e}tale if $\sX$ is Deligne--Mumford. The triviality of $G$ is therefore equivalent to that of all the $G(A)$, with $A = \ov{k}$ being sufficient if $\sX$ is Deligne--Mumford.

\item
The failure of the condition on $\xi$ implies that the groupoid of $A$-points of some $A$-fiber of $f$ has a nonidentity automorphism, and the necessity follows. For the sufficiency, due to \cite{Con07a}*{2.2.7}, it is enough to argue that the assumed condition implies that each $\ov{k}$-fiber $X$ of $f$ is an algebraic space, so  it remains to observe that this condition ensures that $X$ meets the criterion of \ref{RC-a}. 
\qedhere
\eenum
\epf

To infer further representability by schemes, we will often use the following well-known lemma:

\blem \lab{rep-schemes}
For stacks $\sX$ and $\sY$ over a scheme $S$, an $S$-morphism $f\colon \sX \ra \sY$ that is representable by algebraic spaces, separated, and locally quasi-finite is representable by schemes{\upshape;} if, in addition, $f$ is proper, then $f$ is finite.
\elem

\bpf
This follows from \cite{LMB00}*{A.2} (see also \cite{Con07a}*{2.2.6}) and \cite{EGAIV4}*{18.12.4}.
\epf

We are ready to exploit the relationship between the two contractions introduced in \S\ref{contr} to extract further information about the stacks $\EEl_n$.

\bthm \lab{Bnm-input}
For $\sB_n$ as in {\upshape\S\ref{def-Bn}} and any $n, m \in \bZ_{\ge 1}$, consider the commutative diagram 
\[
\xymatrix{
\sB_{nm} \ar[d]_-{c'} \ar[r]^-{a} & \EEl_{nm} \ar[d]_-{c} \\
\sB_n \ar[r]^-{b} & \EEl_n
}
\]
in which $c$ and $c'$ are the contraction morphisms of {\upshape\S\ref{contr}} and $a$ and $b$ forget the subgroup $G$.
\benum
\item \lab{Bnm-a}
The contractions $c$ and $c'$ are flat and of finite presentation. Moreover, $c$ is proper, with finite diagonal, and surjective, whereas $c'$ is representable by schemes, separated, and quasi-finite.

\item \lab{Bnm-b}
The closed substack
\[
\qq \EEl_n^\infty \times_{\EEl_n, c} \EEl_{nm} \subset \EEl_{nm}
\]
is a relative effective Cartier divisor over $\Spec \bZ$ that is a positive integer multiple of $\EEl_{nm}^\infty$.

\item \lab{Bnm-c}
The multiple needed in {\upshape\ref{Bnm-b}} is $m$, i.e.,
\[
\q [\EEl_n^\infty \times_{\EEl_n, c} \EEl_{nm}] = m\cdot [\EEl_{nm}^\infty]
\]
as Cartier divisors on $\EEl_{nm}$. 
\eenum
\ethm

\bpf 
The commutativity of the diagram follows from the identification discussed in \S\ref{contr}.

By \Cref{Bn-input}~\ref{Bn-b}, the maps $a$ and $b$ are representable by schemes, separated, quasi-finite, of finite presentation, flat, and faithfully flat onto $\EEl_{nm}^{nm\text{-ord}}$ and $\EEl_{n}^{n\text{-ord}}$, respectively.
\benum
\item
By \Cref{elln-props}~\ref{EN-a}, the stacks $\EEl_{nm}$ and $\EEl_n$ are $\bZ$-proper with finite diagonal, so $c$ is also proper, with finite diagonal, and of finite presentation. Since the contraction of the standard $nm$-gon with respect to its $n$-torsion is the standard $n$-gon, $c$ is surjective.
Moreover, $c|_{\Ell}$ is the identity, $\Ell$ and $\EEl_{nm}^{nm\text{-ord}}$ cover $\EEl_{nm}$, and, by \Cref{Bn-input}~\ref{Bn-b}, $a$ is faithfully flat onto $\EEl_{nm}^{nm\text{-ord}}$, so the flatness of $c$ will follow once we establish that of $c'$.

It remains to establish the claims about $c'$. For the representability of $c'$ by algebraic spaces, due to \Cref{rep-crit}~\ref{RC-b}, it suffices to observe that if $E$ is the standard $nm$-gon over an algebraically closed field and $G \simeq \bZ/nm\bZ$ is a subgroup of $E^\sm$ that meets every irreducible component of $E$, then, by \Cref{auto-n-gon}, no nonidentity automorphism of $(E, G)$ restricts to the identity map on $(E^\sm)^0$. The separatedness of $c'$ follows from the separatedness of $b \circ c' = c \circ a$ and of $b$, and similarly for the finite presentation of $c'$. For the quasi-finiteness of $c'$ it therefore suffices to observe that a generalized elliptic curve over an algebraically closed field has finitely many subgroups of order $nm$. The representability of $c'$ by schemes follows from \Cref{rep-schemes}.

Finally, since $c'$ is a quasi-finite map between separated Deligne--Mumford stacks that are smooth of relative dimension $1$ over $\bZ$, it is flat by \cite{EGAIV2}*{6.1.5}.

\item
Since $c$ is flat by \ref{Bnm-a} and $\EEl_n^\infty \subset \EEl_{n}$ is a relative effective Cartier divisor over $\Spec \bZ$ by \Cref{elln-props}~\ref{EN-d}, the pullback in question is also a relative effective Cartier divisor over $\Spec \bZ$. Both 
\[
\EEl_n^\infty \times_{\EEl_n, c} \EEl_{nm} \qq  \text{and} \qq \EEl_{nm}^\infty
\]
are supported on the same closed subset of the underlying topological space of $\EEl_{nm}$ and, by \Cref{elln-props}~\ref{EN-c}--\ref{EN-d}, this subset is irreducible and has $\EEl_{nm}^\infty$ as its associated reduced closed substack (see \cite{LMB00}*{5.6.1~(ii)}). Moreover, $\EEl_{nm}$ is regular, so on any of its scheme atlases Cartier divisors identify with Weil divisors. Passage to such an atlas then shows that $\EEl_n^\infty \times_{\EEl_n, c} \EEl_{nm}$ is a positive integer multiple of $\EEl_{nm}^\infty$---the global constancy of the needed factor across the irreducible components of the pullback of $\EEl_{nm}^\infty$ to the atlas follows from the irreducibility of $\EEl_{nm}^\infty$ (to check that the generic points of such irreducible components map to the generic point of $\EEl_{nm}^\infty$, one uses the fact that generizations lift along a flat morphism; see \cite{LMB00}*{5.8}). 

\item
Due to \ref{Bnm-b} and the moduli interpretation, it suffices to find a single generalized elliptic curve $E \xra{\pi} S$ with $nm$-gon degenerate geometric fibers such that its contraction $E' \xra{\pi'} S$ with respect to $E^\sm[n]$ satisfies the equality 
\[
\qqqq \sI_{S^{\infty, \pi'}} = \sI_{S^{\infty, \pi}}^d  \qq \text{of $\sO_S$-ideal sheaves for} \qq d = m
\]
 but does not satisfy this equality for any other $d \in \bZ_{\ge 1}$
(here $\sI_{S^{\infty, \pi}} \subset \sO_S$ is the ideal sheaf that cuts out the degeneracy locus $S^{\infty, \pi} \subset S$, and likewise for $\sI_{S^{\infty, \pi'}}$).  Tate curves supply such an $E$, see \cite{DR73}*{VII.1.11 and VII.1.14}.
\qedhere
\eenum
\epf

We now record some concrete consequences of our analysis of the contraction $c\colon \EEl_{nm} \ra \EEl_n$.

\bcor \lab{raise-it}
For a generalized elliptic curve $E \xra{\pi} S$, let $\sI_{S^{\infty, \pi}} \subset \sO_S$ be the ideal sheaf that cuts out the degeneracy locus $S^{\infty, \pi} \subset S$. If the degenerate geometric fibers of $E \xra{\pi} S$ are $nm$-gons and $c_{E^\sm[n]}(E) \xra{\pi'} S$ is the contraction of $E \xra{\pi} S$ with respect to $E^\sm[n]$, then
\[
\sI_{S^{\infty, \pi'}} = \sI_{S^{\infty, \pi}}^m.
\]

\ecor

\bpf
This is a reformulation of \Cref{Bnm-input}~\ref{Bnm-c}.
\epf

\bcor \lab{pump-it-up} 
For each $n \in \bZ_{\ge 1}$, every generalized elliptic curve $E \ra S$ is fppf locally on $S$ the contraction {\upshape(}with respect to some subgroup{\upshape)} of a generalized elliptic curve $E' \ra S$ for which the number of irreducible components of each degenerate geometric fiber is divisible by $n$. An fppf cover of $S$ over which such an $E'$ exists may be chosen to be locally quasi-finite over $S$.
\ecor

\bpf
We may assume that there is a $d \in \bZ_{\ge 1}$ such that the degenerate geometric fibers of $E$ are $d$-gons (see Remark \ref{gen-loc-proj}). The first claim then follows from flatness, surjectivity, and finite presentation of $\EEl_{nd} \xra{c} \EEl_{d}$. The second claim follows from the first and \cite{EGAIV4}*{17.16.2}.
\epf

We conclude the section by using \Cref{pump-it-up} to analyze the torsion subgroups of a generalized elliptic curve in a formal neighborhood of the degeneracy locus. Similar analysis in the case of Tate curves has been carried out in \cite{DR73}*{VII.1.13--VII.1.15}. 

\bprop \lab{Gn-Hn}
Let $E \xra{\pi} S$ be a generalized elliptic curve with $S = \Spec R$ for a Noetherian $R$ that is complete and separated with respect to the ideal $I \subset R$ that cuts out $S^{\infty, \pi} \subset S$. 
\benum
\item \lab{Gn-Hn-a}
For every $n \in \bZ_{\ge 1}$, the $S$-group $(E^\sm)^0$ has a unique finite locally free $S$-subgroup $B_n \subset (E^\sm)^0$ of order $n$, and  $B_n \simeq \mu_n$ \'{e}tale locally on $S$. If an $m \in \bZ_{\ge 1}$ divides both $n$ and the number of irreducible components of each degenerate geometric fiber of $E$, then $E^\sm[n]$ has a unique finite locally free $S$-subgroup $A_{n, m}$ that meets precisely $m$ irreducible components of every degenerate geometric fiber of $E$, contains every other finite locally free $S$-subgroup of $E^\sm[n]$ with this property, is of order $nm$, and fits into a short exact sequence
\[
\q 0 \ra B_n \ra A_{n, m} \ra C_m \ra 0
\]
with $C_m$ isomorphic to $\bZ/m\bZ$ \'{e}tale locally on $S$.

\item \lab{Gn-Hn-b}
For every $n \in \bZ_{\ge 1}$, over $S - S^{\infty, \pi}$ the group $B_n$ from \ref{Gn-Hn-a} fits into a short exact sequence
\[
\q 0 \ra (B_n)_{S - S^{\infty, \pi}} \ra E_{S - S^{\infty, \pi}}[n] \ra C_n' \ra 0
\]
with $C_n'$ an $(S - S^{\infty, \pi})$-group scheme that is  isomorphic to $\bZ/n\bZ$ \'{e}tale locally on $S - S^{\infty, \pi}$.
\eenum
\eprop

\bpf \hfill
\benum
\item
If $S$ is an infinitesimal thickening of $S^{\infty, \pi}$, then \Cref{can-submult} gives the claim. Therefore, the uniqueness and the existence of $B_n$ and $A_{n, m}$ follow from \cite{EGAIII1}*{5.1.4 and 5.4.1} (the $S$-group structure of $B_n$ may be read off from the action morphism $B_n \times_S E \ra E$, so the nonproperness of $E^\sm$ does not intervene, and likewise for $A_{n, m}$). The \'{e}tale local structure of $B_n$ translates into that of its Cartier dual, so it may be read off on geometric fibers at points in $S^{\infty, \pi}$, and likewise for the \'{e}tale local structure of $C_m$. 

\item
In the case when $n$ divides the number of irreducible components of each degenerate geometric fiber of $E$, the claim follows from \ref{Gn-Hn-a}. In general, $C_n'$ is a finite locally free $(S - S^{\infty, \pi})$-group scheme of order $n$ and it suffices to check that its geometric fibers are isomorphic to $\bZ/n\bZ$. In order to check this at a point $\eta \in S - S^{\infty, \pi}$, we choose a specialization $s \in S^{\infty, \pi}$ of $\eta$ and use \cite{EGAII}*{7.1.9} to find an $S$-scheme $T$ that is the spectrum of a complete discrete valuation ring whose generic (resp.,~closed) point maps to $\eta$ (resp.,~$s$). Due to the uniqueness of $B_n$, the formation of $C'_n$ commutes with base change of $E$ to $T$, so we are reduced to the case when $S = \Spec R$ for some complete discrete valuation ring $R$ and $I \subset R$ is a nonzero power of the maximal ideal. In this case, \Cref{pump-it-up} and \cite{EGAIV4}*{18.5.11 (a)$\Leftrightarrow$(c)} supply a finite faithfully flat $R$-algebra $R'$ such that $E_{R'}$ is the contraction of a generalized elliptic curve $E' \ra \Spec R'$ for which $n$ divides the number of irreducible components of each degenerate geometric fiber. Base change to $R'$ reduces the claim to the settled case of $E'$. \qedhere
\eenum
\epf


\centering \subsection{The coarse moduli space of $\EEl_n$ \nopunct} \lab{coarse-space-section} \hfill

\justify

We seek to prove in \Cref{Elln-coarse} that for any scheme $S$ and any $n \in \bZ_{\ge 1}$ the coarse moduli space of $(\EEl_n)_S$ is isomorphic to $\bP^1_S$, the ``$j$-line.'' Of course, this is hardly surprising, but even in the $n = 1$ case we are not aware of a reference that would treat arbitrary $S$---for $n = 1$, \cite{DR73}*{VI.1.1} settles the basic case $S = \Spec \bZ$, whereas \cite{FO10}*{2.1} handles general locally Noetherian $S$ (the formation of the coarse moduli space need not commute with nonflat base change, so the $S = \Spec \bZ$ case does not automatically imply the general case). We will build on the above result of Deligne and Rapoport through the following lemma. 

The existence of all the coarse moduli spaces that we will consider in this section is guaranteed by \cite{KM97}*{1.3~(1)} (see also \cite{Con05}*{1.1} and \cite{Ryd13}*{6.12}).

\blem \lab{coarse-miracle}
Let $\sX$ be a Deligne--Mumford stack that is separated, flat, and locally of finite type over $\bZ$, and let 
\[
f\colon \sX \ra X
\]
be its coarse moduli space map. If $f_{\bF_p}\colon \sX_{\bF_p} \ra X_{\bF_p}$ is the coarse moduli space map of $\sX_{\bF_p}$ for every prime $p$, then $f_S\colon \sX_S \ra X_S$ is the coarse moduli space map of $\sX_S$ for every scheme $S$. 
\elem

\bpf
The formation of the coarse moduli space $f\colon \sX \ra X$ commutes with flat base change in $X$, and we may work fppf locally on $X_S$ when checking that $f_S\colon \sX_S \ra X_S$ is the coarse moduli space of $\sX_S$. We may therefore assume that $S = \Spec R$ for some ring $R$ and, by \cite{AV02}*{2.2.3 and its proof}, that 
\[
X = \Spec A \qq \text{and}  \qq \sX = [(\Spec B)/G]
\]
for some finite $A$-algebra $B$ equipped with an action of a finite group $G$. In this situation, as is explained in \cite{Con05}*{3.1}, we have $A = B^G$, the coarse moduli space of $\sX_S$ is $\Spec((B \tensor_\bZ R)^G)$, and we seek to prove that the map
\[
j_R\colon B^G \tensor_\bZ R \ra (B \tensor_\bZ R)^G
\]
is an isomorphism granted that it is an isomorphism whenever $R = \bF_p$ for any $p$. 

The $\bZ$-flatness of $\sX$ ensures that $B$ is torsion-free, so the abelian group $B/B^G$ is also torsion-free. Therefore, $B^G \tensor_\bZ R \ra B \tensor_\bZ R$, and hence also $j_R$, is injective for every $\bZ$-module $R$. In order to conclude, we will prove that $j_R$ is also surjective for every $\bZ$-module $R$.

By passage to a filtered direct limit, we may assume that the $\bZ$-module $R$ is finitely generated. Thus, since the case $R = \bZ$ is clear, we may assume that $R = \bZ/n\bZ$ for some $n \in \bZ_{\ge 1}$. To then finally reduce to the assumed $R = \bZ/p\bZ$ case by devissage, it remains to use the commutative diagram
\[
\xymatrix{
0 \ar[r] & B^G \tensor_\bZ R' \ar[r] \ar@{^(->}[d]_-{j_{R'}} & B^G \tensor_\bZ R \ar[r] \ar@{^(->}[d]_-{j_{R}} & B^G \tensor_\bZ R'' \ar[r] \ar@{^(->}[d]_-{j_{R''}} & 0 \\
0 \ar[r] & (B \tensor_\bZ R')^G \ar[r] & (B \tensor_\bZ R)^G \ar[r] & (B \tensor_\bZ R'')^G}
\]
that is in place whenever one has a short exact sequence $0 \ra R' \ra R \ra R'' \ra 0$ of $\bZ$-modules.
\epf

We are ready for the promised conclusion about the coarse moduli space of $(\EEl_n)_S$.

\bprop \lab{Elln-coarse}
For any $n \in \bZ_{\ge 1}$, the coarse moduli space of $\EEl_n$ (resp.,~of the open substack $\Ell \subset \EEl_n$) is isomorphic to $\bP^1_\bZ$ (resp.,~to $\bA^1_\bZ \subset \bP^1_\bZ$, with the map $\Ell \ra \bA^1_\bZ$ being given by the $j$-invariant) and its formation commutes with base change to an arbitrary scheme $S$. In particular, $\EEl_n$ is irreducible and has geometrically irreducible $\bZ$-fibers.
\eprop

\bpf
The last assertion follows from the rest because  the map to the coarse moduli space induces a homeomorphism on topological spaces. 

We begin with the $n = 1$ case, for which the base $S = \Spec \bZ$ has been treated  in \cite{DR73}*{VI.1.1 and VI.1.3} and we only need to prove that the formation of the coarse moduli space of $\EEl_1$ commutes with arbitrary base change. Let 
\[
\sC \subset \EEl_1
\]
be the preimage of the open subscheme of $\mathbb{P}^1_{\bZ}$ obtained by removing the sections $j = 0$ and $j = 1728$. By \cite{Del75}*{5.3~(III)}, the automorphism functor of every generalized elliptic curve classified by $\sC$ is the constant group $\{\pm 1\}$.  Therefore, as is explained in \cite{ACV03}*{\S5.1}, \cite{Rom05}*{\S5}, or \cite{AOV08}*{Appendix~A}, we may ``quotient out'' this constant group from the automorphism functors to obtain the algebraic stack $\sC\!\!\! \fatslash \{\pm 1\}$ that is a ``rigidification''  of $\sC$. By, for instance, \cite{AOV08}*{A.1}, the rigidification map 
\[
\sC \ra \sC\!\!\! \fatslash \{\pm 1\}
\]
induces an isomorphism on coarse moduli spaces. However, by \cite{LMB00}*{8.1.1}, the algebraic stack $\sC\!\!\! \fatslash \{\pm 1\}$ is its own coarse moduli space. Thus, since the formation of $\sC\!\!\! \fatslash \{\pm 1\}$ commutes with arbitrary base change, so does that of the coarse moduli space of $\sC$. In particular, for every prime $p$, the map from the coarse moduli space of $(\EEl_1)_{\bF_p}$ to $\bP^1_{\bF_p}$ is an isomorphism on a dense open subscheme. However, this map is finite locally free due to the normality of its source inherited from the $\bF_p$-smooth $(\EEl_1)_{\bF_p}$, so it is an isomorphism globally. This settles the $n = 1$ case for $S = \Spec \bF_p$, and the general $n = 1$ case then follows from \Cref{coarse-miracle}.

For general $n$, we begin by arguing that the coarse moduli space $Y$ of $\EEl_n$ is $\bZ$-flat and that its formation commutes with arbitrary base change. By the settled $n = 1$ case, this is true on the elliptic curve locus, so we may focus on the open substack $\sC_n \subset \EEl_n$ that is the preimage of $\sC$. By \cite{DR73}*{II.2.8}, every generalized elliptic curve has the automorphism $-1$ that restricts to  inversion on the smooth locus. In particular, the constant group scheme $\{ \pm 1\}$ is a canonical subgroup functor of the automorphism functor of every generalized elliptic curve classified by $\sC_n$, so we may pass to the rigidification $\sC_n\!\! \fatslash \{\pm 1\}$ and need to argue that its coarse moduli space is $\bZ$-flat and of formation compatible with base change. This follows from \cite{AOV08}*{3.3} because the algebraic stack $\sC_n\!\! \fatslash \{\pm 1\}$ is tame by \Cref{auto-n-gon} and \cite{Del75}*{5.3~(III)}.

It remains to prove that the map $f\colon Y \ra \bP^1_\bZ$ between the coarse moduli spaces of $\EEl_n$ and $\EEl_1$ is an isomorphism. By \cite{Ryd13}*{6.12}, the coarse moduli space $Y$ is $\bZ$-proper, so the map in question is proper and quasi-finite, and hence also finite by \Cref{rep-schemes}. Once we prove its flatness, and hence also local freeness, it will remain to inspect the elliptic curve locus to see that it is an isomorphism. Due to the $\bZ$-flatness of $Y$ and \cite{EGAIV3}*{11.3.11}, for the remaining flatness of $f$ we may work $\bZ$-fiberwise, and hence conclude with the help of \cite{EGAIV2}*{6.1.5} after observing that for every field $k$, the reducedness of the $k$-smooth $(\EEl_n)_k$ ensures the reducedness, and hence also the Cohen--Macaulay property, of its $1$-dimensional coarse moduli space $Y_k$.
\epf


\centering \subsection{Algebraization of formal generalized elliptic curves and of their homomorphisms \nopunct} \lab{GAGAmama} \hfill

\justify

The goal of this section is to prove that a formal generalized elliptic curve that is adic over an affine Noetherian formal scheme and whose number of irreducible components of a degenerate geometric fiber is constant may be uniquely algebraized, and likewise for generalized elliptic curve homomorphisms---see \Cref{GAGA} for a precise statement. Such algebraizability does not immediately follow  from Grothendieck's formal GAGA formalism because the loci of smoothness may not be proper over the base, but it nevertheless is not surprising: if this formalism applied to the $\bZ$-proper stack $\EEl_n$ as it does in the scheme case, then the pullback map
\[
\tst \EEl_n(R) \ra \varprojlim_m \EEl_n(R/I^m)
\]
would be an equivalence for every adic Noetherian ring $R$ with an ideal of definition $I$, and \Cref{GAGA}~\ref{GAGA-a} would follow. The key difference from the scheme case is that  a section of $(\EEl_n)_R \ra \Spec R$ is not a closed immersion. Nevertheless, an argument that we have extracted from \cite{Ols06}*{5.4} proves a suitable formal GAGA statement recorded in \Cref{formal-sec} (see also \cite{Aok06}*{\S3.4} and \cite{Aok06e} for a similar argument).

\blem \lab{formal-sec}
Let $R$ be a Noetherian ring that is complete and separated with respect to an ideal $I \subset R$. For every proper $R$-algebraic stack $\sX$ with finite diagonal $\Delta_{\sX/R}$ {\upshape(}for instance, for every proper Deligne--Mumford $R$-stack $\sX${\upshape)}, the functor
\be \lab{FS-map}
\tst \sX(R) \ra \varprojlim_m \sX(R/I^m)
\ee
is an equivalence of categories.
\elem

\bpf
If $x, x' \in \sX(R)$, then the isomorphism functor $\Isom(x, x')$ is a finite $R$-scheme, so
\[
\tst \Isom(x, x')(R) \ra \varprojlim_m \Isom(x, x')(R/I^m)
\]
is bijective by formal GAGA for schemes \cite{EGAIII1}*{5.1.6}. In other words, the functor \eqref{FS-map} is fully faithful. For its essential surjectivity, suppose that 
\[
\{x_m \in \sX(R/I^m)\}_{m \ge 1}
\]
is a compatible sequence of objects. Due to the finiteness of $\Delta_{\sX/R}$, each map
\[
\Spec(R/I^m) \xra{x_m} \sX_{R/I^m}
\]
is representable by schemes and finite. Therefore, $x_m$ corresponds to a coherent $\sO_{\sX_{R/I^m}}$-algebra $\sA_m$. By formal GAGA for Artin stacks, i.e.,~by \cite{Ols06}*{A.1}, the compatible system $\{ \sA_m \}_{m\ge 1}$ comes via base change from a unique coherent $\sO_{\sX}$-algebra $\sA$. It remains to argue that the composition of the finite morphism $X \xra{x} \sX$ corresponding to $\sA$ and the structure morphism  $\sX \ra \Spec R$ is an isomorphism. By construction, $x_{R/I^m} = x_m$ for every $m \ge 1$, so the claim will follow from \cite{EGAIII1}*{5.1.6} once we prove that the proper $R$-algebraic stack $X$ is a finite $R$-scheme.

By \cite{Con07a}*{2.2.5~(2)}, the algebraic space locus of $X$ is open and contains $X_{R/I}$, so it must coincide with $X$. Since the relative dimension of $X$ over $R$ may be computed \'{e}tale locally on $X$, \cite{EGAIV3}*{13.1.3} proves that the relative dimension $0$ locus of $X$ is open, and hence must equal $X$ because it contains $X_{R/I}$. To conclude that $X \ra \Spec R$ is finite one then applies \Cref{rep-schemes}.
\epf

The algebraization \Cref{GAGA}~\ref{GAGA-a} has already been proved in \cite{Con07a}*{2.2.4} by a different argument that does not use formal GAGA for Artin stacks (a similar argument had previously been used in \cite{DR73}*{VII.1.10} to construct Tate curves), but it seems worthwhile to put this result in the context of \Cref{formal-sec}. In contrast, the method of \cite{Con07a}*{2.2.4} does not seem to suffice for the proof of the algebraizability of homomorphisms (beyond the case of isomorphisms), i.e.,~for \Cref{GAGA}~\ref{GAGA-b}. To algebraize homomorphisms we exploit their structure detailed in \S\ref{quotient-section}.

\bthm \lab{GAGA}
Let $R$ be a Noetherian ring, complete and separated with respect to an ideal $I \subset R$. 
\benum
\item \lab{GAGA-a}
For each $n \in \bZ_{\ge 1}$, every compatible under pullback sequence 
\[
\qq \{E_m \ra \Spec (R/I^m)\}_{m \ge 1}
\]
of generalized elliptic curves whose degenerate geometric fibers are $n$-gons is isomorphic to the sequence obtained via base change from a unique generalized elliptic curve $E \ra \Spec R$.

\item \lab{GAGA-b}
For generalized elliptic curves $E \ra \Spec R$ and $E' \ra \Spec R$, every compatible sequence
\[
\qq \{f_m\colon E_{R/I^m} \ra E'_{R/I^m} \}_{m \ge 1}
\]
of generalized elliptic curve homomorphisms {\upshape(}defined in Definition~{\upshape\ref{homo-def})} comes via base change from a unique generalized elliptic curve homomorphism $f\colon E \ra E'$.
\eenum
\ethm

\bpf \hfill
\benum
\item
\Cref{formal-sec} applied to $\EEl_n$ proves the claim (for the uniqueness, Remark \ref{gen-loc-proj} ensures that the degenerate geometric fibers of $E$ are $n$-gons).

\item
We begin with the case when all the $f_m$ are isomorphisms (\Cref{formal-sec} does not apply because $E$ need not correspond to an object of $\EEl_n$ for any $n$). Due to Remark \ref{gen-loc-proj}, there is no geometric point $\ov{s}$ of $\Spec R$ for which $E_{\ov{s}}$ and $E'_{\ov{s}}$ are both degenerate but have distinct numbers of irreducible components, so \Cref{extend-iso}~\ref{EI-a} shows that the isomorphism functor $\Isom(E, E')$ is a finite $R$-scheme. Therefore, by \cite{EGAIII1}*{5.1.6}, the sequence 
\[
\qqq \tst (f_m) \in \varprojlim_m \Isom(E, E')(R/I^m) \q \text{is induced by a desired unique} \q f \in \Isom(E, E')(R).
\]

In the general case, by \cite{EGAIII1}*{5.4.1}, the $f_m$ algebraize to a unique $R$-morphism 
\[
\qq f \colon E \ra E',
\] 
and our task is to show that $f$ is a generalized elliptic curve homomorphism. Since idempotents of $R/I$ lift uniquely to $R$ (see \cite{EGAIV4}*{18.5.16~(ii)}), we may use \Cref{Zar-loc-isog} to write  
\[
\qq R = R' \times R'' \qq \text{and} \qq I = I' \times I''
\]
in such a way that $(f_1)_{R'/I'}$ is the zero homomorphism and $(f_1)_{R''/I''}$ is an isogeny. Then $R'$ (resp.,~$R''$) is complete and separated with respect to $I'$ (resp.,~$I''$) and each $(f_m)_{R'/I'^m}$ (resp.,~$(f_m)_{R''/I''^m}$) is the zero homomorphism (resp.,~an isogeny). Thus, $f_{R'}$ must be the zero homomorphism, and we are reduced to the case when all the $f_m$ are isogenies.

Let $K_m \subset E_{R/I^m}$ be the kernel of the isogeny $f_m$. The group law of $K_m$ is the restriction of the action morphism 
\[
\qq K_m \times E_{R/I^m} \ra E_{R/I^m},
\]
so \cite{EGAIII1}*{5.1.4 and 5.4.1} supply a finite locally free $R$-subgroup $K \subset E^\sm$ that algebraizes all the $K_m$. \Cref{isog-prop}~\ref{IP-b} and the settled case when the $f_m$ are isomorphisms then provide the identification $E/K \cong E'$, so $f$ is identified with the isogeny $E \ra E/K$ and hence is a homomorphism.
\qedhere
\eenum
\epf


\section{Modular descriptions of modular curves} \lab{moduli-problems}

With the compactifications $\EEl_n$ at our disposal, we are ready to exhibit the moduli interpretations and the regularity of several classical modular curves, such as $\sX(n)$ or $\sX_1(n)$ (see \S\ref{Con-Con} for an overview of our method and of previous work).  We begin in section \ref{section-arb-level} by reviewing the construction and the properties of modular curves of arbitrary congruence level. The moduli interpretations of $\sX(n)$ and $\sX_1(n)$ given in sections \ref{Gamma-n-case} and \ref{Gamma-1-case} use Drinfeld structures on generalized elliptic curves, so in section \ref{KM-fest} we extend a number of properties of such structures from the elliptic curve case studied by Katz and Mazur. In section \ref{axiomatic}, we synthesize the arguments used for $\sX(n)$ and $\sX_1(n)$ in the form of an axiomatic result, which we use in section \ref{Gamma-1-Nn} to treat further modular curves $\wt{\sX}_1(n; n')$, $\sX_1(n; n')$, and $\sX_0(n; n')$ for suitable $n$ and $n'$. The analysis of $\sX_1(n; n')$ is used in section \ref{Hecke} to give a modular construction of some Hecke correspondences for $\sX_1(n)$.


\centering \subsection{Modular curves of congruence level \nopunct} \lab{section-arb-level} \hfill

\justify

The main goal of this section is to review the definition given by Deligne and Rapoport in \cite{DR73}*{IV.3.3} of (stacky) modular curves over $\bZ$ of congruence level. The definition is via a normalization procedure, and for general levels there is no known description of these $\bZ$-curves as moduli spaces of generalized elliptic curves equipped with additional structure (one of the principal goals of this paper is to give such a description in the case of $\Gamma_0(n)$ level). The normalization procedure rests on the case of ``no level,'' with which we begin.

\bpp[The case of no additional level] \lab{no-level}
In this case, the modular curve in question is the $\bZ$-stack $\EEl_1$ that parametrizes generalized elliptic curves with integral geometric fibers (see \Cref{def-ell-n}). In the context of level structures, we will denote $\EEl_1$ by $\sX_{\GL_2(\wh{\bZ})}$, by $\sX_{\Gamma(1)}$, or simply by $\sX(1)$, and we will denote its elliptic curve locus $\Ell$ by similar notation with $\sX$ replaced by $\sY$, e.g.,~by 
\[
\sY(1) \subset \sX(1).
\]
By  \Cref{elln-props}~\ref{EN-a}--\ref{EN-b} (i.e.,~by \cite{DR73}*{III.2.5~(i), III.1.2~(iii), and IV.2.2}), the stack $\sX(1)$ is Deligne--Mumford and the morphism $\sX(1) \ra \Spec \bZ$ is proper and smooth of relative dimension $1$. 
\epp

\bpp[The case of an arbitrary congruence level $H$] \lab{genl-level}
The \emph{level} is an open (and hence finite index) subgroup $H$ of $\GL_2(\wh{\bZ})$. Its associated modular curve $\sX_H$ is a Deligne--Mumford $\bZ$-stack that, loosely speaking, compactifies the stack $\sY_H[\f{1}{\mathrm{level}}]$ which represents the ``level $H$ moduli problem'' on elliptic curves over schemes on which bad primes that depend on the level are invertible. More precisely, given $H$, one fixes an $n \in \bZ_{\ge 1}$ for which 
\[
\Ker(\GL_2(\wh{\bZ}) \surjects \GL_2(\bZ/n\bZ)) \subset H \qq \text{and sets} \qq \ov{H} \ce \im(H \ra \GL_2(\bZ/n\bZ)).
\]
One then lets $\sY_H[\f{1}{n}]$ be the $\bZ[\f{1}{n}]$-stack that, for variable $\bZ[\f{1}{n}]$-schemes $S$, parametrizes elliptic curves $E \ra S$ equipped with an $S$-point of the finite \'{e}tale $S$-scheme 
\[
\ov{H}\setminus\Isom(E[n], (\bZ/n\bZ)^2).
\]
Finally, one defines $\sX_H$ to be the Deligne--Mumford $\sX(1)$-stack obtained by 
\[
\tst \text{normalizing} \q \sX(1) \q \text{with respect to the ``forgetful'' finite \'{e}tale morphism} \q \sY_H[\f{1}{n}] \ra \sY(1)_{\bZ[\f{1}{n}]}.
\]
One lets $\sY_H$ be the preimage of $\sY(1)$ in $\sX_H$. It is proved in \cite{DR73}*{IV.3.6} that different choices of $n$ lead to canonically isomorphic $\sX_H$.

The map 
\be \lab{XH-str}
\sX_H \ra \sX(1)
\ee
is representable, finite, and also surjective because $\sX(1)$ is irreducible. Moreover, by \cite{EGAIV2}*{6.1.5} (which applies because of ``going down'' and the normality of $\sX_H$), the map \eqref{XH-str} is flat, so it is locally free of rank $[\GL_2(\wh{\bZ}) : H]$ and  $\sX_H$ is of relative dimension $1$ over $\bZ$ at every point. By \cite{DR73}*{IV.6.7}, the proper and flat structure morphism $\sX_H \ra \Spec \bZ$ is even smooth over $\bZ[\f{1}{n}]$. If $H' \subset H$, then the finite \'{e}tale $\sY(1)$-morphism $\sY_{H'}[\f{1}{n}] \ra \sY_{H}[\f{1}{n}]$ obtained from the $S$-morphisms
\[
\ov{H'}\setminus\Isom(E[n], (\bZ/n\bZ)^2) \ra \ov{H}\setminus\Isom(E[n], (\bZ/n\bZ)^2)
\]
gives rise to the finite $\sX(1)$-morphism 
\[
\sX_{H'} \ra \sX_H.
\]
Thus, due to the following lemma and \Cref{Xn-sch}, all the $\sX_H$ are schemes for small enough $H$.
\epp

\blem \lab{Xh-sch}
If the modular curve $\sX_H$ has an open substack $U \subset \sX_H$ whose geometric points have no nontrivial automorphisms, then $U$ is a scheme that is quasi-projective over $\Spec \bZ$.
\elem

\bpf
By \Cref{rep-crit}~\ref{RC-a}, $U$ is an algebraic space. Moreover, the coarse moduli space morphism $\sX(1) \ra \bP^1_\bZ$ is separated and quasi-finite, so $U \ra \bP^1_\bZ$ is also separated and quasi-finite, and hence $U$ is a scheme by \Cref{rep-schemes}. Finally, the morphism $U \ra \bP^1_\bZ$ is quasi-projective by \cite{EGAIV3}*{8.11.2} or by Zariski's main theorem \cite{EGAIV3}*{8.12.6}, so $U \ra \Spec \bZ$ is also quasi-projective.
\epf

\brem
Due to \Cref{Xh-sch} and \cite{Con07a}*{2.2.5 (2)}, each $\sX_H$ has a unique largest open subscheme. This subscheme contains exactly those geometric points of $\sX_H$ whose automorphism functors are trivial.
\erem

One suspects that $\sX_H$ is the ``correct'' modular curve of level $H$, in part because there is no other choice granted that one believes that such a modular curve should be representable and finite over $\sX(1)$, normal, and agree with $\sY_H[\f{1}{n}]$ over $\sY(1)_{\bZ[\f{1}{n}]}$. One of the bottlenecks limiting practical usefulness of the stacks $\sX_H$ is the lack of descriptions of their functors of points (without inverting the level) in terms of generalized elliptic curves equipped with additional data. In the cases where such descriptions have been found, one has been able to analyze $\sX_H$ more thoroughly, e.g.,~to prove that $\sX_H$ is regular (and not just normal). Such regularity is useful in practice (but is not known in general)---for instance, through \cite{EGAIV2}*{6.1.5} it would ensure
flatness of the maps $\sX_H \ra \sX_{H'}$ mentioned above. Similarly, the proof of the $\bZ[\f{1}{n}]$-smoothness of $(\sX_H)_{\bZ[\f{1}{n}]}$ given in \cite{DR73}*{IV.6.7} rests on the modular description of $(\sX_H)_{\bZ[\f{1}{n}]}$  presented in \emph{loc.~cit.}~for \emph{any} $H$ (however, this description is not explicit enough to \emph{a priori} recover the ``obvious'' candidate descriptions for classical choices of $H$).

Modular descriptions of $\sX_H$ are known for most ``classical'' $H$, and we will reprove some of them in sections \ref{Gamma-n-case}--\ref{Gamma-1-Nn} below.


\centering \subsection{Drinfeld level structures on generalized elliptic curves via congruences \nopunct} \lab{KM-fest} \hfill
\justify

In order to efficiently handle all residue characteristics, the modular descriptions of various $\sX_H$ that will be discussed in subsequent sections will use Drinfeld level structures on generalized elliptic curves. In the elliptic curve case, the necessary properties of such structures follow from the work of Katz and Mazur presented in \cite{KM85}, and the goal of this section is to extend them to the generalized elliptic curve case. Some such extensions have already been obtained in \cite{Con07a}, but our method seems simpler, more direct, and applies in a wider range of situations. The key idea is to exploit ``mod $n$ congruences'' with elliptic curves: the properties of various ``mod $n$ Drinfeld level structures'' tend to be fppf local and to depend solely on the $n$-torsion $E^\sm[n]$, so for many purposes we may first use \Cref{pump-it-up} to reduce to the case when $E^\sm[n]$ is finite locally free of rank $n^2$ and then apply the following lemma to  further reduce to the elliptic curve case.

\blem \lab{n-tor-iso}
For every $n \in \bZ_{\ge 1}$ and every generalized elliptic curve $E \ra S$ for which $n$ divides the number of irreducible components of each degenerate geometric fiber, there is an fppf cover $S' \ra S$ and an elliptic curve $E' \ra S'$ for which 
\[
E^\sm_{S'}[n] \simeq E'[n].
\] 
\elem

\bpf
We may work \'{e}tale locally on $S$, so limit arguments allow us to assume that $S$ is local and strictly Henselian. We may then also assume that the special fiber of $E$ is degenerate, so the connected-\'{e}tale sequence (together with \Cref{can-submult}) shows that $E^\sm[n]$ is an extension of $\bZ/n\bZ$ by $\mu_n$. After passage to an fppf cover this extension splits and our task reduces to showing that fppf locally on $\Spec \bZ$ there is an elliptic curve $E'$ with $E'[n] \cong \mu_n \times \bZ/n\bZ$.

Via limit arguments, it suffices to find such an $E'$ over each strict Henselization $(R, \fm)$ of $\Spec \bZ$ at every closed point. The conclusion then follows from choosing an ordinary elliptic curve over $R/\fm$, lifting its Weierstrass equation to $R$, and using the connected-\'{e}tale sequence again.
\epf

To make sense of Drinfeld level structures as alluded to above, we recall the following key definition:

\bd \lab{def-Dr-str}
For a finite abelian group $A$ and a generalized elliptic curve $E \ra S$, a \emph{Drinfeld $A$-structure on $E$} is a homomorphism $\gA\colon A \ra E^\sm(S)$ for which the relative effective Cartier divisor 
\[
\tst D_\gA \ce \sum_{a\in A} [\gA(a)] \subset E^\sm
\]
is an $S$-subgroup scheme. If this $S$-subgroup $G \subset E^\sm$ is given in advance, then we say that $\gA$ is a \emph{Drinfeld $A$-structure on $G$}. 
\ed

\brem \lab{away-level}
By \cite{KM85}*{1.5.3}, if $\#A$ is invertible on $S$, then a Drinfeld $A$-structure $\gA$ on $E$ amounts to an isomorphism induced by $\gA$ between the constant $S$-group $\underline{A}_S$ and some $S$-subgroup of $E^\sm$.
\erem

\begin{conv} \lab{main-conv}
In the sequel we will sometimes deal with Drinfeld $\bZ/nm\bZ$- or $(\bZ/nm\bZ)^2$-structures for fixed $n, m \in \bZ_{\ge 1}$ and will want to obtain $\bZ/n\bZ$- or $(\bZ/n\bZ)^2$-structures by restricting to the $n$-torsion subgroups. To make sense of this we need to choose noncanonical isomorphisms
\[
\bZ/n\bZ \simeq (\bZ/nm\bZ)[n] \qq \text{and} \qq (\bZ/n\bZ)^2 \simeq (\bZ/nm\bZ)^2[n].
\]
The particular choices will never matter for the results, but for definiteness we always choose the isomorphisms induced by multiplication by $m$ on $\bZ$ or on $\bZ^2$.
\end{conv}

In the results below, the ``compare with'' references point to the elliptic curve cases treated by Katz and Mazur. We begin by detailing the properties of restrictions to subgroups of various Drinfeld structures on generalized elliptic curves. Parts \ref{R-a} and \ref{R-c} of \Cref{restrict} have been proved in \cite{Con07a}*{2.3.2} by a different method that also eventually reduces to the elliptic curve case.

\bprop \lab{restrict}
Let $n, m \in \bZ_{\ge 1}$, and let $E \ra S$ be a generalized elliptic curve.
\benum
\item {\upshape (Compare with \cite{KM85}*{5.5.2~(1) and 5.5.7~(1)}).} \lab{R-a}
If $\gA$ is a Drinfeld $(\bZ/nm\bZ)^2$-structure on $E^\sm[nm]$, then $\gA|_{(\bZ/nm\bZ)^2[n]}$ is a Drinfeld $(\bZ/n\bZ)^2$-structure on $E^\sm[n]$ and $\gA|_{\bZ/nm\bZ \times \{ 0\}}$ is a Drinfeld $\bZ/nm\bZ$-structure on $E$. 

\item {\upshape (Compare with \cite{KM85}*{5.5.8~(1)}).} \lab{R-b}
If $\gA\colon (\bZ/nm\bZ)^2 \ra E^\sm(S)$ is a group homomorphism, every prime divisor of $m$ divides $n$, and $\gA|_{(\bZ/nm\bZ)^2[n]}$ is a Drinfeld $(\bZ/n\bZ)^2$-structure on $E^\sm[n]$, then $\gA$ is a Drinfeld $(\bZ/nm\bZ)^2$-structure on $E^\sm[nm]$ (so, in particular, the number of irreducible components of each degenerate geometric fiber of $E$ is divisible by $nm$). 

\item  \lab{R-c} {\upshape (Compare with \cite{KM85}*{5.5.7~(2)}).}
If $\gA$ is a Drinfeld $\bZ/nm\bZ$-structure on $E$, then $\gA|_{(\bZ/nm\bZ)[n]}$ is a Drinfeld $\bZ/n\bZ$-structure on $E$.

\item \lab{R-d} {\upshape (Compare with \cite{KM85}*{5.5.8~(2)}).}
If $\gA \colon \bZ/nm\bZ \ra E^\sm(S)$ is a group homomorphism, every prime divisor of $m$ divides $n$, and $\gA|_{(\bZ/nm\bZ)[n]}$ is a Drinfeld $\bZ/n\bZ$-structure on $E$, then $\gA$ is a Drinfeld $\bZ/nm\bZ$-structure on $E$.

\item {\upshape (Compare with \cite{KM85}*{5.5.2~(2)}).} \lab{R-e}
For brevity, set $N \ce nm$. If $\gA$ is a Drinfeld $(\bZ/N\bZ)^2$-structure on $E^\sm[N]$ and $G \subset E^\sm$ is the subgroup $\sum_{i \in \bZ/N\bZ \times\{ 0\}} [\gA(i)]$ supplied by {\upshape\ref{R-a}}, then 
\[
\q \gA|_{\{0\} \times \bZ/N\bZ} \colon \{0\} \times \bZ/N\bZ \ra (E/G)^\sm(S)
\]
is a Drinfeld $\bZ/N\bZ$-structure on $E^\sm[N]/G \subset (E/G)^\sm$.
\eenum
\eprop

\bpf
It suffices to work fppf locally on $S$, so we may use \Cref{pump-it-up} to reduce to the case when the number of irreducible components of each degenerate geometric fiber of $E$ is divisible by $nm$ (in parts \ref{R-a} and \ref{R-e} we are in this case at the outset). We may then apply \Cref{n-tor-iso} to assume further that there is an elliptic curve $E' \ra S$ with $E'[nm] \simeq E^\sm[nm]$. By \cite{KM85}*{1.10.6 and 1.10.11}, the properties under consideration depend solely on the $S$-group scheme $E^\sm[nm]$ equipped with the homomorphism $\gA$ and not on the embedding of $E^\sm[nm]$ into a smooth $S$-group scheme of relative dimension $1$ (such as $E^\sm$ or $E'$). Thus, the claims result from their elliptic curve cases.
\epf

Cyclic subgroups of generalized elliptic curves will be important for us, so we recall their definition.

\bd  \lab{def-cyclic}
For a generalized elliptic curve $E \ra S$, a finite locally free $S$-subgroup $G \subset E^\sm$ is \emph{cyclic of order $n$} if fppf locally on $S$ there is  a Drinfeld $\bZ/n\bZ$-structure on $G$. For a $G$ that is cyclic of order $n$, a section $g \in G(S)$ is a \emph{generator of $G$} (or \emph{generates $G$}) if the homomorphism $\gA\colon \bZ/n\bZ \ra E^\sm(S)$ defined by $\gA(1) = g$ is a Drinfeld $\bZ/n\bZ$-structure on $G$. An isogeny of constant degree $n$ between generalized elliptic curves over $S$ is \emph{cyclic} if its kernel is cyclic of order $n$.
\ed

We turn to the properties of cyclic subgroups of generalized elliptic curves. Parts \ref{CC-a}, \ref{CC-d}, and \ref{CC-e} of \Cref{cyclicity} have also been reduced to the elliptic curve case in \cite{Con07a}*{2.3.7, 2.3.8, and 2.3.5} by a different method. 

\bprop \lab{cyclicity}
Let $E \ra S$ be a generalized elliptic curve, $G \subset E^\sm$ an $S$-subgroup that is finite locally free of rank $n$ over $S$, and $G^\times \subset G$ the $S$-subsheaf of generators of $G$ {\upshape (}by \cite{KM85}*{1.6.5}, the $S$-subsheaf $G^\times$ is a closed $S$-subscheme of $G$ of finite presentation{\upshape)}.
\benum
\item \lab{CC-a} {\upshape (The Katz--Mazur cyclicity criterion; compare with \cite{KM85}*{6.1.1~(1)}).} 
The subgroup $G$ is cyclic of order $n$ if and only if $G^\times$ is finite locally free of rank $\phi(n)$ over $S$. In particular, $G$ is cyclic of order $n$ if and only if it becomes cyclic of order $n$ over an fpqc cover of $S$. If $n$ is invertible on $S$ and $G$ is cyclic of order $n$, then $G^\times \ra S$ is \'{e}tale.

\item \lab{CC-b} {\upshape (Compare with \cite{KM85}*{6.1.1~(2)}).} 
If $g \in G(S)$ is a generator of $G$, then 
\[
\qq \tst G^\times = \sum_{i \in (\bZ/n\bZ)^\times} [i \cdot g] \qq  \text{as effective Cartier divisors on $E^\sm$.}
\]

\item \lab{CC-c} {\upshape (Compare with \cite{KM85}*{6.4.1}).}
There is a finitely presented closed subscheme $T \subset S$ such that the base change $G_{S'}$ to an $S$-scheme $S'$ is cyclic if and only if $S' \ra S$ factors through $T$.

\item \lab{CC-d} {\upshape (Compare with \cite{KM85}*{6.8.7}).}
If $n$ is squarefree, then $G$ is cyclic.

\item \lab{CC-f} {\upshape (Compare with \cite{KM85}*{5.5.4~(3)}).}
If $G$ is cyclic of order $n$ and the number of irreducible components of each degenerate geometric fiber of $E \ra S$ is divisible by $n$, then the subgroup $E^\sm[n]/G$ of $E/G$ is cyclic of order $n$.

\item \lab{CC-e} {\upshape (Compare with \cite{KM85}*{6.7.2}).}
If $G$ is cyclic and $g, g' \in G(S)$ are generators of $G$, then for every positive divisor $d$ of $n$ both $\f{n}{d}\cdot g$ and $\f{n}{d}\cdot g'$ are generators of the same $S$-subgroup 
\[
\qq G_d \subset G
\]
that is cyclic of order $d$. In particular, if $G$ is cyclic, then the fppf local on $S$ subgroup of $G$ defined in this way descends to a canonical cyclic $S$-subgroup $G_d \subset G$ of order $d$.
\eenum
\eprop

\bpf
Cyclicity is an fppf local condition, so we may work fppf locally on $S$. We may therefore use \Cref{pump-it-up} and \Cref{n-tor-iso} to assume that the number of irreducible components of each degenerate geometric fiber of $E \ra S$ is divisible by $n$ and that there is an elliptic curve $E' \ra S$ such that $E^\sm[n] \simeq E'[n]$. Thus, since, by \cite{KM85}*{1.10.6 and its generalization 1.10.1}, the properties under consideration depend solely on the $S$-group scheme $E^\sm[n]$ and its subgroup $G$, the claims follow from their elliptic curve cases (in \ref{CC-a}, if $n$ is invertible on $S$, then a cyclic $G$ of order $n$ becomes isomorphic to $\bZ/n\bZ$ over an \'{e}tale cover of $S$, so that $G^\times$ becomes isomorphic to the constant subscheme $(\bZ/n\bZ)^\times \subset \bZ/n\bZ$).
\epf

\bd \lab{std-subgp-def}
For a generalized elliptic curve $E \ra S$ and a cyclic $S$-subgroup $G \subset E^\sm$ of order $n$, the $S$-subgroup $G_d$ defined in \Cref{cyclicity}~\ref{CC-e} is \emph{the standard cyclic subgroup of $G$ of order $d$}. Isogenies $f_1\colon E \ra E'$ and $f_2\colon E' \ra E''$ of constant degrees between generalized elliptic curves over $S$ are \emph{cyclic in standard order} if $\Ker(f_2 \circ f_1)$ is cyclic and $\Ker f_1$ is its standard cyclic subgroup (so that, in particular, $f_1$ and $f_2$ are both cyclic by \Cref{std-subgp}~\ref{SS-c} below).
\ed

In \Cref{std-subgp,std-fact} we extend various results of \cite{KM85}*{\S6.7} about standard cyclic subgroups and standard order factorizations of cyclic isogenies to the case of generalized elliptic curves  (\S\ref{quotients} provides a robust extension of the notion of an isogeny). Some of these extensions will be important for the analysis of $\sX_{\Gamma_0(n)}$ carried out in Chapter \ref{Gamma-0-case}.

\bprop \lab{std-subgp}
Let $E \ra S$ be a generalized elliptic curve, let $G \subset E^\sm$ be a cyclic $S$-subgroup of order $n$, let $d$ and $d'$ be positive divisors $n$, and let 
\[
G_d \subset G
\]
denote the standard cyclic subgroup of order $d$.
\benum 
\item \lab{SS-a} {\upshape (Compare with \cite{KM85}*{6.7.4}).}
If $d \mid d'$, then $G_d$ is identified with the standard cyclic subgroup of $G_{d'}$ of order $d$.

\item \lab{SS-b}
Interpreting the intersection as that of fppf subsheaves of $G$, we have
\[
\qq G_{d}\cap G_{d'} = G_{\gcd(d, d')}.
\]

\item \lab{SS-d} 
If $G$ meets precisely $m$ irreducible components of every degenerate geometric fiber of $E$, then $G_d$ meets precisely $\f{m}{\gcd(m, \f{n}{d})}$ irreducible components of every degenerate geometric fiber of $E$.

\item \lab{SS-f} {\upshape (Compare with \cite{KM85}*{6.7.5}).}
Letting $G_d^\times$ denote the $S$-scheme parametrizing the generators of $G_d$ {\upshape(}so that, by Proposition {\upshape\ref{cyclicity}~\ref{CC-a}}, $G_d^\times$ is a closed subscheme of $G_d$ and is finite locally free of rank $\phi(d)$ over $S${\upshape)}, we have
\[
\tst \qq G = \sum_{d\mid n} G_d^\times \qq \text{as effective Cartier divisors on~$E^\sm$.}
\]

\item \lab{SS-c} {\upshape (Compare with \cite{KM85}*{6.7.4}).}
The quotient 
\[
\qq G/G_d \subset (E/G_d)^\sm
\]
is a cyclic $S$-subgroup of order $\f{n}{d}$, the image of any generator of $G$ generates $G/G_d$, and if $d\mid d'$, then the standard cyclic subgroup of $G/G_d$ of order $\f{d'}{d}$ is identified with $G_{d'}/G_d$.

\item \lab{SS-e} {\upshape (Compare with \cite{KM85}*{6.7.11 (2)}).}
If $n$ and $\f{n}{d}$ have the same prime divisors, then $g \in G(S)$ generates $G$ if and only if its image generates $G/G_d$, and, in particular, $g$ generates $G$ if and only if $g + h$ generates $G$ for some {\upshape(}equivalently, for any{\upshape)} $h \in G_d(S)$.
\eenum
\eprop

\bpf
Part \ref{SS-a} follows from the definitions because we may work fppf locally to assume that $G$ has a generator. Part \ref{SS-b} follows from \ref{SS-a}: since $G_{\gcd(d, d')}$ lies inside both $G_{d}$ and $G_{d'}$, it suffices to observe that $G_{d}/G_{\gcd(d, d')}$ and $G_{d'}/G_{\gcd(d, d')}$ have coprime orders and hence intersect trivially inside $G/G_{\gcd(d, d')}$. Part \ref{SS-d} follows from the definition of $G_d$. To prove part \ref{SS-f}, we pass to an fppf cover of $S$ over which $G$ admits a generator and apply \Cref{cyclicity}~\ref{CC-b}.

For the remaining \ref{SS-c} and \ref{SS-e}, we work fppf locally on $S$ and use \Cref{pump-it-up} and \Cref{n-tor-iso} to assume that $G$ has a generator, that the number of irreducible components of each degenerate geometric fiber of $E$ is divisible by $n$, and that there is an elliptic curve $E' \ra S$ with $E^\sm[n] \simeq E'[n]$. By \cite{KM85}*{1.10.6}, the properties under consideration in \ref{SS-c} and \ref{SS-e} depend solely on the $S$-group $G$ and not on its embedding into $E^\sm$ or $E'$, so \ref{SS-c} and \ref{SS-e} follow from their elliptic curve cases.
\epf

\bprop\lab{std-fact}
Let 
\[
f_1\colon E_0 \ra E_1, \qq f_2\colon E_1 \ra E_2, \qq \text{and} \qq f \ce f_2 \circ f_1\colon E_0 \ra E_2
\]
 be isogenies of constant degrees $d_1$, $d_2$, and $d_1d_2$ between generalized elliptic curves over $S$.
\benum
\item \lab{SF-a} {\upshape (Compare with \cite{KM85}*{6.7.8}).}
If $f$ is cyclic and $\Ker f_2$ is \'{e}tale over $S$, then $f_1$ and $f_2$ are cyclic in standard order.

\item \lab{SF-b} {\upshape (Compare with \cite{KM85}*{6.7.10}).}
If $d_1$ and $d_2$ are coprime, then $f$ is cyclic if and only if both $f_1$ and $f_2$ are cyclic, in which case $f_1$ and $f_2$ are cyclic in standard order.

\item\lab{SF-c} {\upshape (Compare with \cite{KM85}*{6.7.12}).}
If $f_1$ and $f_2$ are cyclic, $d_1$ and $d_2$ have the same prime divisors, and $g \in (\Ker f)(S)$ is such that $d_2 \cdot g$ generates $\Ker f_1$ and $f_1(g)$ generates $\Ker f_2$, then $f_1$ and $f_2$ are cyclic in standard order and $g$ generates $\Ker f$.

\item\lab{SF-d} {\upshape (Compare with \cite{KM85}*{6.7.15}).}
If $\{E_{i - 1} \xra{f_i} E_i\}_{i = 3}^n$ are further isogenies of constant degrees $d_i$ between generalized elliptic curves over $S$ such that $d_1, \dotsc, d_n$ all have the same prime divisors and such that for each $1 \le i \le n - 1$ the isogenies $f_i$ and $f_{i + 1}$ are cyclic in standard order, then $\Ker(f_n \circ \dotsb \circ f_1)$ is cyclic and each $\Ker(f_i \circ \dotsb \circ f_1)$ is its standard cyclic subgroup.
\eenum
\eprop

\bpf
For notational convenience, we set $n \ce 2$ in \ref{SF-a}, \ref{SF-b}, and \ref{SF-c}. By  \Cref{isog-prop} and \cite{KM85}*{1.10.6}, the properties under consideration may be expressed in terms of the $S$-group scheme $\Ker(f_n \circ \dotsc \circ f_1)$ equipped with its $S$-subgroups $\Ker(f_i \circ \dotsc \circ f_1)$. Thus, since the claims are fppf local on $S$, \Cref{pump-it-up} and \Cref{n-tor-iso} allow us to assume that the number of irreducible components of each degenerate geometric fiber of $E_0$ is divisible by $\prod_{i = 1}^n d_i$ and that there is an elliptic curve $E' \ra S$ with 
\[
\tst E^\sm_0[\prod_{i = 1}^n d_i] \simeq E'[\prod_{i = 1}^n d_i].
\]
This reduces  to the elliptic curve cases treated by Katz--Mazur in \emph{op.~cit.}
\epf

We wish to prove in \Cref{Con07-fix}~\ref{CF-b} a generalization of the claim of \cite{Con07a}*{2.4.5} that is important for the definition of $\Gamma_1(N; n)$-structures given there. The argument given in \emph{loc.~cit.}~seems to require further input: the ``universal deformation technique'' invoked towards the end of the proof does not seem to apply directly because it is based on \cite{DR73}*{III.1.2~(iii)} that requires the number of irreducible components of the closed fiber to be prime to the residue characteristic and the $\bZ/N\bZ$-structure $P$ may interfere with this requirement.

\bprop \lab{Con07-fix}
Let $E \ra S$ be a generalized elliptic curve, and let $n, m \in \bZ_{\ge 1}$.
\benum
\item \lab{CF-a}
If $G \subset E^\sm$ and $H \subset E^\sm$ are $S$-subgroups that are cyclic of orders $n$ and $m$, respectively, and $\gA$ and $\gB$ are fppf local on $S$ Drinfeld $\bZ/n\bZ$- and $\bZ/m\bZ$-structures on $G$ and $H$, then
\[
\tst\q \sum_{\substack{i \in \bZ/n\bZ \\ j \in \bZ/m\bZ}} [\gA(i) + \gB(j)]
\]
is an effective Cartier divisor on $E^\sm$ that does not depend on the choices of $\gA$ and $\gB$ and descends to a well-defined relative effective Cartier divisor on $E^\sm$ over $S$ denoted by $[G + H]$.

\item \lab{CF-b}
Set $d \ce \gcd(n, m)$ and suppose that the number of irreducible components of each degenerate geometric fiber of $E \ra S$ is divisible by $d$. If $G \subset E^\sm$ and $H \subset E^\sm$ are $S$-subgroups that are cyclic of orders $n$ and $m$, respectively, and $[G_d + H_d] = E^\sm[d]$, then $[G + H]$ is a finite locally free $S$-subgroup scheme of $E^\sm$ of order $nm$ and killed by $\lcm(n, m)$, and any Drinfeld $\bZ/n\bZ$-structure on $G$ induces a Drinfeld $\bZ/n\bZ$-structure on $[G + H]/H \subset (E/H)^\sm$.
\eenum
\eprop

\bpf
For \ref{CF-a}, the cases when either $\gA$ or $\gB$ is fixed suffice, so one only needs to observe that translation by an $S$-point is an automorphism of the $S$-scheme $E^\sm$ and hence commutes with the formation of the sum of effective Cartier divisors---for example, the left hand side of
\[
\tst \gA(i) + H = \sum_{j \in \bZ/m\bZ} [ \gA(i) + \gB(j) ]
\]
does not depend on $\gB$.

For \ref{CF-b}, we work fppf locally on $S$ and use \Cref{pump-it-up} to assume that the number of irreducible components of each degenerate geometric fiber of $E \ra S$ is divisible by $nm$ and that there are Drinfeld $\bZ/n\bZ$- and $\bZ/m\bZ$-structures $\gA$ and $\gB$ on $G$ and $H$. We then imitate the argument of \cite{Con07a}*{top of p.~231} given in the elliptic curve case. Namely, we use \cite{KM85}*{1.7.2 and 1.10.6} to ``factor into prime powers'' to reduce to the case when $n = p^r$ and $m = p^s$ for some prime $p$ and $r \le s$ (the $r \ge s$ case of the last aspect of the claim will be argued separately in the last paragraph of this proof). We assume that $r \ge 1$ (otherwise $[G + H] = H$) and, after replacing $S$ by an fppf cover, we choose a homomorphism $\wt{\gA}\colon \bZ/p^s\bZ \ra E(S)$ with $p^{s - r}\wt{\gA}(1) = \gA(1)$. By \Cref{restrict}~\ref{R-b}, 
\[
\wt{\gA} + \gB \colon (\bZ/p^s\bZ)^2 \ra E^\sm(S)
\]
is a Drinfeld $(\bZ/p^s\bZ)^2$-structure on $E[p^s]$, so, by \Cref{restrict}~\ref{R-e}, 
\[
\wt{\gA}\colon \bZ/p^s\bZ \ra (E/H)^\sm(S)
\]
is a Drinfeld $\bZ/p^s\bZ$-structure on $E/H$. Then, by \Cref{restrict}~\ref{R-c}, 
\[
\gA\colon \bZ/p^r\bZ \ra (E/H)^\sm(S)
\]
is a Drinfeld $\bZ/p^r\bZ$-structure on a subgroup $K \subset (E/H)^\sm$. Finally, by \cite{KM85}*{1.11.3}, the scheme $[G + H]$ is  the preimage of $K$ in $E$, so is a subgroup, as desired. Moreover, $[G + H]$ is killed by $p^s$ because the quotient $[G + H]/E[p^r]$ is killed by its order, i.e.,~by $p^{s - r}$, whereas $E[p^r]$ is killed by $p^r$. By construction, $\gA$, whose particular choice is irrelevant for the argument, induces a Drinfeld $\bZ/p^r\bZ$-structure on $[G + H]/H$.

It remains to prove that any $\gA$ also induces a Drinfeld $\bZ/p^r\bZ$-structure on $[G + H]/H \subset (E/H)^\sm$ when $r \ge s$ and $s \ge 1$. For this, by \Cref{restrict}~\ref{R-e}, $\gA|_{(\bZ/p^r\bZ)[p^s]}$ induces a Drinfeld $\bZ/p^s\bZ$-structure on $E/H$, so, by \Cref{restrict}~\ref{R-d}, $\gA$ induces a Drinfeld $\bZ/p^r\bZ$-structure on some $S$-subgroup $K' \subset (E/H)^\sm$, and it remains to apply \cite{KM85}*{1.11.3} again to deduce that the preimage of $K'$ in $E$ must equal $[G + H]$.
\epf

One of the cornerstones of our approach to the study of various moduli stacks of Drinfeld $A$-structures on generalized elliptic curves is a direct reduction of many questions to the $A = (\bZ/n\bZ)^2$ case. To make reductions of this sort feasible we will need the following result:

\bprop \lab{from-top}
Let $E \ra S$ be a generalized elliptic curve, let $n, m \in \bZ_{\ge 1}$, let $S'$ be a variable $S$-scheme, and recall Convention {\upshape\ref{main-conv}}.
\benum
\item \lab{FT-a}
If the number of irreducible components of each degenerate geometric fiber of $E \ra S$ is divisible by $nm$ and $\gA$ is a Drinfeld $(\bZ/n\bZ)^2$-structure on $E^\sm[n]$, then the functor
\[
\qq S' \mapsto \{ \text{Drinfeld $(\bZ/nm\bZ)^2$-structures $\gB$ on $E^\sm_{S'}[nm]$ such that $\gB|_{(\bZ/nm\bZ)[n]} = \gA_{S'}$} \}
\]
is representable by a finite locally free $S$-scheme of rank $\f{\#\GL_2(\bZ/nm\bZ)}{\#\GL_2(\bZ/n\bZ)}$ that is \'{e}tale if $nm$ is invertible on $S$.

\item \lab{FT-b} {\upshape (Compare with \cite{KM85}*{5.5.3}).}
If $E \ra S$ is a generalized elliptic curve for which $n$ divides the number of irreducible components of each degenerate geometric fiber and $\gA$ is a Drinfeld $\bZ/n\bZ$-structure on $E$, then  the functor
\[
\qq S' \mapsto \{ \text{Drinfeld $(\bZ/n\bZ)^2$-structures $\gB$ on $E^\sm_{S'}[n]$ such that $\gB|_{\bZ/n\bZ \times \{0\}} = \gA_{S'}$} \}
\]
is representable by a finite locally free $S$-scheme of rank $n\cdot \phi(n)$.

\item \lab{FT-c} {\upshape (Compare with \cite{KM85}*{5.5.3}).}
If the number of irreducible components of each degenerate geometric fiber of $E \ra S$ is divisible by $n$ and, for some $S$-subgroup $G \subset E$,  
\[
\qq\gA\colon \bZ/n\bZ \ra E^\sm(S) \qq \text{ and } \qq \gB \colon \bZ/n\bZ \ra (E/G)^\sm(S)
\]
are Drinfeld $\bZ/n\bZ$-structures on $G$ and on $E^\sm[n]/G$, respectively, then the functor
\[
\ba
\qq S' \mapsto \{ \text{Drin}&\text{feld $(\bZ/n\bZ)^2$-structures $\gG$ on $E^{\sm}_{S'}[n]$ such that} \\
&\gA_{S'} = \gG|_{\bZ/n\bZ \times \{ 0 \}} \q \text{and} \q \gB_{S'} = \gG|_{\{ 0 \} \times \bZ/n\bZ} \colon \bZ/n\bZ \ra (E/G)^\sm(S')\}
\ea
\]
is representable by a finite locally free $S$-scheme of rank $n$. 

\item \lab{FT-d}
Set $d \ce \gcd(n, m)$ and $N \ce \lcm(n, m)$. If the number of irreducible components of each degenerate geometric fiber of $E \ra S$ is divisible by $N$ and $\gA$ and $\gB$ are, respectively, Drinfeld $\bZ/n\bZ$- and $\bZ/m\bZ$-structures on $E$ such that 
\[
\qq \gA|_{(\bZ/n\bZ)[d]} + \gB|_{(\bZ/m\bZ)[d]} \colon (\bZ/d\bZ)^2 \ra E^\sm(S)
\]
is a Drinfeld $(\bZ/d\bZ)^2$-structure on $E^\sm[d]$, then the functor 
\[
\ba
\qqq S' \mapsto \{ \text{Drinfeld $(\bZ/N\bZ)^2$-struc}&\text{tures $\gG$ on $E^\sm_{S'}[N]$ such that} \\
&\gA_{S'} = \gG|_{(\bZ/N\bZ \times \{ 0 \})[n]} \q \text{and} \q \gB_{S'} = \gG|_{(\{0\} \times \bZ/N\bZ)[m]} \}
\ea
\]
is representable by a finite locally free $S$-scheme of rank $\f{N\cdot \phi(N)}{d \cdot \phi(d)}$.
\eenum
\eprop

\bpf
All the functors in question are fppf sheaves, so we may work fppf locally on $S$. Setting $N \ce nm$ (resp.,~$N \ce n$) in part \ref{FT-a} (resp.,~in parts \ref{FT-b} and \ref{FT-c}) for notational convenience, we may therefore apply \Cref{n-tor-iso} to assume that there is an elliptic curve $E' \ra S$ with 
\[
E'[N] \simeq E^\sm[N].
\]
By \cite{KM85}*{1.10.6}, all the properties and functors under consideration depend solely on the $S$-scheme $E^\sm[N]$ (and its subgroup $G$ in \ref{FT-c}), so we may pass to $E'$ to reduce to the elliptic curve case. This already settles \ref{FT-b} and \ref{FT-c}, and in order to also obtain \ref{FT-a} it remains to combine \cite{EGAIV2}*{6.1.5} with \cite{KM85}*{5.1.1}, which ensures that for every $\ell \in \bZ_{\ge 1}$, the moduli stack parametrizing Drinfeld $(\bZ/\ell\bZ)^2$-structures on elliptic curves is finite locally free of rank $\#\GL_2(\bZ/\ell\bZ)$ over $\Ell$, \'{e}tale over $\Ell_{\bZ[\f{1}{\ell}]}$, and regular. 

For the remaining elliptic curve case of \ref{FT-d}, we use \cite{KM85}*{1.7.2} to ``factor into prime powers'' and reduce to the case when 
\[
n = p^r \qq \text{and} \qq m = p^s \qq\text{for some prime $p$.}
\]
Without loss of generality $r \ge s$, so the case $s = 0$ is settled by \ref{FT-b}. In the case $s \ge 1$, by \Cref{restrict}~\ref{R-b} (i.e.,~by \cite{KM85}*{5.5.8~(1)}), the functor in question is identified with the functor parametrizing $Q \in E(S')$ such that $p^{r - s} Q = \gB_{S'}(1)$. This functor is an $E[p^{r - s}]$-torsor, so it is representable by a finite locally free $S$-scheme of rank $p^{2(r - s)} = \f{p^r\cdot \phi(p^r)}{p^s \cdot \phi(p^s)}$.
\epf

When proving the algebraicity of moduli stacks of Drinfeld structures on generalized elliptic curves we will sometimes rely on the representability of functors parametrizing various such structures on a fixed curve. The key case of this representability is \Cref{from-bottom}~\ref{FB-a} recorded below---further cases may be deduced from it with the help of \Cref{cyclicity}~\ref{CC-a}. It will be important to have such representability when the structures being parametrized are assumed to be ample, so we first review the notion of ampleness.

\bd
A finite locally free $S$-subgroup $G \subset E^\sm$ of a generalized elliptic curve $E \ra S$ is \emph{ample} if $G$ is $S$-ample as a relative effective Cartier divisor on $E$, equivalently, if $G$ meets every irreducible component of every geometric fiber of $E \ra S$. For a finite abelian group $A$, a Drinfeld $A$-structure $\gA$ on $E$ is \emph{ample} if the $S$-subgroup $D_\gA \ce \sum_{a\in A} [\gA(a)] \subset E^\sm$ is ample.
\ed

\brem \lab{ample-rem}
The role of ampleness of $\gA$ in the study of various stacks that classify Drinfeld $A$-structures on generalized elliptic curves is twofold: it facilitates descent considerations (e.g.,~the ones in the definition of a stack) by endowing $E \ra S$ with a canonical $S$-ample line bundle $\sO_E(D_\gA)$, and it also kills undesirable automorphisms that would hinder the representability of various ``forget the level'' contraction morphisms (e.g.,~if $\gA$ is ample and $S$ is a geometric point, then one sees from \Cref{auto-n-gon} that only the identity automorphism of $(E, \gA)$ fixes $(E^\sm)^0$).
\erem

\bprop \lab{from-bottom} 
Let $E \ra S$ be a generalized elliptic curve, let $S'$ be a variable $S$-scheme, and recall the notation $G_d$ and $[G + H]$ introduced in Definition~{\upshape\ref{std-subgp-def}} and Proposition~{\upshape\ref{Con07-fix}~\ref{CF-a}}.
\benum
\item \lab{FB-a} 
Fix $n, m \in \bZ_{\ge 1}$, and set $d \ce \gcd(n, m)$ and $N \ce \lcm(n, m)$. The functor 
\[
\ba
\qqq \cF\colon S' \mapsto \{ \text{cyclic $S'$-subgroups $G, H \subset E^\sm_{S'}$ of orders $n$ and $m$ with $[G_d + H_d] = E^\sm_{S'}[d]$} \}
\ea
\]
{\upshape(}resp.,~its analogue which, in addition, requires $[G + H]$ to be ample{\upshape)} is representable by a finitely presented, separated, quasi-finite, flat $S$-scheme $F$ that is \'{e}tale if $nm$ is invertible on $S$. If $N$ divides the number of irreducible components of each degenerate geometric fiber of $E \ra S$, then $F$ {\upshape(}defined without the ampleness requirement{\upshape)} is finite locally free of rank $\#\GL_2(\bZ/N\bZ) \cdot \f{d \cdot \phi(d)}{N \cdot \phi(N)\cdot \phi(n) \cdot \phi(m)} $ over $S$.

\item \lab{FB-b} {\upshape (Compare with \cite{KM85}*{6.8.1}).}
For every $n \in \bZ_{\ge 1}$, the functor 
\[
\qq \cI\colon S' \mapsto \{ \text{finite locally free $S'$-subgroups $G \subset E^\sm_{S'}$ of rank $n$} \}
\]
{\upshape(}resp.,~its analogue which, in addition, requires $G$ to be ample{\upshape)} is representable by a finitely presented, separated, quasi-finite, flat $S$-scheme $I$ that is \'{e}tale if $n$ is invertible on $S$. If $n$ divides the number of irreducible components of each degenerate geometric fiber of $E \ra S$, then $I$ {\upshape(}defined without the ampleness requirement{\upshape)} is finite locally free over $S$ and its rank is constant and equals the number of subgroups of $(\bZ/n\bZ)^2$ of order $n$.
\eenum
\eprop

\brem
In \ref{FB-a}, an important special case is $m = 1$, when $\cF$ parametrizes cyclic subgroups of order $n$. In \ref{FB-b}, due to \Cref{isog-prop}~\ref{IP-b}, $\cI$ parametrizes $n$-isogenies with source $E$.
\erem

\bpf[Proof of Proposition~{\upshape\ref{from-bottom}}]
Due to \cite{EGAIV3}*{9.6.4} and limit arguments that reduce to a Noetherian base, the additional ampleness requirement cuts out quasi-compact open subfunctors of $\cF$ and $\cI$, so the ampleness variant of the claims will follow once we establish the rest.

To ease notation, we set $N \ce n$ in \ref{FB-b}. By \cite{EGAIV4}*{18.12.12}, quasi-finite and separated morphisms are quasi-affine, so effectivity of fppf descent for relatively quasi-affine schemes enables us to work fppf locally on $S$. We may therefore apply \Cref{pump-it-up} to assume that $E^\sm$ is an open $S$-subgroup of the smooth locus of another generalized elliptic curve $E' \ra S$ for which $N$ divides the number of irreducible components of each degenerate geometric fiber. The functor $\cF$ (resp.,~$\cI$) is an open subfunctor of the corresponding functor $\cF'$ (resp.,~$\cI'$) for $E'$, and the open immersion $\cF \subset \cF'$ (resp.,~$\cI \subset \cI'$) is quasi-compact due to limit arguments, so it suffices to settle the claims for $E'$ in place of $E$. We may then use \Cref{n-tor-iso} to assume that there is an elliptic curve $E'' \ra S$ with 
\[
E''[N] \simeq E^{\prime \sm}[N].
\]
Since $E'$ and $E''$ give isomorphic functors $\cI$, this reduces \ref{FB-b} to its elliptic curve case \cite{KM85}*{6.8.1}. 

For \ref{FB-a}, we let $\cF_N'$ denote the functor that parametrizes Drinfeld $(\bZ/N\bZ)^2$-structures $\gA$ on $E^{\prime \sm}_{S'}[N]$. By \Cref{from-top}~\ref{FT-a}, $\cF'_N$ is representable by a finite locally free $S$-scheme of rank $\#\GL_2(\bZ/N\bZ)$ that is \'{e}tale if $N$ is invertible on $S$. By \Cref{restrict}~\ref{R-a} and \ref{R-c}, there is a well-defined morphism 
\[
\cF'_N \ra \cF'
\]
that sends $\gA$ to the pair of subgroups on which $\gA|_{(\bZ/N\bZ \times \{ 0\})[n]}$ and $\gA|_{(\{ 0\} \times \bZ/N\bZ)[m]}$ are Drinfeld $\bZ/n\bZ$- and $\bZ/m\bZ$-structures, respectively. By \Cref{cyclicity}~\ref{CC-a} and \Cref{from-top}~\ref{FT-d}, $\cF'_N \ra \cF'$ is representable by schemes and finite locally free of rank $\f{N \cdot \phi(N) \cdot \phi(n) \cdot \phi(m)}{d \cdot \phi(d)}$. Therefore, the desired claim about $\cF'$ follows from \cite{SGA3Inew}*{V,~4.1} (combined with \cite{EGAIV2}*{2.2.11~(ii)} and \cite{EGAIV4}*{17.7.5 and 17.7.7}). 
\epf


\centering \subsection{A modular description of $\sX_{\Gamma(n)}$ \nopunct} \lab{Gamma-n-case} \hfill

\justify

The main goal of this section is to give a modular description of $\sX_{\Gamma(n)}$, where $n \in \bZ_{\ge 1}$ and
\[
\Gamma(n) \ce \Ker(\GL_2(\wh{\bZ}) \surjects \GL_2(\bZ/n\bZ))
\]
(see \S\ref{genl-level} for the definition of $\sX_{\Gamma(n)}$; see also \S\ref{conv}). This description and the proof of its correctness follow already from the results of \cite{Con07a}, which also show the regularity and other properties of $\sX_{\Gamma(n)}$. We reprove both the description and some of the properties of $\sX_{\Gamma(n)}$ by exploiting a direct relationship with the compactification $\EEl_n$ studied in Chapter \ref{compactify-ell}. The resulting proofs seem more direct and more versatile---for instance, we will see in \S\ref{Gamma-1-case} that virtually the same strategy also handles the $H = \Gamma_1(n)$ case, which is significantly more complex for the methods of \emph{op.~cit.} Another pleasant feature of this approach is that it eliminates the crutch of analytic uniformizations---for instance, in the proof of the ``ampleness'' of $\sX(n)^\infty \subset \sX(n)$ given in \Cref{Xn-Elln}~\ref{XE-b}, the only input that is needed from the theory over $\bC$ is the fact that the coarse moduli space of $(\EEl_1)_{\bC}$ is $\bP^1_\bC$ (this comes in through our reliance on \cite{DR73}*{VI.1.1} in the proof of \Cref{Elln-coarse}).

We begin by giving the definition of the modular stack $\sX(n)$ that classifies generalized elliptic curves endowed with an ample level $n$ structure, and proceed to establish enough of its properties to arrive at the identification $\sX(n) = \sX_{\Gamma(n)}$.

\bpp[The stack $\sX(n)$] \lab{def-Xn}
This is the $\bZ$-stack that, for a fixed $n \in \bZ_{\ge 1}$, and for variable schemes $S$, parametrizes the pairs
\[
(E \xra{\pi} S,\ \gA\colon (\bZ/n\bZ)^2 \ra E^\sm(S))
\]
consisting of a generalized elliptic curve $E \xra{\pi} S$ whose degenerate geometric fibers are $n$-gons and an (automatically ample) Drinfeld $(\bZ/n\bZ)^2$-structure $\gA$ on $E^\sm[n]$. The notation agrees with that of \S\ref{no-level} because $\sX(1) = \EEl_1$. We let 
\[
\sX(n)^\infty \subset \sX(n) \qq \text{and} \qq \sY(n) \subset \sX(n)
\]
be the closed substack cut out by the degeneracy loci $S^{\infty, \pi}$ and its open complement (the elliptic curve locus), respectively. Due to \Cref{away-level}, for variable $\bZ[\f{1}{n}]$-schemes $S$, the base change $\sY(n)_{\bZ[\f{1}{n}]}$ parametrizes elliptic curves $E \ra S$ equipped with an $S$-isomorphism $\gA\colon \underline{(\bZ/n\bZ)}^2_S \isomto E[n]$.
\epp

The results of section \ref{KM-fest} lead to the following direct relationship between $\sX(n)$ and $\EEl_n$.

\bprop \lab{Xn-Elln} 
Consider the $\bZ$-morphism $f\colon \sX(n) \ra \EEl_n$ that forgets $\gA$.
\benum
\item \lab{XE-a}
The morphism $f$ is representable, finite, and locally free of degree equal to $\#\GL_2(\bZ/n\bZ)${\upshape;} moreover, $f$ is \'{e}tale over $\bZ[\f{1}{n}]$. In particular, $\sX(n)$ is a Cohen--Macaulay, reduced algebraic $\bZ$-stack that is proper, flat, and of relative dimension $1$ over $\Spec \bZ$ at every point{\upshape;} moreover, $\sX(n)$ is smooth over $\bZ[\f{1}{n}]$. 

\item \lab{XE-b}
The closed substack $\sX(n)^\infty \subset \sX(n)$ is the preimage of the closed substack $\EEl_n^\infty \subset \EEl_n$ and is a reduced relative effective Cartier divisor over $\Spec \bZ$ that meets every irreducible component of every geometric fiber of $\sX(n) \ra \Spec \bZ$ and is smooth over $\bZ[\f{1}{n}]$.
\eenum
\eprop

\bpf \hfill
\benum
\item
The asserted properties of $f$ follow from \Cref{from-top}~\ref{FT-a}. The asserted properties of $\sX(n)$ other than the reducedness then result from \Cref{elln-props}~\ref{EN-a} (and \cite{EGAIV2}*{6.4.2} for the Cohen--Macaulay aspect). By \cite{EGAIV2}*{5.8.5}, the reducedness amounts to the combination of (R$_0$) and (S$_1$). The Cohen--Macaulay aspect implies (S$_1$), whereas (R$_0$) follows from the $\bZ$-flatness and $\bZ[\f{1}{n}]$-smoothness.

\item
In the given moduli interpretation, the map $\sX(n) \ra \EEl_n$ does not change the underlying generalized elliptic curves, so an $S$-point of $\sX(n)$ factors through $\sX(n)^\infty$ if and only if its image in $\EEl_n$ factors through $\EEl_n^\infty$. In other words, 
\[
\qq \sX(n)^\infty = \sX(n) \times_{\EEl_n} \EEl_n^\infty,
\]
as desired. 
All the remaining claims then follow from \ref{XE-a} and from their counterparts for $\EEl_n$ supplied by \Cref{elln-props}~\ref{EN-c}--\ref{EN-d} and \Cref{Elln-coarse} (for the reducedness of $\sX(n)^\infty$ one uses the (R$_0$)$+$(S$_1$) criterion as in the proof of \ref{XE-a}). 
\qedhere
\eenum
\epf

\bpp[The contraction morphisms] \lab{contr-Xn}
Due to \Cref{restrict}~\ref{R-a}, the contraction morphism
\[
\q \sX(nm) \xra{c} \sX(n) \qq \text{is well defined by} \qq (E, \gA) \mapsto (c_{E^\sm[n]}(E), \gA|_{(\bZ/nm\bZ)^2[n]})
\]
(see \Cref{main-conv}) for every $n, m \in \bZ_{\ge 1}$. This morphism is compatible with its analogue for $\EEl_n$ discussed in \S\ref{contr} in the sense that there is the commutative diagram
\[
\xymatrix{
\sX(nm) \ar[d]_-{c} \ar[r]^{f_{nm}} & \EEl_{nm} \ar[d] \\
\sX(n) \ar[r]^-{f_n} & \EEl_n
}
\]
whose horizontal maps forget the level structures $\gA$.
\epp

\bprop \lab{props-c-Xn}
For every $n, m \in \bZ_{\ge 1}$, the contraction $c\colon \sX(nm) \ra \sX(n)$ is representable, finite, and locally free of rank $\f{\#\GL_2(\bZ/nm\bZ)}{\#\GL_2(\bZ/n\bZ)}$. In particular, each $\sX(n)$ is Deligne--Mumford.
\eprop

\bpf 
Since $\sX(1)$ is Deligne--Mumford, the last assertion follows from the rest (applied with $n = 1$). The representability of $c$ by algebraic spaces follows from \Cref{rep-crit}~\ref{RC-b} and \Cref{auto-n-gon}.

The contraction $c$ inherits properness and finite presentation from $\sX(nm) \ra \Spec \bZ$, and so is quasi-finite due to its moduli interpretation. Therefore, by \Cref{rep-schemes}, the map $c$ is representable by schemes and finite. It remains to prove that $c$ is flat---once this is done, the asserted rank may be read off on the elliptic curve locus by using \Cref{Xn-Elln}~\ref{XE-a}. 

The flatness of the base change 
\[
\EEl_{nm} \times_{\EEl_n} \sX(n) \xra{a} \sX(n)
\]
follows from that of $\EEl_{nm} \ra \EEl_n$ supplied by \Cref{Bnm-input}~\ref{Bnm-a}. On the other hand, 
\[
\EEl_{nm} \times_{\EEl_n} \sX(n)
\]
parametrizes generalized elliptic curves endowed with a Drinfeld $(\bZ/n\bZ)^2$-structure on $E^\sm[n]$ subject to the constraint that the degenerate geometric fibers are $nm$-gons, so the map 
\[
\sX(nm) \xra{b} \EEl_{nm} \times_{\EEl_n} \sX(n)
\]
is flat by \Cref{from-top}~\ref{FT-a}. In conclusion, the composite $c = a \circ b$ is also flat.
\epf

We are ready for the promised identification $\sX(n) = \sX_{\Gamma(n)}$.

\bthm \lab{Xn-agree} 
The Deligne--Mumford stack $\sX(n)$ is regular and is identified with the stack $\sX_{\Gamma(n)}$ of {\upshape\S\ref{genl-level}} {\upshape(}see the proof for the description of the identification\upshape{)}.
\ethm

\bpf
By \cite{KM85}*{5.1.1}, the open substack $\sY(n) \subset \sX(n)$ is regular. By combining this with the conclusions of \Cref{Xn-Elln}, we see that $\sX(n)$ satisfies both (R$_1$) and (S$_2$), i.e., is normal. Therefore, due to the conclusions of \Cref{props-c-Xn}, $\sX(n)$ is identified with the normalization of $\sX(1)$ in $\sY(n)_{\bZ[\f{1}{n}]}$. However, the moduli interpretations of the $\sY(1)$-stacks $\sY(n)_{\bZ[\f{1}{n}]}$ and $\sY_{\Gamma(n)}[\f{1}{n}]$ coincide (see \S\ref{genl-level} and \S\ref{def-Xn}), so $\sX(n)$ is identified with the normalization of $\sX(1)$ in $\sY_{\Gamma(n)}[\f{1}{n}]$, i.e.,~with $\sX_{\Gamma(n)}$. To then extend the regularity of $\sY(n)$ supplied by \cite{KM85}*{5.1.1} to the regularity of the entire $\sX(n)$, we recall that it follows from from \cite{DR73}*{4.13} that $\sX_{\Gamma(n)}$ is regular  away from the supersingular points in characteristics dividing $n$. 
\epf

In the sequel we will identify $\sX(n)$ and $\sX_{\Gamma(n)}$. We conclude the section by recording all the cases in which $\sX(n)$ is a scheme (see \cite{DR73}*{IV.2.9} for such a result over $\bZ[\f{1}{n}]$). 

\bprop \lab{Xn-sch}
The stack $\sX(n)$ is a {\upshape(}necessarily projective{\upshape)} scheme over $\bZ$ unless $n = p^s$ or $n = 2p^s$ for some prime $p$ and some $s \in \bZ_{\ge 1}$.
\eprop

\bpf
If $n = p^s$ or $n = 2p^s$, then every supersingular elliptic curve $E$ over $\ov{\bF}_p$ equipped with a Drinfeld $(\bZ/n\bZ)^2$-structure on $E[n]$ has multiplication by $-1$ as an automorphism, so $\sX(n)$ cannot be a scheme. Outside of these cases, $n = n' n''$ for relatively prime $n' \ge 3$ and $n'' \ge  3$, so, due to \cite{KM85}*{2.7.2~(1)} and \Cref{auto-n-gon}, the geometric points of $\sX(n)$ have no nontrivial automorphisms, and hence $\sX(n)$ is a projective $\bZ$-scheme by \Cref{Xh-sch}.
\epf


\centering \subsection{A modular description of $\sX_{\Gamma_1(n)}$ \nopunct} \lab{Gamma-1-case} \hfill

\justify

The main goal of this section is to give a modular description of $\sX_{\Gamma_1(n)}$, where $n \in \bZ_{\ge 1}$ and
\[
\Gamma_1(n) \ce \left\{ \p{ \begin{smallmatrix} a & b \\ c & d \end{smallmatrix} } \in \GL_2(\wh{\bZ}) \q \text{such that} \q a \equiv 1 \bmod n \ \ \text{and}\ \  c \equiv 0 \bmod n \right\}
\]
(see \S\ref{genl-level} for the definition of $\sX_{\Gamma_1(n)}$; see also \S\ref{conv}). The overall strategy is similar to the case of $\Gamma(n)$ treated in the previous section: through relations with the compactifications $\EEl_m$ we infer enough
properties of the stack $\sX_1(n)$ that classifies generalized elliptic curves endowed with an ample Drinfeld $\bZ/n\bZ$-structure to arrive at the identification $\sX_1(n) = \sX_{\Gamma_1(n)}$. As in the case of $\Gamma(n)$, this identification and the finer properties of $\sX_1(n)$, such as regularity, follow already from the results of \cite{Con07a}, but the alternative proofs given below seem simpler. In particular, when proving the regularity of $\sX_1(n)$ we do not use any computations with schemes of $\Gamma_1(n)$-structures on Tate curves or with universal deformation rings, but instead directly deduce such regularity from the regularity of $\sX(n)$.

\bpp[The stack $\sX_1(n)$] \lab{X1n-def}
This is the $\bZ$-stack that, for a fixed $n \in \bZ_{\ge 1}$ and for variable schemes $S$, parametrizes the pairs
\[
(E \xra{\pi} S,\, \gA\colon \bZ/n\bZ \ra E^\sm(S))
\]
consisting of a generalized elliptic curve $E \xra{\pi} S$ and an ample Drinfeld $\bZ/n\bZ$-structure $\gA$ on $E$. As before, we let 
\[
\sX_1(n)^\infty \subset \sX_1(n) \qq \text{and} \qq \sY_1(n) \subset \sX_1(n)
\]
be the closed substack cut out by the degeneracy loci $S^{\infty, \pi}$ and its open complement (the elliptic curve locus), respectively. 

For a positive divisor $m$ of $n$, we let 
\[
\sX_1(n)_{(m)} \subset \sX_1(n)
\]
be the open substack that classifies those $(E, \gA)$ for which the degenerate geometric fibers of $E \ra S$ are $m$-gons (the openness follows from Remark \ref{gen-loc-proj}), and we set 
\[
\sX_1(n)_{(m)}^\infty \ce  \sX_1(n)_{(m)} \cap \sX_1(n)^\infty.
\]
When $m$ varies, the open substacks $\sX_1(n)_{(m)}$ cover $\sX_1(n)$, and we will use them to prove the algebraicity of $\sX_1(n)$.
\epp

\bprop \lab{X1-alg}
Consider the $\bZ$-morphism $f_{(m)}\colon \sX_1(n)_{(m)} \ra \EEl_m$ that forgets $\gA$.
\benum
\item \lab{X1-a}
The morphism $f_{(m)}$ is representable by schemes, quasi-finite, separated, flat, and of finite presentation{\upshape;} moreover, $f_{(m)}$ is \'{e}tale over $\bZ[\f{1}{n}]$. In particular, $\sX_1(n)$ is an algebraic $\bZ$-stack with a quasi-compact and separated diagonal and is flat, of finite presentation, and of relative dimension $1$ over $\Spec \bZ$ at every point{\upshape;} moreover, $\sX_1(n)$ is smooth over $\bZ[\f{1}{n}]$.

\item \lab{X1-b}
The closed substack $\sX_1(n)_{(m)}^\infty \subset \sX_1(n)_{(m)}$ is the preimage of $\EEl_m^\infty \subset \EEl_m$. In particular, $\sX_1(n)^\infty \subset \sX_1(n)$ is a reduced relative effective Cartier divisor over $\Spec \bZ$ that is smooth over $\bZ[\f{1}{n}]$.
\eenum
\eprop

\bpf \hfill
\benum
\item
The asserted properties of $f_{(m)}$ follow from \Cref{from-bottom}~\ref{FB-a} and \Cref{cyclicity}~\ref{CC-a}. Since the $\sX_1(n)_{(m)}$ cover $\sX_1(n)$, the asserted properties of $\sX_1(n)$ follow from those of $f_{(m)}$ and from \Cref{elln-props}~\ref{EN-a}.

\item
For the first assertion, it suffices to observe that in the given moduli interpretation, the map $f_{(m)}$ does not change the underlying generalized elliptic curve. The remaining assertions then follow from the first, \ref{X1-a}, and \Cref{elln-props}~\ref{EN-c}--\ref{EN-d}, using the (R$_0$)$+$(S$_1$) criterion together with \cite{EGAIV2}*{6.4.2} to establish the claimed reducedness. \qedhere
\eenum
\epf

\bpp[The relation to $\sX(n)$]
There is a forgetful contraction morphism
\[
g\colon \sX_1(n) \ra \sX(1),
\]
and, due to \Cref{restrict}~\ref{R-a}, also an $\sX(1)$-morphism  
\[
h\colon \sX(n) \ra \sX_1(n), \qq (E, \gA) \mapsto (c_{\gA|_{\bZ/n\bZ \times \{ 0\}}}(E), \gA|_{\bZ/n\bZ \times \{ 0\}})
\]
that contracts $E$ with respect to the unique finite locally free subgroup of $E^\sm$ on which $\gA|_{\bZ/n\bZ \times \{ 0\}}$ is a Drinfeld $\bZ/n\bZ$-structure. 
\epp

We will extract further information about $\sX_1(n)$ by studying $h$. The main difficulty is that $h$ changes $E$, which makes its key properties, such as flatness, less transparent. To overcome this, we will further exploit the compactifications $\EEl_m$.

\bthm \lab{X1-grand} \hfill
\benum
\item \lab{X1G-a}
The morphism $h\colon \sX(n) \ra \sX_1(n)$ is representable, finite, and locally free of rank $n \cdot \phi(n)$. In particular, $\sX_1(n) \ra \Spec \bZ$ is proper, $\sX_1(n)$ is regular, and $\sX_1(n)^\infty$ meets every irreducible component of every geometric $\bZ$-fiber of $\sX_1(n)$.

\item \lab{X1G-b}
The contraction $g\colon \sX_1(n) \ra \sX(1)$ is representable, finite, and locally free of rank $\f{\#\GL_2(\bZ/n\bZ)}{n\cdot \phi(n)}$.

\item \lab{X1G-c}
The stack $\sX_1(n)$ is Deligne--Mumford and is identified with the stack $\sX_{\Gamma_1(n)}$ of {\upshape\S\ref{genl-level}}; more precisely, both $\sX_1(n)$ and $\sX_{\Gamma_1(n)}$ are the normalizations of $\sX(1)$ in $\sY_1(n)_{\bZ[\f{1}{n}]} \cong \sY_{\Gamma_1(n)}[\f{1}{n}]$.
\eenum
\ethm

\bpf \hfill
\benum
\item
The representability of $h$ by algebraic spaces follows from \Cref{rep-crit}~\ref{RC-b} and \Cref{auto-n-gon}. 
Let $\sX(n)_{(m)} \subset \sX(n)$ be the $h$-preimage  of $\sX_1(n)_{(m)}$, let $h_{(m)} \colon \sX(n)_{(m)} \ra \sX_1(n)_{(m)}$ be the restriction of $h$, and let $f_{(m)} \colon \sX_1(n)_{(m)} \ra \EEl_m$ be the forgetful map studied in \Cref{X1-alg}. By \eqref{G-Gpr-id}, the composition $f_{(m)}\circ h_{(m)}$ agrees with the composition
\[
\qq \sX(n)_{(m)} \ra \EEl_n \xra{c} \EEl_m
\]
in which the first map forgets the Drinfeld $(\bZ/n\bZ)^2$-structure. Therefore, the universal property of the fiber product gives the commutative diagram
\[
\xymatrix{
\sX(n)_{(m)} \ar[r]^-{h'} \ar[rd]_-{h_{(m)}} & \sX_1(n)_{(m)}  \times_{\EEl_m} \EEl_n \ar[d]^{h''}\ar[r] & \EEl_n\ar[d]_-{c} \\
& \sX_1(n)_{(m)} \ar[r]^-{f_{(m)}} & \EEl_m
}
\]
in which the square is Cartesian. By \Cref{from-top}~\ref{FT-b}, the map $h'$ is representable and finite locally free of rank $n\cdot \phi(n)$. By \Cref{Bnm-input}~\ref{Bnm-a}, the base change $h''$ of $c$ is proper, flat, and surjective. The representable map $h_{(m)}$ is therefore proper, flat, surjective, and, due to its moduli interpretation, also quasi-finite. Since $h$ inherits these properties, we see from \Cref{rep-schemes} that $h$ is representable by schemes and finite locally free. Its rank is determined  on the elliptic curve locus, so equals $n\cdot \phi(n)$.

The remaining claims follow from the combination of \Cref{Xn-Elln}, \Cref{Xn-agree}, and \cite{EGAIV2}*{6.5.3~(i)}, once we establish the $\bZ$-separatedness of $\sX_1(n)$. For this, since the diagonal map $\Delta_{\sX_1(n)/\bZ}$ is separated and of finite type by \Cref{X1-alg}~\ref{X1-a}, its properness follows from the commutative diagram
\[
\qq\xymatrix{
\sX(n) \ar[d]^{h}\ar[rr]^-{\Delta_{\sX(n)/\bZ}} && \sX(n) \times_{\bZ} \sX(n)  \ar[d]^{h\times h} \\
\sX_1(n) \ar[rr]^-{\Delta_{\sX_1(n)/\bZ}} && \sX_1(n) \times_{\bZ} \sX_1(n)
}
\]
and the properness of $(h\times h) \circ \Delta_{\sX(n)/\bZ}$.

\item
Since $\sX_1(n) \ra \Spec \bZ$ is proper, $g$ is also proper. Moreover, $g$ is representable by algebraic spaces and quasi-finite due to its moduli interpretation, \Cref{rep-crit}~\ref{RC-b}, and \Cref{auto-n-gon}. Thus, due to \Cref{rep-schemes}, $g$ is representable by schemes and finite. The remaining assertions follow by considering the composite $\sX(n) \xra{h} \sX_1(n) \xra{g} \sX(1)$ and combining \ref{X1G-a} with \Cref{props-c-Xn}.

\item
Thanks to \ref{X1G-b}, the Deligne--Mumford property is inherited from $\sX(1)$. For the rest, due to the regularity of $\sX_1(n)$ and the finiteness of $\sX_1(n) \ra \sX(1)$, we need to identify the stack $\sY_1(n)_{\bZ[\f{1}{n}]}$ with the stack $\sY_{\Gamma_1(n)}[\f{1}{n}]$ that, for variable $\bZ[\f{1}{n}]$-schemes $S$, parametrizes pairs consisting of an elliptic curve $E \ra S$ and an $S$-point of the finite \'{e}tale $S$-scheme
\[
\{ \p{\begin{smallmatrix} 1 & * \\ 0 & * \end{smallmatrix}} \subset \GL_2(\bZ/n\bZ) \} \setminus \Isom(E[n], (\bZ/n\bZ)^2).
\]
The datum of such an $S$-point amounts to the datum of an isomorphism between $\bZ/n\bZ$ and a subgroup of $E$, so the sought identification results from \Cref{away-level}.
\qedhere
\eenum
\epf


\centering \subsection{An axiomatic criterion for recognizing correctness of a modular description \nopunct} \lab{axiomatic} \hfill

\justify

The arguments of the preceding section that supplied the identification $\sX_1(n) = \sX_{\Gamma_1(n)}$ and proved the regularity of $\sX_{\Gamma_1(n)}$ illustrate a general method that will similarly handle more complicated cases in the sequel. Therefore, in order to avoid repetitiveness, we wish to present the following axiomatic result that ensures that for any open subgroup $H \subset \GL_2(\wh{\bZ})$ any ``good enough'' candidate stack $\sX_H'$ agrees with the $\sX_H$ defined in \S\ref{genl-level} and that $\sX_H$ is automatically regular whenever such a good candidate is present. Of course, the main difficulty of this approach to the regularity of $\sX_H$ lies in finding a suitable $\sX_H'$. In all the cases presented in the sequel, the candidate $\sX_H'$ will be defined by a modular description of its functor of points and \Cref{axiom} will act as a criterion for recognizing that this modular description actually yields $\sX_H$.

\bthm \lab{axiom}
Let $H \subset \GL_2(\wh{\bZ})$ be an open subgroup, let $n \in \bZ_{\ge 1}$ be such that $\Gamma(n) \subset H$, and let $\sX_H'$ be a $\bZ$-stack. 
\benum
\item \lab{axiom-a}
If there is a cover 
\[
\qq\tst \sX_H' = \bigcup_{m \mid n} (\sX_H')_{(m)}\qq \text{by open substacks} \qq (\sX_H')_{(m)}\subset \sX_H'
\] 
each of which admits a representable by algebraic spaces, separated, finite type morphism 
\[
\qq (\sX_H')_{(m)} \ra \EEl_{d(m)}
\]
for some $d(m) \in \bZ_{\ge 1}$, then $\sX_H'$ is algebraic, has a quasi-compact and separated diagonal $\Delta_{\sX_H'/\bZ}$, and is of finite type over $\bZ$.

\item \lab{axiom-b}
 If $\sX_H'$ is algebraic, has a quasi-compact and separated diagonal, is of finite type over $\bZ$, and
\benuma \addtocounter{enumii}{0}
\item \lab{axiom-2}
there is a proper, flat, and surjective $\bZ$-morphism $\sX(n) \xra{h} \sX_H'$,
\eenum
then $\sX_H'$ is regular, $\sX_H' \ra \Spec \bZ$ is a proper, flat surjection, and $(\sX_H')_{\bZ[\f{1}{n}]}$ is $\bZ[\f{1}{n}]$-smooth.

\item \lab{axiom-c}
If $\sX_H'$ is algebraic, $\bZ$-proper, and satisfies {\upshape\ref{axiom-2}} together with 
\benuma \addtocounter{enumii}{1}
\item \lab{axiom-3}
there is a representable by algebraic spaces $\bZ$-morphism $\sX_H' \xra{g} \sX(1)$ that over $\bZ[\f{1}{n}]$ is identified with the morphism $\sY_H[\f{1}{n}] \ra \sY(1)_{\bZ[\f{1}{n}]}$ of {\upshape\S\ref{genl-level}}, and

\item \lab{axiom-4}
the composition $g \circ h \colon \sX(n) \ra \sX(1)$ is identified with the contraction of {\upshape\S\ref{contr-Xn}},
\eenum
then $\sX_H'$ is Deligne--Mumford and the morphism $g$ induces the identification 
\[
\qq \sX_H = \sX_H';
\] 
more precisely, then both $\sX_H$ and $\sX_H'$ are the normalizations of $\sX(1)$ in $\sY_H[\f{1}{n}]$.
\eenum
\ethm

\brem
The flatness of $h$ is one of the most stringent requirements. For the $\sX_H'$ that we will construct this flatness will be supplied by the results of Katz and Mazur through congruences with elliptic curves (see \Cref{from-top}~\ref{FT-b} and the proof of \Cref{X1-grand}~\ref{X1G-a} for an example).
\erem

\bpf[Proof of Theorem~{\upshape\ref{axiom}}] \hfill
\benum
\item
The algebraicity of each $(\sX_H')_{(m)}$ follows from that of $\EEl_{d(m)}$ supplied by \Cref{elln-props}~\ref{EN-a} (see \cite{LMB00}*{4.5~(ii)}). This suffices for the algebraicity of $\sX_H'$ because the diagonal $\Delta_{\sX_H'/\bZ}$ factors as the composition
\[
\qq \tst \sX_H' = \bigcup_{m\mid n} (\sX_H')_{(m)} \ra \bigcup_{m\mid n} (\sX_H')_{(m)} \times_{\bZ} (\sX_H')_{(m)} \subset \sX_H' \times_\bZ \sX_H'
\]
in which the inclusion is representable by open immersions. Since the inclusion is also quasi-compact and each $(\sX_H')_{(m)}$ is separated over $\bZ$, i.e.,~each $\Delta_{(\sX_H')_{(m)}/\bZ}$ is proper, it also follows that $\Delta_{\sX_H'/\bZ}$ is quasi-compact and separated.

\item 
 In the commutative diagram 
\[
\qq\xymatrix{ 
\sX(n) \ar[rr]^-{\Delta_{\sX(n)/\bZ}} \ar[d]^{h} && \sX(n) \times_\bZ \sX(n) \ar[d]^{h \times h} \\
\sX_H' \ar[rr]^-{\Delta_{\sX_H'/\bZ}} && \sX_H' \times_\bZ \sX_H'
}
\]
the composite $(h \times h) \circ \Delta_{\sX(n)/\bZ}$ is proper, $\Delta_{\sX_H'/\bZ}$ is separated and of finite type, and $h$ is surjective, so $\Delta_{\sX_H'/\bZ}$ is proper. In other words, $\sX_H' \ra \Spec \bZ$ is separated, so $\sX_H'$ inherits $\bZ$-properness from $\sX(n)$. Due to the flatness and surjectivity of $h$, the flatness, regularity, and smoothness aspects for $\sX_H'$ follow from the corresponding aspects for $\sX(n)$ supplied by \Cref{Xn-Elln}~\ref{XE-a} and \Cref{Xn-agree}.

\item 
The Deligne--Mumford property follows from the representability of $g$. The map $g$ inherits properness from $\sX_H' \ra \Spec \bZ$ and quasi-finiteness from $g \circ h$, so $g$ is finite by \Cref{rep-schemes}. Moreover, $\sX_H'$ is normal by \ref{axiom-b}, so, due to the requirement \ref{axiom-3}, $g$ identifies $\sX_H'$ with the normalization of $\sX(1)$ with respect to $\sY_H[\f{1}{n}] \ra \sY(1)_{\bZ[\f{1}{n}]}$. On the other hand, by definition, this normalization is $\sX_H$ (see \S\ref{genl-level}). \qedhere
\eenum
\epf

\beg
\Cref{axiom} is useful for proving that ``obvious'' candidate modular descriptions for various mixtures of standard moduli problems are correct. When treating ``mixture situations,'' one cannot simply ``reduce to individual constituents'' via fiber products (unlike on the elliptic curve locus): such ``reductions'' fail already in situations where no mixtures are involved, for instance, 
\[
\sX(15) \not \cong \sX(3) \times_{\sX(1)} \sX(5), \qq \text{even though} \qq \sY(15) \cong \sY(3) \times_{\sY(1)} \sY(5),
\]
as one sees by inspecting the ramification at the cusps (e.g.,~$\bC\llb q^{\f{1}{15}}\rrb \not\cong \bC\llb q^{\f{1}{3}}\rrb \tensor_{\bC\llb q\rrb} \bC\llb q^{\f{1}{5}}\rrb$).

The concrete example of a ``mixture situation'' for which we wish to illustrate \Cref{axiom} has 
\[
H = \Gamma(d) \cap \Gamma_1(\ell) \qq \text{with coprime} \qq d, \ell \in \bZ_{\ge 1}.
\]
For this $H$, due to the factorizations of Drinfeld structures discussed in \cite{KM85}*{1.7.2}, the ``obvious'' candidate $\sX_H'$ is the stack that, for variable schemes $S$, parametrizes ample Drinfeld $((\bZ/d\bZ)^2 \times \bZ/\ell\bZ)$-structures $\gA$ on generalized elliptic curves $E \ra S$ subject to the requirement that $\gA|_{(\bZ/d\bZ)^2 \times \{ 0 \}}$ is a Drinfeld $(\bZ/d\bZ)^2$-structure on $E^\sm[d]$ (so $d$ divides the number of irreducible components of each degenerate geometric fiber of $E \ra S$). 

For this $\sX_H'$, we let the maps $h$ and $g$ in \Cref{axiom} be the forgetful contractions with $n = d\ell$ and let 
\[
(\sX_H')_{(m)} \subset \sX_H'
\]
be the open substack parametrizing those $E \ra S$ whose degenerate geometric fibers are $m$-gons. The requirements of \Cref{axiom}~\ref{axiom-a} are met due to \cite{KM85}*{1.7.2} and Propositions \ref{restrict}~\ref{R-a}, \ref{cyclicity}~\ref{CC-a}, and \ref{from-bottom}~\ref{FB-a} (with $(n, m) = (d\ell, d)$ in the latter). The requirement \ref{axiom-b}~\ref{axiom-2} is checked with the help of a diagram analogous to the one in the proof of \Cref{X1-grand}~\ref{X1G-a}, the key point being that the induced map
\[
\sX(n)_{(m)} \ra (\sX_H')_{(m)} \times_{\EEl_m} \EEl_n
\]
from the $h$-preimage $\sX(n)_{(m)}$ of $(\sX_H')_{(m)}$ is finite locally free of rank $\ell \cdot \phi(\ell)$ due to \Cref{from-top}~\ref{FT-b}. The requirement \ref{axiom-c}~\ref{axiom-3} is checked as in the proof of \Cref{X1-grand}~\ref{X1G-c} by using the fact that the image of $H$ in $\GL_2(\bZ/n\bZ)$ is the pointwise stabilizer of $(\bZ/d\bZ)^2 \times \bZ/\ell\bZ$ in $(\bZ/n\bZ)^2$. Finally, the requirement \ref{axiom-c}~\ref{axiom-4} follows from the definitions of $g$ and $h$. 

In conclusion, 
\[
\sX_H' = \sX_{\Gamma(d) \cap \Gamma_1(\ell)}
\]
and $\sX_{\Gamma(d) \cap \Gamma_1(\ell)}$ is regular (such regularity at the cusps is not an automatic consequence of the regularity of $\sX_{\Gamma(d)}$ and $\sX_{\Gamma_1(\ell)}$).
\eeg


\centering \subsection{A modular description of $\sX_{\Gamma_1(n;\, n')}$ and $\sX_{\Gamma_0(n;\, n')}$ for suitable $n$ and $n'$ \nopunct} \lab{Gamma-1-Nn} \hfill

\justify

Let $n$ and $n'$ be positive integers, and let 
\[
\Gamma_1(n; n') \subset \GL_2(\wh{\bZ})
\]
be the preimage of the subgroup of $\GL_2(\bZ/nn'\bZ)$ that stabilizes the subgroup $\{0\} \times (\bZ/nn'\bZ)[n']$ in $(\bZ/nn'\bZ)^2$ and that fixes $(\bZ/nn'\bZ)[n] \times \{0\}$ pointwise. Our goal is to prove that the ``obvious'' candidate modular description for $\sX_{\Gamma_1(n; n')}$ presented in \S\ref{X1-Nn-def} is correct under the assumption that
\[
\ord_p(n') \le \ord_p(n) + 1
\]
for every prime $p$. The importance of $\sX_{\Gamma_1(n;\, n')}$ stems from its role in defining Hecke correspondences for $\sX_1(n)$ (see section \ref{Hecke}), but there also are the following reasons for treating $H = \Gamma_1(n; n')$.
\begin{itemize}
\item
The techniques used below to study $\sX_{\Gamma_1(n;\, n')}$ simultaneously expose properties of the stack $\sX_0(n)^\naive$ that parametrizes generalized elliptic curves equipped with an ample cyclic subgroup of order $n$. Although in general $\sX_0(n)^\naive$ does not agree with $\sX_{\Gamma_0(n)}$, its properties will nevertheless be crucial for the study of $\sX_{\Gamma_0(n)}$ in Chapter \ref{Gamma-0-case}.

\item
Under the additional assumption that $\ord_p(n') \le \ord_p(n)$ for all $p \mid \gcd(n, n')$, the correctness of the candidate modular description of $\sX_{\Gamma_1(n;\, n')}$ also follows from the results of \cite{Con07a} but it seems worthwhile to simplify the proofs of \emph{op.~cit.}~with the help of the general \Cref{axiom}. In fact, Conrad does not assume that $\ord_p(n') \le 1$ for $p \nmid n$, but outside this case the forgetful contraction morphism from the algebraic stack $\cM_{\Gamma_1(n;\, n')}$ that he constructs in \emph{op.~cit.}~to $\sX(1)$ is not representable (even over $\bC$), so $\cM_{\Gamma_1(n;\, n')}$ cannot agree with $\sX_{\Gamma_1(n;\, n')}$ (a related pathology is that $\cM_{\Gamma_1(n;\, n')}$ is not Deligne--Mumford in characteristics $p\nmid n$ with $p^2 \mid n'$).
\end{itemize}

 In order to also recover and generalize the results of \cite{Con07a} in the cases when $\ord_p(n') > 1$ for some prime $p \nmid n$, we initially drop \emph{all} requirements on $n$ and $n'$, define a certain stack $\sX_1(n; n')$ that agrees with the stack $\cM_{\Gamma_1(n;\, n')}$ considered in \emph{op.~cit.}~(in the cases in which $\cM_{\Gamma_1(n;\, n')}$ was defined), prove that $\sX_1(n; n')$ is algebraic, $\bZ$-proper, and regular (among other properties), and only then impose assumptions on $n$ and $n'$ in order to arrive at the agreement with $\sX_{\Gamma_1(n;\, n')}$.

\bpp[The stack $\sX_1(n; n')$] \lab{X1-Nn-def}
This is the $\bZ$-stack that, for fixed $n, n' \in \bZ_{\ge 1}$ with $d \ce \gcd(n, n')$ and for variable schemes $S$, parametrizes the triples
\[
(E \xra{\pi} S,\, \gA\colon \bZ/n\bZ \ra E^\sm(S),\, H)
\]
consisting of a generalized elliptic curve $E \xra{\pi} S$, a Drinfeld $\bZ/n\bZ$-structure $\gA$ on some $S$-subgroup $G \subset E^\sm$, and a cyclic $S$-subgroup $H \subset E^\sm$ of order $n'$ subject to the requirements that 
\be \lab{X1Nn-req}
[G_d + H_d] = E^\sm[d] \qq \text{and} \qq [G + H] \ \ \text{is ample}
\ee
(we implicitly use \Cref{std-subgp-def} and \Cref{Con07-fix}~\ref{CF-a} to make sense of $[G_d + H_d]$ and $[G + H]$). The effectivity of descent needed for $\sX_1(n; n')$ to be a stack is ensured by the ampleness of $[G + H]$ as in \Cref{ample-rem}. The requirement $[G_d + H_d] = E^\sm[d]$ implies that the number of irreducible components of each degenerate geometric fiber of $E$ is divisible by $d$, so \Cref{Con07-fix}~\ref{CF-b} ensures that $[G + H]$ is a finite locally free $S$-subgroup of $E^\sm$  of rank $nn'$ that is killed by $\lcm(n, n')$. 

We let 
\[
\sX_1(n; n')^\infty \subset \sX_1(n; n') \qq \text{and} \qq \sY_1(n; n') \subset \sX_1(n; n')
\]
 be the closed substack cut out by the degeneracy loci $S^{\infty, \pi}$ and its open complement (the elliptic curve locus), respectively. Similarly to the case of $\sX_1(n)$ (discussed in \S\ref{X1n-def}), for every positive divisor $m$ of $\lcm(n, n')$, we let 
 \[
 \sX_1(n; n')_{(m)} \subset \sX_1(n; n')
 \]
be the open substack over which the degenerate geometric fibers of $E$ are $m$-gons.
\epp

\bpp[Variants $\wt{\sX}_1(n; n')$ and $\sX_0(n; n')$] \lab{var}
Slight modifications of the definition of $\sX_1(n; n')$ give the following related stacks:
\begin{itemize}
\item
The stack $\wt{\sX}_1(n; n')$ obtained by replacing the datum $H$  by the datum of a Drinfeld $\bZ/n'\bZ$-structure $\gB$ on some $S$-subgroup $H \subset E^\sm$ subject to \eqref{X1Nn-req};

\item
The stack $\sX_0(n; n')$ obtained by replacing the datum $\gA$  by the datum of a cyclic $S$-subgroup $G \subset E^\sm$ of order $n$ subject to \eqref{X1Nn-req}.
\end{itemize}
Due to \Cref{cyclicity}~\ref{CC-a}, the forgetful maps
\be \lab{var-forget}
\wt{\sX}_1(n; n') \ra \sX_1(n; n') \qq \text{and} \qq \sX_1(n; n') \ra \sX_0(n; n')
\ee
are representable by schemes, finite locally free of ranks $\phi(n')$ and $\phi(n)$, respectively, and, over $\bZ[\f{1}{n'}]$ and $\bZ[\f{1}{n}]$, respectively, \'{e}tale. As before, for every positive divisor $m$ of $\lcm(n, n')$ we let 
\[
\wt{\sX}_1(n; n')_{(m)} \subset \wt{\sX}_1(n; n') \qq  \text{and} \qq \sX_0(n; n')_{(m)} \subset \sX_0(n; n')
\] 
be the open substacks over which the degenerate geometric fibers of $E$ are $m$-gons,  let 
\[
\wt{\sX}_1(n; n')^\infty \subset \wt{\sX}_1(n; n') \qq \text{and} \qq \sX_0(n; n')^\infty \subset \sX_0(n; n')
\]
be the degeneracy loci, and let 
\[
\wt{\sY}_1(n; n') \subset \wt{\sX}_1(n; n') \qq \text{and} \qq \sY_0(n; n') \subset \sX_0(n; n')
\]
be the elliptic curve loci.

For suitably constrained $n$ and $n'$, the stacks $\wt{\sX}_1(n; n')$ and $\sX_0(n; n')$ were also considered in \cite{Con07a} (in the notation $\cM_{\wt{\Gamma}_1(N;\, n)}$ and $\cM_{\Gamma_0(N;\, n)}$). There $\wt{\sX}_1(n; n')$ was often used as an intermediary in the proofs of the properties of $\sX_1(n; n')$, whereas $\sX_0(n; n')$ was mentioned on page 273 in relation to modifications that one needs to make to the method of \emph{op.~cit.}~to also construct Hecke correspondences for $\sX_0(n)$. We will see below that the proofs of the properties of $\sX_1(n; n')$ will also prove the corresponding properties of $\wt{\sX}_1(n; n')$ and $\sX_0(n; n')$.
\epp

\bpp[Contraction maps from $\sX(nn')$]
There is a forgetful contraction map
\be \lab{Xtilde-for}
\sX(nn') \ra \wt{\sX}_1(n; n')
\ee
that sends a Drinfeld $(\bZ/nn'\bZ)^2$-structure $\gG$ to $\gA \ce \gG|_{(\bZ/nn'\bZ)[n] \times \{ 0\}}$ and $\gB \ce \gG|_{ \{ 0\} \times (\bZ/nn'\bZ)[n']}$ (see \Cref{restrict}~\ref{R-a} and \ref{R-c} and Convention~\ref{main-conv}) and contracts the underlying generalized elliptic curve accordingly. Similar forgetful contraction maps 
\[
\sX(nn') \ra \sX_1(n; n') \qq  \text{and} \qq \sX(nn') \ra \sX_0(n; n')
\]
 are the compositions of \eqref{Xtilde-for} with the forgetful maps from \eqref{var-forget}.
\epp

We are ready to address the basic properties of the stack $\sX_1(n; n')$ and its variants.

\bthm \lab{Con07-feast}
Fix $n, n' \in \bZ_{\ge 1}$ and let $\sX \in \{ \wt{\sX}_1(n; n'), \sX_1(n; n'), \sX_0(n; n')\}$.
\benum
\item \lab{feast-a}
The $\bZ$-stack $\sX$ is algebraic, regular, proper, flat, and of relative dimension $1$ over $\Spec \bZ$ at every point{\upshape;} moreover, $\sX$ is smooth over $\bZ[\f{1}{nn'}]$. The diagonal $\Delta_{\sX/\bZ}$ is finite.

\item \lab{feast-c}
The forgetful contraction map $\sX(nn') \ra \sX$ is representable by schemes and is finite locally free of constant positive rank.

\item \lab{feast-b}
The closed substack $\sX^\infty \subset \sX$ is a reduced relative effective Cartier divisor over $\Spec \bZ$ that meets every irreducible component of every geometric $\bZ$-fiber of $\sX$ and is smooth over $\bZ[\f{1}{nn'}]$.
\eenum
\ethm

\bpf \hfill
\benum
\item
By \Cref{from-bottom}~\ref{FB-a} and the finiteness of the maps \eqref{var-forget}, for every positive divisor $m$ of $\lcm(n, n')$ the forgetful map $\sX_{(m)} \ra \EEl_m$ is representable, separated, and of finite type, so, by \Cref{axiom}~\ref{axiom-a}, $\sX$ is algebraic and has a quasi-compact and separated diagonal.

Except for the relative dimension and the diagonal aspects, the rest of the claim follows from \Cref{axiom}~\ref{axiom-b} once we prove that the forgetful contraction $\sX(nn') \ra \wt{\sX}_1(n; n')$ is proper, flat, and surjective. For this, we first let $\sX(nn')_{(m)}$ for every positive divisor $m$ of $\lcm(n, n')$ be the preimage of $\wt{\sX}_1(n; n')_{(m)}$. Due to \Cref{Bnm-input}~\ref{Bnm-a}, it then suffices to note that, by \Cref{from-top}~\ref{FT-a} and \ref{FT-d}, the induced map
\[
\qq \sX(nn')_{(m)} \ra \wt{\sX}_1(n; n')_{(m)} \times_{\EEl_m} \EEl_{nn'},
\]
both components of which are forgetful, is finite locally free of constant positive rank.

The relative dimension aspect will follow from the corresponding aspect for $\sX(nn')$ once we prove that the surjective map $\sX(nn') \ra \wt{\sX}_1(n; n')$ is finite locally free. In fact, due to \Cref{rep-schemes} and the previous paragraph, representability by algebraic spaces and quasi-finiteness would suffice. The representability is inherited from $\sX(nn') \ra \sX(1)$ and the quasi-finiteness follows from the moduli interpretation.

The diagonal $\Delta_{\sX/\bZ}$ is proper due to the $\bZ$-separatedness of $\sX$ and is quasi-finite due to \Cref{elln-props}~\ref{EN-a}, so its finiteness follows from \Cref{rep-schemes}.

\item 
Due to the proof of \ref{feast-a} and the fact that the forgetful contractions \eqref{var-forget} are representable and finite locally free, only the constancy of the rank requires attention and we may focus on $\sX_0(n; n')$. Moreover, since $\sY_0(n; n')$ is dense in  $\sX_0(n; n')$, we may work on the elliptic curve locus. Therefore, since the rank of $\sY(nn') \ra \sY(1)$ is constant, the conclusion follows from \Cref{from-bottom}~\ref{FB-a} which proves that $\sY_0(n; n') \ra \sY(1)$ is finite locally free of constant positive rank.

\item
The assertion about the geometric fibers follows from the corresponding assertion for $\sX(n n')^\infty \subset \sX(nn')$ supplied by \Cref{Xn-Elln}~\ref{XE-b}, so it suffices to prove that for each positive divisor $m$ of $\lcm(n, n')$ the restriction $\sX^\infty_{(m)} \subset \sX_{(m)}$ of $\sX^\infty \subset \sX$ is a reduced relative effective Cartier divisor over $\Spec \bZ$ that is smooth over $\bZ[\f{1}{nn'}]$. To do so, it suffices to note that $\sX^\infty_{(m)}$ is the pullback of $\EEl_m^{\infty}$, to apply \Cref{elln-props}~\ref{EN-c}--\ref{EN-d} and \Cref{from-bottom}~\ref{FB-a}, to use the properties of the forgetful maps \eqref{var-forget}, and to use the (R$_0$)$+$(S$_1$) criterion for reducedness.
\qedhere
\eenum
\epf

In principle it is possible to determine the largest Deligne--Mumford open substacks of $\wt{\sX}_1(n; n')$, $\sX_1(n; n')$, and $\sX_0(n; n')$ (such open substacks make sense \emph{a priori} due to \Cref{DM-make-sense}): one needs to inspect the defining modular descriptions to determine those geometric points whose automorphism functors are not \'{e}tale. To illustrate the procedure, in \Cref{DM-fix} we exhibit large Deligne--Mumford open substacks of $\wt{\sX}_1(n; n')$, $\sX_1(n; n')$, and $\sX_0(n; n')$ (the actual Deligne--Mumford loci of $\sX_1(n; n')$ and $\sX_0(n; n')$ may be larger). For the stack $\cM_{\Gamma_1(N; n)}$ considered in \cite{Con07a}, \Cref{DM-fix}~\ref{DMF-b} improves on \cite{Con07a}*{3.1.7} by proving that the Deligne--Mumford locus includes all the cusps in characteristics $p \mid N$ (even when $p^2\mid n$).

\bprop \lab{DM-fix}
Fix $n, n' \in \bZ_{\ge 1}$ and set $d \ce \gcd(n, n')$.
\benum
\item \lab{DMF-a}
The stack $\wt{\sX}_1(n; n')$ is Deligne--Mumford. In fact, the forgetful contraction morphism
\[
\qq \wt{\sX}_1(n; n') \ra \sX(1)
\]
is representable by algebraic spaces.

\item \lab{DMF-b}
The open substack of $\sX_1(n; n')$ obtained by removing the closed substacks $\sX_1(n; n')^\infty_{\bF_p}$ for the primes $p$ with $\ord_p(n') \ge \ord_p(n) + 2$ is Deligne--Mumford. If $\ord_p(n') \le \ord_p(n) + 1$ for every prime $p$, then the forgetful contraction morphism 
\[
\qq\sX_1(n; n') \ra \sX(1)
\]
is representable by algebraic spaces.

\item \lab{DMF-c}
The open substack of $\sX_0(n; n')$ obtained by removing the closed substacks $\sX_0(n; n')^\infty_{\bF_p}$ for the primes $p$ with $\abs{\ord_p(n) - \ord_p(n')} \ge 2$ is Deligne--Mumford. If $\abs{\ord_p(n) - \ord_p(n')} \le 1$ for every prime $p$, then the forgetful contraction morphism 
\[
\qq \sX_0(n; n') \ra \sX(1)
\]
is representable by algebraic spaces.
\eenum
\eprop

\bpf 
We recall from \Cref{auto-n-gon} that the automorphism functor of the standard $m$-gon generalized elliptic curve is $\mu_m \times \bZ/2\bZ$. To test the Deligne--Mumford property of an open substack of $\wt{\sX}_1(n; n')$, $\sX_1(n; n')$, or $\sX_0(n; n')$, we will use the criterion of having unramified automorphism functors at geometric points (see \Cref{DM-make-sense}). To test the representability of contraction morphisms, we will use \Cref{rep-crit}~\ref{RC-b}. These preliminary remarks already settle part \ref{DMF-a}.
\benum \addtocounter{enumi}{1}
\item
Our task is to show that if $p$ is a prime, $E$ is the standard $m$-gon with $p \mid m$ over an algebraically closed field $\ov{k}$, and $(E, \gA, H)$ is an object of $\sX_1(n; n')(\ov{k})$ with $\ord_p(n') \le \ord_p(n) + 1$, then $\mu_p \subset \Aut(E)$ does not fix both $\gA$ and $H$. By decomposing into primary parts with the help of \cite{KM85}*{1.7.2} and by contracting away from the $p$-primary part of $[G + H]$, we loose no generality by assuming that $n$, $n'$, and $m$ are powers of $p$ and $m > 1$.

Suppose that $\mu_p$ fixes both $\gA$ and $H$. Then $\gA$ cannot be ample, so $H$ is ample, $H \cap (E^\sm)^0$  contains $\mu_p \subset (E^\sm)^0$, and $\ord_p(n') \ge 2$. Therefore, the standard cyclic subgroup $H_p \subset H$ of order $p$ is contained in $(E^\sm)^0$ and hence equals $\mu_p$. Moreover, due to the requirement $\ord_p(n') \le \ord_p(n) + 1$, we have $n > 1$, so, by \Cref{restrict}~\ref{R-a}, the requirement $[G_d + H_d] = E^\sm[d]$ implies that $[G_p + H_p] = E^\sm[p]$. The latter forces $G_p$ to project isomorphically onto the $p$-torsion subgroup of the component group of $E^\sm$, so $G$ injects into this component group. Since $H$ is ample and $H \cap (E^\sm)^0 \neq 0$, this violates the requirement $\ord_p(n') \le \ord_p(n) + 1$ unless $G$ is ample, that is, unless $\gA$ is ample, which is a contradiction.

\item
Our task is to show that if $p$ is a prime, $E$ is the standard $m$-gon with $p \mid m$ over an algebraically closed field $\ov{k}$, and $(E, G, H)$ is an object of $\sX_0(n; n')(\ov{k})$ with $\abs{\ord_p(n) - \ord_p(n')} \le 1$, then $\mu_p \subset \Aut(E)$ does not fix both $G$ and $H$. As in the proof of \ref{DMF-b}, we assume that $n$, $n'$, and $m$ are powers of $p$ and $m > 1$.

Suppose that $\mu_p$ fixes both $G$ and $H$. By the conclusion of \ref{DMF-b}, $\mu_p$ cannot fix any Drinfeld $\bZ/n\bZ$-structure (resp.,~$\bZ/n'\bZ$-structure) on $G$ (resp.,~$H$), so $G$ and $H$ must both be ample, and hence must both contain $\mu_p \subset (E^\sm)^0$. Then $G_p = H_p = \mu_p$ inside $(E^\sm)^0$, which is a contradiction to the requirement $[G_p + H_p] = E^\sm[p]$ inherited from $[G_d + H_d] = E^\sm[d]$.
\qedhere
\eenum
\epf

With \Cref{DM-fix} in hand, we are ready for identifications with suitable modular curves $\sX_H$.

\bthm \lab{X1Nn-finale}
Fix $n, n' \in \bZ_{\ge 1}$.
\benum
\item \lab{X1Nnf-a}
Let $\wt{\Gamma}_1(n; n')$ be the preimage in $\GL_2(\wh{\bZ})$ of the subgroup of $\GL_2(\bZ/nn'\bZ)$ that fixes the subgroups $(\bZ/nn'\bZ)[n] \times\{ 0\}$ and $\{ 0\} \times (\bZ/nn'\bZ)[n']$ pointwise in $(\bZ/nn'\bZ)^2$. The forgetful contraction $\wt{\sX}_1(n; n') \ra \sX(1)$ induces the identification 
\[
\qq \wt{\sX}_1(n; n') = \sX_{\wt{\Gamma}_1(n;\, n')}.
\]

\item\lab{X1Nnf-b}
Let $\Gamma_1(n; n')$ be the preimage in $\GL_2(\wh{\bZ})$ of the subgroup of $\GL_2(\bZ/nn'\bZ)$ that fixes the subgroup $(\bZ/nn'\bZ)[n] \times\{ 0\}$ pointwise and stabilizes the subgroup $\{ 0\} \times (\bZ/nn'\bZ)[n']$ in $(\bZ/nn'\bZ)^2$. If $\ord_p(n') \le \ord_p(n) + 1$ for every prime $p$, then the forgetful contraction $\sX_1(n; n') \ra \sX(1)$ induces the identification 
\[
\qq \sX_1(n; n') = \sX_{\Gamma_1(n;\, n')}.
\]

\item \lab{X1Nnf-c}
Let $\Gamma_0(n; n')$ be the preimage in $\GL_2(\wh{\bZ})$ of the subgroup of $\GL_2(\bZ/nn'\bZ)$ that stabilizes the subgroups $(\bZ/nn'\bZ)[n] \times\{ 0\}$ and $\{ 0\} \times (\bZ/nn'\bZ)[n']$ in $(\bZ/nn'\bZ)^2$. If $\abs{\ord_p(n') - \ord_p(n)} \le 1$ for every prime $p$, then the forgetful contraction $\sX_0(n; n') \ra \sX(1)$ induces the identification 
\[
\qq \sX_0(n; n') = \sX_{\Gamma_0(n;\, n')}.
\]
\eenum
\ethm

\bpf
By \Cref{DM-fix}, the imposed assumptions on $n$ and $n'$ ensure that the forgetful contraction morphisms to $\sX(1)$ are representable by algebraic spaces. Therefore, due to \Cref{Con07-feast} and \Cref{axiom}~\ref{axiom-c}, we only need to show that, for variable $\bZ[\f{1}{nn'}]$-schemes $S$, the $\sY(1)_{\bZ[\f{1}{nn'}]}$-stacks
\[
\wt{\sY}_1(n; n')_{\bZ[\f{1}{nn'}]}, \qq \sY_1(n; n')_{\bZ[\f{1}{nn'}]}, \qq \text{and} \qq \sY_0(n; n')_{\bZ[\f{1}{nn'}]}
\]
 parametrize elliptic curves $E \ra S$ equipped with an $S$-point of
\[
\ba
\ov{\wt{\Gamma}_1(n; n')} \setminus \Isom(E[nn'], (\bZ/nn'\bZ)^2), \qqq \ov{\Gamma_1(n; n')} \setminus \Isom(E[nn'], (\bZ/nn'\bZ)^2), \\
\text{and}\qq\ov{\Gamma_0(n; n')} \setminus \Isom(E[nn'], (\bZ/nn'\bZ)^2), \qqqqqqqq 
\ea
\]
respectively, where overlines denote images in $\GL_2(\bZ/nn'\bZ)$. For this, it suffices to inspect the defining modular descriptions of $\wt{\sX}_1(n; n')$, $\sX_1(n; n')$, and $\sX_0(n; n')$ and to use the definitions of $\wt{\Gamma}_1(n; n')$, $\Gamma_1(n; n')$, and $\Gamma_0(n; n')$ given in the statements of \ref{X1Nnf-a}, \ref{X1Nnf-b}, and \ref{X1Nnf-c}.
\epf


\centering \subsection{A modular construction of Hecke correspondences for $\sX_1(n)$ \nopunct} \lab{Hecke} 

\justify

We wish to explain how the results of sections \ref{quotient-section}, \ref{Gamma-1-case}, and \ref{Gamma-1-Nn} give rise to a Hecke correspondence 
\[
\pi_1, \pi_2\colon \sX_{\Gamma_1(n;\, p)}  \rightrightarrows \sX_{\Gamma_1(n)}
\]
for every $n \in \bZ_{\ge 1}$ and every squarefree $p \in \bZ_{\ge 1}$ that may or may not be coprime with $n$. 

In terms of the moduli interpretations given in \S\ref{X1n-def} and \S\ref{X1-Nn-def} and proved in \Cref{X1-grand}~\ref{X1G-c} and \Cref{X1Nn-finale}~\ref{X1Nnf-b}, the maps are given by
\[
\pi_1((E, \gA, H)) = (c_\gA(E), \gA) \qq \text{and} \qq \pi_2((E, \gA, H)) = (E/H, \gA),
\]
and are well defined due to the last aspect of \Cref{Con07-fix}~\ref{CF-b} (we let $c_\gA(E)$ denote the contraction of $E$ with respect to the unique subgroup on which $\gA$ is a Drinfeld $\bZ/n\bZ$-structure). To argue that we have exhibited a correspondence, it suffices to prove the following lemma:

\blem
The maps $\pi_1$ and $\pi_2$ are representable, finite locally free, and surjective.
\elem

\bpf
Since $\pi_1$ is the $\sX(1)$-morphism induced by the inclusion $\Gamma_1(n; p) \subset \Gamma_1(n)$, its finiteness follows from the finiteness of $\sX_H \ra \sX_{H'}$ observed in the last paragraph of \S\ref{genl-level}. By \Cref{X1-grand}~\ref{X1G-a}, $\sX_{\Gamma_1(n)}$ is regular, so the flatness of $\pi_1$ follows from \cite{EGAIV2}*{6.1.5}. The surjectivity of $\pi_1$ may be checked over $(\sY_{\Gamma_1(n)})_{\bQ}$.

For the representability of $\pi_2$, due to \Cref{rep-crit}~\ref{RC-b} and the representability of $\sX_{\Gamma_1(n;\, p)} \ra \sX(1)$, it suffices to observe that if $E$ is a generalized elliptic curve over an algebraically closed field and $H \subset E^\sm$ is a finite subgroup, then every automorphism $i$ of $E$ that stabilizes $H$ and induces the identity map on $E/H$ must fix $(E^\sm)^0$ because the endomorphism $\id_{E^\sm} - i|_{E^\sm}$ of $E^\sm$ factors through $H$. The properness of $\pi_2$ follows from the $\bZ$-properness of $\sX_{\Gamma_1(n;\, p)}$ and $\sX_{\Gamma_1(n)}$, so its quasi-finiteness may be checked on geometric fibers. Finiteness of $\pi_2$ is then supplied by \Cref{rep-schemes}, and its flatness follows from \cite{EGAIV2}*{6.1.5}. Finally, the surjectivity of $\pi_2$ may be checked over $(\sY_{\Gamma_1(n)})_{\bQ}$.
\epf

In the case when $p$ is a prime, the Hecke correspondence above has already been constructed in \cite{Con07a}*{4.4.3} by a different method: due to the lack of the theory of quotients of generalized elliptic curves by arbitrary finite locally free subgroups, \emph{loc.~cit.}~first defines $\pi_2$ by the same formula on the elliptic curve locus and then argues that the resulting map extends uniquely to the entire $\sX_{\Gamma_1(n; p)}$. The construction above seems simpler and more direct, and it also produces the map $\xi$ of \cite{Con07a}*{4.4.3}: if $e$ and $e'$ are the identity sections of $E \ra S$ and $E/H \ra S$, then there is a map
\[
(e')^* (\Omega^1_{(E/H)/S} )\ra e^* (\Omega^1_{E/S})
\]  
whose formation is compatible with base change in $S$.




\section{A modular description of $\sX_{\Gamma_0(n)}$} \lab{Gamma-0-case}

For an integer $n \in \bZ_{\ge 1}$ and the subgroup
\[
\Gamma_0(n) \ce \left\{ \p{ \begin{smallmatrix} a & b \\ c & d \end{smallmatrix} } \in \GL_2(\wh{\bZ}) \ \vert\  c \equiv 0 \bmod n \right\},
\]
the goal of this chapter is to exhibit the modular curve $\sX_{\Gamma_0(n)}$ defined via normalization (see \S\ref{genl-level}) as a moduli stack parametrizing generalized elliptic curves equipped with a ``$\Gamma_0(n)$-structure,'' which on the elliptic curve locus is the datum of a subgroup that is cyclic of order $n$ in the sense of \Cref{def-cyclic}. The proof  of the correctness of this moduli interpretation in \Cref{X0n-main} will simultaneously deduce the regularity of $\sX_{\Gamma_0(n)}$ from that of $\sY_{\Gamma_0(n)}$ proved by Katz and Mazur. We begin with a naive modular description that recovers $\sX_{\Gamma_0(n)}$ only for squarefree $n$ and then proceed to refine the naive description to a description that works for any $n$.

Throughout Chapter \ref{Gamma-0-case} we fix an integer $n \in \bZ_{\ge 1}$.

\begin{pp-tweak}[The stack $\sX_0(n)^\naive$] \lab{naive-X0n-def}
This is the $\bZ$-stack that, for variable schemes $S$, parametrizes the pairs
\[
(E \xra{\pi} S,\, G)
\]
consisting of a generalized elliptic curve $E \xra{\pi} S$ and an ample $S$-subgroup $G \subset E^\sm$ that is cyclic of order $n$ (in the sense of \Cref{def-cyclic}). We call such a $G$ a \emph{naive $\Gamma_0(n)$-structure on $E$}. 

We let 
\[
\sY_0(n)^\naive \subset \sX_0(n)^\naive
\]
 be the open substack that parametrizes those pairs for which $E$ is an elliptic curve. For each positive divisor $m$ of $n$, we let 
 \[
 \sX_0(n)^\naive_{(m)} \subset \sX_0(n)^\naive
 \]
  be the open substack that parametrizes those pairs for which the degenerate geometric fibers of $E$ are $m$-gons.

In the notation of \S\ref{var}, one has 
\[
\sX_0(n)^\naive = \sX_0(n; 1),
\]
so, by \Cref{Con07-feast}~\ref{feast-a}, the stack $\sX_0(n)^\naive$ is algebraic, proper and flat over $\Spec \bZ$, and regular with finite diagonal $\Delta_{\sX_0(n)^\naive/\bZ}$. By \Cref{Con07-feast}~\ref{feast-c} (and its proof), the morphism
\[
\sX(n) \ra \sX_0(n)^\naive
\]
that sends a Drinfeld $(\bZ/n\bZ)^2$-structure $\gA$ to the subgroup on which $\gA|_{\bZ/n\bZ \times \{0\}}$ is a Drinfeld $\bZ/n\bZ$-structure and contracts the underlying generalized elliptic curve with respect to this subgroup is finite locally free of rank $n\cdot \phi(n)^2$.

If $n$ is squarefree, then \Cref{X1Nn-finale}~\ref{X1Nnf-c} proves that the contraction 
\[
\sX_0(n)^\naive \ra \sX(1) \qq \text{is identified with the structure morphism} \qq \sX_{\Gamma_0(n)} \ra \sX_0(1).
\]
This identification fails when $n$ is is divisible by $p^2$ for some prime $p$: variants of the example given in \S\ref{non-rep} show that for such $n$ the contraction $\sX_0(n)^\naive \ra \sX(1)$ is not representable.
\end{pp-tweak}

\begin{pp-tweak}[The notation $d(m)$] \lab{notation-cm}
For a positive divisor $m$ of $n$, we set
\[
\tst d(m) \ce \f{m}{\gcd(m, \f{n}{m})},
\]
so that $d(m)$ depends \emph{both} on $m$ and on the integer $n$ that is fixed throughout. 

To explain the role of the function $m \mapsto d(m)$ in the context of $\Gamma_0(n)$-structures on generalized elliptic curves, let $E$ be the standard $m$-gon over an algebraically closed field and suppose that $E$ is equipped with an ample cyclic subgroup $G \subset E^\sm$ of order $n$. Then $G \cap (E^\sm)^0 = \mu_{\f{n}{m}}$ and $\mu_m \subset \Aut(E)$ is the subgroup of those automorphisms that induce the identity map on the contraction of $E$ with respect to the zero section (see \Cref{auto-n-gon}). The further subgroup of $\Aut(E)$ that in addition stabilizes $G$ is therefore $\mu_m \cap \mu_{\f{n}{m}}  = \mu_{\gcd(m, \f{n}{m})}$ (intersection in $(E^\sm)^0$), and this subgroup acts trivially on precisely $d(m)$ of the $m$ irreducible components of $E$. 

When refining $G$ to a $\Gamma_0(n)$-structure on such an $E$, we will only remember the contraction $c_{E^\sm[d(m)]}(E)$ that is a $d(m)$-gon together with the standard cyclic subgroup $G_{\f{n}{m} \cdot d(m)}$ of order $\f{n}{m} \cdot d(m)$. In addition, we will require the datum of a compatible ample cyclic $G'$ of order $n$ on every $E'$ that contracts to (a base change) of $c_{E^\sm[d(m)]}(E)$ and that has $m$-gon degenerate geometric fibers. Different $m$ may give the same $d(m)$, so there is no way to recover $m$ from $c_{E^\sm[d(m)]}(E)$ alone; to overcome this, we will incorporate $m$ into the data that comprises a $\Gamma_0(n)$-structure.
\end{pp-tweak}

For the precise definition of a $\Gamma_0(n)$-structure given in \S\ref{X0n-def}, we need the following preparations.

\begin{pp-tweak}[The stack of ``decontractions''] \lab{decon}
Fix a positive divisor $m$ of $n$ and suppose that we have a generalized elliptic curve $E \xra{\pi} S$ and an open subscheme $S_{\pi, (m)} \subset S$ that contains the elliptic curve locus $S - S^{\infty, \pi}$ and such that the degenerate geometric fibers of $E_{S_{\pi, (m)}}$ are $d(m)$-gons. (Such an $S_{\pi, (m)}$ will be part of the data of a $\Gamma_0(n)$-structure on $E$.) The base change $E_{S_{\pi, (m)}}$ determines a map $S_{\pi, (m)} \ra \EEl_{d(m)}$, so we may consider the fiber product algebraic stack 
\[
S_{\pi, (m)} \times_{\EEl_{d(m)}} \EEl_m,
\]
which parametrizes ``decontractions'' of $E_{S_{\pi, (m)}}$, or, more precisely, which, for variable $S_{\pi, (m)}$-schemes $S'$, parametrizes the pairs
\[
(E' \xra{\pi'} S',\, \iota' \colon E_{S'} \isomto c_{E'^{\sm}[d(m)]}(E'))
\]
consisting of a generalized elliptic curve $E' \xra{\pi'} S'$ whose degenerate geometric fibers are $m$-gons and a specified $S'$-isomorphism $\iota'$. We denote the universal object of $S_{\pi, (m)} \times_{\EEl_{d(m)}} \EEl_m$ by
\[
(\cE_{\pi, (m)}, \iota_{\pi, (m)}).
\]
The base change of $S_{\pi, (m)} \times_{\EEl_{d(m)}} \EEl_m$ (resp.,~of $\cE_{\pi, (m)}$) to $S - S^{\infty, \pi}$ is identified with $S - S^{\infty, \pi}$ (resp.,~with $E_{S - S^{\infty, \pi}}$), and the same holds over the entire $S_{\pi, (m)}$ if $d(m) = m$.

We will endow the universal ``decontraction'' $\cE_{\pi, (m)}$ with additional structures. The algebraic stack $\cE_{\pi, (m)}$ is typically not a scheme, but there are two ways to think about such structures concretely:
\begin{itemize}
\item
As compatible with isomorphisms and base change structures on $E'$ for each $(E' \xra{\pi'} S', \iota')$;

\item
As compatible under the pullbacks 
\[
\qqq S_{\pi, (m)} \times_{\EEl_{d(m)}} X_1 \rightrightarrows S_{\pi, (m)} \times_{\EEl_{d(m)}} X_0
\]
structures on the ``decontractions'' over the indicated bases, where $X_1 \rightrightarrows X_0 \ra \EEl_m$ is a once and for all fixed  scheme presentation of the algebraic stack $\EEl_m$, so that 
\[
\qqq S_{\pi, (m)} \times_{\EEl_{d(m)}} X_1 \rightrightarrows S_{\pi, (m)} \times_{\EEl_{d(m)}} X_0 \ra S_{\pi, (m)} \times_{\EEl_{d(m)}} \EEl_m
\]
is a scheme presentation of the algebraic stack $S_{\pi, (m)} \times_{\EEl_{d(m)}} \EEl_m$ (by \Cref{elln-props}~\ref{EN-a}, the algebraic stacks $\EEl_m$ and $\EEl_{d(m)}$ have finite diagonals, so $X_0 \times_{\EEl_m} X_0$ and similar fiber products that would \emph{a priori} be algebraic spaces are schemes).
\end{itemize}

The second way has the advantage of avoiding set-theoretic difficulties that would need to be addressed in order to make the first way completely rigorous.
\end{pp-tweak}

The contractions of the generalized elliptic curves parametrized by the stack $S_{\pi, (m)} \times_{\EEl_{d(m)}} \EEl_m$ are identified. In particular, the degenerate geometric fibers of these curves have canonically isomorphic component groups because the identity component of such a fiber may be used to fix the ``direction'' of the $m$-gon. This observation lies behind the following lemma:


\begin{lemma-tweak} \lab{identify}
Let $E \xra{\pi} S$ and $E' \xra{\pi'} S$ be generalized elliptic curves whose degenerate geometric fibers are $m$-gons and let $\iota\colon c(E) \isomto c(E')$ be an $S$-isomorphism between their contractions with respect to the identity sections.
\benum
\item  \lab{I-a}
If $S$ is a geometric point, then there is a unique identification 
\[
\qq E^\sm/(E^\sm)^0 = E'^\sm/(E'^\sm)^0
\]
of the component groups that is induced by any isomorphism $E \simeq E'$ that is compatible with $\iota$.

\item \lab{I-b}
If $S^{\mathrm{red}} = (S^{\infty, \pi})^{\mathrm{red}}$ {\upshape(}so that also $S^{\mathrm{red}} = (S^{\infty, \pi'})^{\mathrm{red}}${\upshape)}, then there is a unique $S$-identification 
\[
\qq (E^\sm)[m]/(E^\sm)^0[m] = (E'^\sm)[m]/(E'^\sm)^0[m]
\]
whose base change to any geometric $S$-point $\ov{s}$ is induced by any $\ov{s}$-isomorphism $E_{\ov{s}} \simeq E'_{\ov{s}}$ compatible with $\iota_{\ov{s}}$. Any $S$-isomorphism $i\colon E \simeq E'$  compatible with $\iota$ induces this identification.

\item \lab{I-c}
For $g \in E^\sm(S)$ and $g' \in E'^\sm(S)$, the set of $s \in S$ for which $g$ and $g'$ meet the same {\upshape(}in the sense of {\upshape\ref{I-a})} irreducible components of $E_{\ov{s}}$ and $E'_{\ov{s}}$ forms an open subscheme of $S$ that is also closed if $S^{\mathrm{red}} = (S^{\infty, \pi})^{\mathrm{red}}$.
\eenum
\end{lemma-tweak}

\bpf \hfill
\benum
\item
If either $E$ or $E'$ is smooth, then $\iota$ itself induces the desired identification. We may therefore assume that both $E$ and $E'$ are degenerate. Then, by Remark~\ref{gen-loc-proj}, both $E$ and $E'$ are isomorphic to the standard $m$-gon discussed in Remark~\ref{n-gon-gen-ell}.
Moreover, any two isomorphisms $E \simeq E'$ that are compatible with $\iota$ differ by an automorphism of $E'$ that is the identity map on $(E'^\sm)^0$. It remains to observe that, by \Cref{auto-n-gon}, any automorphism of $E'$ that is the identity map on $(E'^\sm)^0$ induces the identity map on $E'^\sm/(E'^\sm)^0$.

\item
If $S$ is a geometric point, then
\[
\qq (E^\sm)[m]/(E^\sm)^0[m] = E^\sm/(E^\sm)^0,
\]
and likewise for $E'$, so the claim follows from \ref{I-a}. In general, by \Cref{can-submult}, both 
\[
\qq (E^\sm)[m]/(E^\sm)^0[m] \qq \text{and} \qq (E'^\sm)[m]/(E'^\sm)^0[m]
\]
are \'{e}tale, so we may and do assume that $S = S^{\mathrm{red}}$. In this case, by Remark \ref{gen-loc-proj}, $i$ exists fppf locally on $S$. Moreover, any $i$ satisfies the defining property, so we only need to check that two different $i$ induce the same identification. For this, the case of a local strictly Henselian $S$ suffices and reduces to the settled case of a geometric point.

\item
We may assume that $S = S^{\infty, \pi} = S^{\infty, \pi'}$ and $S$ is reduced and may work fppf locally on $S$. We therefore use Remark \ref{gen-loc-proj} to fix an $S$-isomorphism  $i\colon E \isomto E'$ that is compatible with $\iota$ and to assume that $E$ is the standard $m$-gon. In this case, the label of the component of $E^\sm$ that meets $g$ is locally constant on $S$, and likewise for $\iota\i(g')$.
\qedhere
\eenum
\epf

\begin{pp-tweak}[Coherence of a cyclic subgroup of the universal ``decontraction''] \lab{coherence}
In the notation of \S\ref{decon}, part of the data of a $\Gamma_0(n)$-structure will be an ample cyclic $(S_{\pi, (m)} \times_{\EEl_{d(m)}} \EEl_m)$-subgroup
\[
\cG_{(m)} \subset \cE_{\pi, (m)}^\sm
\]
of order $n$, or, in more concrete terms, for every $(E' \xra{\pi'} S', \iota')$ an ample cyclic $S'$-subgroup $G' \subset E'^\sm$ of order $n$ that is compatible with base change and with isomorphisms of pairs $(E', \iota')$ (for the notion of cyclicity, see \Cref{def-cyclic}). 

In order to isolate a well-behaved class of such $\cG_{(m)}$, we say that $\cG_{(m)}$ is \emph{coherent} if:
\begin{quote}
For every $S_{\pi, (m)}$-scheme $S'$ and every pair of objects
\[
\q (E'_1 \xra{\pi'_1} S', \iota'_1) \qq \text{and} \qq (E'_2 \xra{\pi'_2} S', \iota'_2)
\]
of $(S_{\pi, (m)} \times_{\EEl_{d(m)}} \EEl_m)(S')$, the pullbacks $G_1' \subset E_1'^\sm$ and $G_2' \subset E_2'^\sm$ of $\cG_{(m)}$ fpqc locally on $S'$ have generators $g_1'$ and $g_2'$ that meet the same (in the sense of \Cref{identify}~\ref{I-a}) irreducible components of the geometric fibers of $E_1'$ and $E_2'$ and satisfy
\[
\tst  (\iota_1')\i(\f{n}{m} \cdot g_1') = (\iota_2')\i(\f{n}{m} \cdot g_2').
\]
\end{quote}
(The last equality takes place in $E$ and makes sense because $\f{n}{m} \cdot g_1'$ lies in the contraction $c_{E_1'^\sm[d(m)]}(E_1')$ by \Cref{std-subgp}~\ref{SS-d}, and likewise for $\f{n}{m} \cdot g_2'$.) Equivalently, the coherence of $\cG_{(m)}$ is a condition of the existence of compatible fpqc local generators of the pullbacks of $\cG_{(m)}$ along the two projections 
\[
(S_{\pi, (m)} \times_{\EEl_{d(m)}} \EEl_m) \times_{S_{\pi, (m)}} (S_{\pi, (m)} \times_{\EEl_{d(m)}} \EEl_m) \rightrightarrows S_{\pi, (m)} \times_{\EEl_{d(m)}} \EEl_m,
\]
where compatibility amounts to the conditions imposed on $g_1'$ and $g_2'$ above. 

In what follows, the purpose of the coherence condition is to ensure that $\cG_{(m)}$ is uniquely determined by its pullback to \emph{any} $(E' \xra{\pi'} S', \iota')$ with $S' = S_{\pi, (m)}$, provided that such an $(E', \iota')$ exists. \Cref{hop-G} will justify this, and its aspect \ref{hG-iii-pr} will show that no generality is lost if one strengthens the coherence condition by fixing an fpqc local generator $g_1'$ of $G_1'$ in advance.

Any $\cG_{(m)}$ is coherent if $S_{\pi, (m)} \times_{\EEl_{d(m)}} \EEl_m = S_{\pi, (m)}$, and also if $n$ is a unit on $S_{\pi, (m)}$ as we now~show.
\end{pp-tweak}

\begin{lemma-tweak} \lab{compatible-gens}
If $n$ is invertible on $S_{\pi, (m)}$, then every ample cyclic $(S_{\pi, (m)} \times_{\EEl_{d(m)}} \EEl_m)$-subgroup $\cG_{(m)} \subset \cE_{\pi, (m)}^\sm$ of order $n$ is coherent.
\end{lemma-tweak}

\bpf 
We will show that for every pair $(E'_1 \xra{\pi'_1} S', \iota'_1)$ and $(E'_2 \xra{\pi'_2} S', \iota'_2)$ as in the definition of coherence, desired generators $g_1'$ and $g_2'$ of $G_1'$ and $G_2'$ exist even \'{e}tale locally on $S'$. For this, due to \Cref{identify}~\ref{I-c}, we may assume that $S'$ is local strictly Henselian and that the special fibers $(E'_1)_{s'}$ and $(E'_2)_{s'}$ are degenerate. Moreover, since $(E_1')^\sm[n]$ and $(E_2')^\sm[n]$ are \'{e}tale and $G_1'$ and $G_2'$ are constant, we may assume further that $S'$ is a geometric point. In the case of a geometric point, it suffices to transport any choice of a $g_1'$ across any $S'$-isomorphism $(E_1', \iota_1') \simeq (E_2', \iota_2')$.
\epf

The following key lemma analyses the coherence condition beyond the case when $n$ is a unit by exhibiting a universal property satisfied by pullbacks of a coherent $\cG_{(m)}$. This property compensates for the loss of a direct reduction to geometric points that governed the case of an invertible $n$.

\begin{lemma-tweak} \lab{hop-G}
Let $m$ be a positive divisor of $n$, let $d \in \bZ_{\ge 1}$ be a multiple of $m$, let $E \xra{\pi} S$ and $E' \xra{\pi'} S$ be generalized elliptic curves whose degenerate geometric fibers are $d$-gons, and let
\[
\iota \colon c_{E^\sm[d(m)]}(E) \isomto c_{E'^\sm[d(m)]}(E')
\]
be an $S$-isomorphism. For every cyclic $S$-subgroup $G \subset E^\sm$ of order $n$ that meets precisely $m$ irreducible components of every degenerate geometric fiber of $E$, there is a unique cyclic $S$-subgroup $G' \subset E'^\sm$ of order $n$ such that 
\benumr
\item  \lab{hG-i}
over $S - S^{\infty, \pi} = S - S^{\infty, \pi'}$ there is an equality $\iota(G_{S - S^{\infty, \pi}}) = G'_{S - S^{\infty, \pi'}}${\upshape;} and

\item \lab{hG-iii}
fpqc locally on $S$ there exist generators $g$ of $G$ and $g'$ of $G'$ that meet the same irreducible components of the geometric fibers of $E$ and $E'$ {\upshape(}in the sense of Lemma~{\upshape\ref{identify}~\ref{I-a})} and satisfy 
\[
\tst \q \iota(\f{n}{m} \cdot g) = \f{n}{m} \cdot g'.
\]
{\upshape(}So $G'$ meets precisely $m$ irreducible components of every degenerate geometric fiber of $E'$.{\upshape)}
\eenum

Moreover, this unique $G'$ is such that
\benumr \addtocounter{enumi}{2}
\item \lab{hG-iii-pr}
for every $S$-scheme $T$ and every generator $\wt{g}$ of $G_{T}$, fpqc locally on $T$ there exists a generator $\wt{g}'$ of $G'_T$ such that $\wt{g}$ and $\wt{g}'$ meet the same irreducible components of the geometric fibers of $E$ and $E'$ and satisfy
\[
\tst\q \iota(\f{n}{m} \cdot \wt{g}) = \f{n}{m} \cdot \wt{g}';
\]

\item \lab{hG-iv-pr}
the standard cyclic subgroups $G_{\f{n}{m}\cdot d(m)} \subset G$ and $G_{\f{n}{m}\cdot d(m)}' \subset G'$ of order $\f{n}{m}\cdot d(m)$ satisfy
\[
\q \iota(G_{\f{n}{m}\cdot d(m)}) = G'_{\f{n}{m}\cdot d(m)}.
\]
\eenum
\end{lemma-tweak}

\begin{rem-tweak}
The equalities displayed in \ref{hG-iii}--\ref{hG-iv-pr} make sense due to \Cref{std-subgp}~\ref{SS-d}.  
\end{rem-tweak}

\bpf[Proof of Lemma~{\upshape\ref{hop-G}}]
We have broken the argument up into six steps.

\emph{Step {\upshape1.}~The claim of {\upshape\ref{hG-iv-pr}} follows from the rest.}
The subgroups $\iota(G_{\f{n}{m}\cdot d(m)})$ and $G'_{\f{n}{m}\cdot d(m)}$ of $E'^\sm$ are cyclic of order $\f{n}{m} \cdot d(m)$, agree with $\iota((G_{\f{n}{m}\cdot d(m)})_{S - S^{\infty, \pi}})$ over $S - S^{\infty, \pi'}$, and fpqc locally on $S$ have generators $\iota(\f{m}{d(m)} \cdot g)$ and $\f{m}{d(m)} \cdot g'$ whose $\f{n}{m}$-multiples equal $\iota(\f{n}{m} \cdot (\f{m}{d(m)} \cdot g))$. Therefore, $\iota(G_{\f{n}{m}\cdot d(m)})$ and $G'_{\f{n}{m}\cdot d(m)}$ must be equal because they satisfy \ref{hG-i} and \ref{hG-iii} when $n$, $m$, and $G$ are replaced by $\f{n}{m} \cdot d(m)$, $d(m)$, and $G_{\f{n}{m}\cdot d(m)}$, respectively ($G_{\f{n}{m}\cdot d(m)}$ meets precisely $d(m)$ irreducible components of every degenerate geometric fiber of $E$ due to \Cref{std-subgp}~\ref{SS-d}).

\emph{Step {\upshape2.}~The claim of {\upshape\ref{hG-iii-pr}}.}
We may assume that $T = S$ and may work fpqc locally on $S$, so we fix $g$, $g'$, and $\wt{g}$ over $S$. In order to find a desired fpqc local $\wt{g}'$, we work Zariski locally on $S$ and use limit arguments together with \Cref{identify}~\ref{I-c} to reduce to the case when $S = \Spec R$ for some Noetherian $R$. Then we pass to an fpqc cover to assume that $R$ is complete and separated with respect to the ideal $I$ that cuts out $S^{\infty, \pi}$ (equivalently, with respect to the ideal that cuts out $S^{\infty, \pi'}$; see \Cref{raise-it}). 

By \Cref{Gn-Hn}~\ref{Gn-Hn-a}, $E^\sm[n]$ (resp.,~$E'^\sm[n]$) has the largest finite locally free $S$-subgroup $A_{n, m}$ (resp.,~$A_{n, m}'$) that meets precisely $m$ irreducible components of every degenerate geometric fiber of $E$ (resp.,~$E'$), so $G \subset A_{n, m}$ and $G' \subset A_{n, m}'$. Moreover, \Cref{Gn-Hn}~\ref{Gn-Hn-a} supplies extensions
\[
\xymatrix{
0 \ar[r] & B_n \ar[r] \ar@{=}[d] & A_{n, m} \ar[r] & C_m \ar[r]\ar@{=}[d]  & 0 \\
0 \ar[r] & B_n \ar[r]  & A'_{n, m} \ar[r] & C_m \ar[r]  & 0
}
\]
of $S$-group schemes, where the identification of $B_n$ is via $\iota$ and the identification of $C_m$ is via \Cref{identify}~\ref{I-b} (applied over $R/I^j$ for every $j\ge 1$ to the contractions of $E_{R/I^j}$ and $E'_{R/I^j}$ with respect to the $m$-torsion). As may be checked on degenerate geometric fibers, the generators $g \in G(S)$ and $g' \in G'(S)$ project to the same section of $C_m$ that gives an isomorphism $C_m \simeq \bZ/m\bZ$. 

The homomorphism $G \ra C_m$ is finite locally free and, by \Cref{std-fact}~\ref{SF-a}, its kernel is the standard cyclic subgroup $G_{\f{n}{m}} \subset G$ of order $\f{n}{m}$. By replacing $g$ and $g'$ by $u\cdot g$ and $u \cdot g'$ for a suitable $u \in \underline{(\bZ/n\bZ)^\times}(S)$, we reduce to the case when $g$ and $\wt{g}$ have the same image in $C_m$. Then $g - \wt{g} \in G_{\f{n}{m}}$, so $\f{n}{m} \cdot g = \f{n}{m} \cdot \wt{g}$, which means that we may choose $\wt{g}'$ to be $g'$.

For the rest of the proof, we focus on the remaining claim about the existence and uniqueness of $G'$.

\emph{Step {\upshape3.}~Reduction to the case when $n$ is a prime power.}
The group $G$, as well as any candidate $G'$, decomposes as a product of its $p$-primary parts for various primes $p$ dividing $n$. By \cite{KM85}*{1.7.2}, cyclicity of $G$ or of $G'$ is equivalent to the cyclicity of the primary factors, and the datum of a generator of $G$ or of $G'$ corresponds to the datum of a generator of each primary factor. Therefore, for the existence and the uniqueness of the sought $G'$ we may assume that $n$ is a prime power.

For the rest of the proof, we assume that $n = p^r$ and $m = p^s$ for some prime $p$ and $r, s \in \bZ_{\ge 0}$.

\emph{Step {\upshape4.}~The case $s = 0$.}
For the existence,  $\iota(G)$ fulfills the requirements \ref{hG-i}--\ref{hG-iii}. The uniqueness reduces to the case of an Artinian local $S$ and then follows from \Cref{Gn-Hn}~\ref{Gn-Hn-a}.

For the rest of the proof, we assume that $s \ge 1$, so that $\f{n}{m} \neq n$.

\emph{Step {\upshape5.}~Uniqueness of $G'$.}
Due to the claim concerning \ref{hG-iii-pr} (i.e.,~due to Step 2), we may assume that the two candidates $G_1', G_2' \subset E'^{\sm}$ have generators $g_1'$ and $g_2'$ that meet the same irreducible components of the geometric fibers of $E'$ and satisfy $\f{n}{m} \cdot  g_1' = \f{n}{m} \cdot g_2'$. Furthermore, we may assume that the base $S$ is Noetherian, then local, then complete, and finally Artinian, and that $E'$ is nonsmooth over $S$. Then, since $g_1' - g_2' \in (E'^\sm)^0(S)$ and $\f{n}{m} \cdot g_1' = \f{n}{m} \cdot g_2'$, we have 
\[
\tst g_2' = g_1' + h \qq \text{for some} \q h \in (E'^\sm)^0[\f{n}{m}](S).
\]
By \Cref{can-submult} and \Cref{std-fact}~\ref{SF-a}, the $S$-group $(E'^\sm)^0[\f{n}{m}]$ is the standard cyclic subgroup of $G_1'$ of order $\f{n}{m}$, so \Cref{std-subgp}~\ref{SS-e} ensures that $g_1' + h$ generates $G_1'$, which means that $G_1' = G_2'$.

\emph{Step {\upshape6.}~Existence of $G'$.}
Due to the uniqueness of $G'$, for its existence we may work fpqc locally on $S$, so we fix a generator $g$ of $G$. Moreover, as in Step 2 we reduce to the case when $S = \Spec R$ for a Noetherian $R$ that is complete and separated with respect to the ideal $I \subset R$ that cuts out $S^{\infty, \pi}$ and use \Cref{Gn-Hn}~\ref{Gn-Hn-a} to obtain the diagram of extensions displayed in Step 2.

By \Cref{Gn-Hn}~\ref{Gn-Hn-a}, $E'^\sm[m] \subset A_{n, m}'$, so $E'^\sm[\f{m}{d(m)}] \subset A_{n, m}'$, too, and hence the image of $A_{n, m}'$ under the multiplication by $\f{m}{d(m)}$ map of $E'^\sm$ is a finite locally free $S$-subgroup of $A'_{\f{n}{m} \cdot d(m), d(m)}$ of order $\p{\f{n}{m} \cdot d(m)} \cdot d(m)$. This image therefore equals $A'_{\f{n}{m} \cdot d(m), d(m)}$, so, since $\iota(\f{m}{d(m)} \cdot g)$ lies in $A'_{\f{n}{m} \cdot d(m), d(m)}$, after replacing $S$ by a finite locally free cover we may choose a $g' \in A'_{n, m}(S)$ with 
\[
\tst \f{m}{d(m)}\cdot g' = \iota(\f{m}{d(m)} \cdot g).
\]
Since $E'^\sm[\f{m}{d(m)}]$ is an extension of $(C_m)[\f{m}{d(m)}]$ by $(B_n)[\f{m}{d(m)}]$, after a further finite locally free cover of $S$ we may adjust $g'$ by a lift to $(E'^\sm[\f{m}{d(m)}])(S)$ of the difference of the images of $g$ and $g'$ in $C_m$ to arrange that $g$ and $g'$ have the same image in $C_m$ and hence meet the same irreducible components of the geometric fibers of $E$ and $E'$.

By \Cref{restrict}~\ref{R-d}, $g'$ generates a cyclic $S$-subgroup $G' \subset E'^\sm$ of order $n$. Since $\f{m}{d(m)}\mid \f{n}{m}$, the group $G'$ satisfies \ref{hG-iii}. Thus, to complete Step 6, and hence also the proof of \Cref{hop-G}, it suffices to show that 
\[
\iota(G_{S - S^{\infty, \pi}}) = G'_{S - S^{\infty, \pi'}}.
\]
We have $G \subset A_{n, m}$ and $G' \subset A'_{n, m}$ with $g$ and $g'$ projecting to the same section of $C_m$. Moreover, by \Cref{Gn-Hn}~\ref{Gn-Hn-b} and the diagram displayed in Step 2, both $\iota((A_{n, m})_{S - S^{\infty, \pi'}})$ and $(A'_{n, m})_{S - S^{\infty, \pi'}}$ are the preimages in $E'_{S - S^{\infty, \pi'}}[n]$ of the unique $(S - S^{\infty, \pi'})$-subgroup of $(E'^\sm[n]/B_n)_{S - S^{\infty, \pi'}}$ of order $m$, so 
\[
\iota \q \text{identifies} \q A_{n, m} \q \text{and} \q  A'_{n, m} \qq \text{over $S - S^{\infty, \pi'}$}.
\]
We claim that under this identification via $\iota$, the image of $g_{S - S^{\infty, \pi}}$ in $A_{n, m}/B_n$ agrees with the image of $g'_{S - S^{\infty, \pi'}}$ in $A'_{n, m}/B_n$. Since $A'_{n, m}/B_n$ is finite \'{e}tale, it suffices to check the claimed agreement on the geometric fibers at the points in $S - S^{\infty, \pi'}$, so the technique used in the proof of \Cref{Gn-Hn}~\ref{Gn-Hn-b} reduces the proof of the claimed agreement to the case when $R$ is a discrete valuation ring and $E$ and $E'$ have smooth generic fibers but nonsmooth closed fibers. In this case, by \Cref{extend-iso}~\ref{EI-b}, $\iota$ extends to a unique isomorphism $E \simeq E'$, which then must induce the identification of the groups $C_m$ for $E$ and $E'$. Thus, in this case the claimed agreement follows from the agreement of the images of $g$ and $g'$ in $C_m$.

Returning to the proof of $\iota(G_{S - S^{\infty, \pi}}) = G'_{S - S^{\infty, \pi'}}$, via the above reasoning, we conclude that $g'_{S - S^{\infty, \pi'}} - \iota(g_{S - S^{\infty, \pi}})$ lies in $B_n$. Moreover, since $\f{m}{d(m)}\mid \f{n}{m}$, the construction of $g'$ ensures that  
\[
\tst \f{n}{m} \cdot g'_{S - S^{\infty, \pi'}} = \f{n}{m} \cdot \iota(g_{S - S^{\infty, \pi}}).
\]
Therefore, there is an $h \in ((B_n)[\f{n}{m}])(S - S^{\infty, \pi'})$ such that
 \[
 g'_{S - S^{\infty, \pi'}} = \iota(g_{S - S^{\infty, \pi}}) + h.
 \]
By the uniqueness aspect of the first assertion of \Cref{Gn-Hn}~\ref{Gn-Hn-a} and by \Cref{std-subgp}~\ref{SS-d}, $(B_n)[\f{n}{m}]$ is the standard cyclic subgroup of $G$ of order $\f{n}{m}$, so $\iota(g_{S - S^{\infty, \pi}}) + h$ generates $\iota(G_{S - S^{\infty, \pi}})$ by \Cref{std-subgp}~\ref{SS-e}. The sought equality $\iota(G_{S - S^{\infty, \pi}}) = G'_{S - S^{\infty, \pi'}}$ follows.
\epf

We are ready for the definition of a $\Gamma_0(n)$-structure on a generalized elliptic curve.

\begin{pp-tweak}[$\Gamma_0(n)$-structures] \lab{Gamma-0-def}
For a generalized elliptic curve $E \xra{\pi} S$, a \emph{$\Gamma_0(n)$-structure on $E$} is a tuple
\[
(G,\, \{ S_{\pi, (m)} \}_{m\mid n},  \, \{\cG_{(m)} \}_{m\mid n})
\]
consisting of the following data.
\benuma
\item \lab{X0n-data-2}
A cyclic $(S - S^{\infty, \pi})$-subgroup $G \subset E_{S - S^{\infty, \pi}}$ of order $n$ (in the sense of \Cref{def-cyclic}).

\item \lab{X0n-data-3}
For each positive divisor $m$ of $n$, an open subscheme $S_{\pi, (m)} \subset S$ such that
\begin{enumerate}[labelindent=5pt,leftmargin=*,label={(2.\arabic*)}]
\item \lab{X0n-data-31}
$S =  \bigcup_{m}  S_{\pi, (m)}$;

\item \lab{X0n-data-32}
if $m \neq m'$, then $S_{\pi, (m)} \cap S_{\pi, (m')} = S - S^{\infty, \pi}$;

\item \lab{X0n-data-33}
the degenerate geometric fibers of $E_{S_{\pi, (m)}}$ are $d(m)$-gons, where $d(m) = \f{m}{\gcd(m, \f{n}{m})}$.
\eenum

\item \lab{X0n-data-5}
For each positive divisor $m$ of $n$, in the notation of \S\ref{decon}, an ample cyclic $(S_{\pi, (m)} \times_{\EEl_{d(m)}} \EEl_m)$-subgroup 
\[
\qq \cG_{(m)} \subset \cE^\sm_{\pi, (m)}
\]
of order $n$ such that

\begin{enumerate}[labelindent=5pt,leftmargin=*,label={(3.\arabic*)}]
\item \lab{X0n-data-51}
on the elliptic curve locus, 
\[
\qq (\cG_{(m)})_{S - S^{\infty, \pi}} = \iota_{\pi, (m)}(G);
\] 

\item \lab{X0n-data-53}
the cyclic subgroup $\cG_{(m)}$ is coherent in the sense of \S\ref{coherence}.
\eenum
\eenum
\end{pp-tweak}

\brems
\remi \lab{same-thing}
If $E  \ra S$ is smooth, then the data \ref{X0n-data-3}--\ref{X0n-data-5} are uniquely determined by \ref{X0n-data-2} and a $\Gamma_0(n)$-structure on $E$ is nothing else than a cyclic $S$-subgroup of order $n$.

\remi
If $n$ is invertible on $S$, then, by \Cref{compatible-gens}, the requirement \ref{X0n-data-53} is superfluous.

\remi \lab{sqfree-same-thing}
If $n$ is squarefree, then $d(m) = m$ for every $m$, so that $S_{\pi, (m)}$ is the open subscheme of $S$ obtained by removing all the $S^{\infty, \pi, m'}$ with $m' \neq m$, the ``decontraction'' $\cE_{\pi, (m)}$ is $E_{S_{\pi, (m)}}$ itself, and a $\Gamma_0(n)$-structure on $E$ is nothing else than an ample cyclic $S$-subgroup of $E^\sm$ order $n$.

In general, the datum $\{ S_{\pi, (m)} \}_{m\mid n}$ of \ref{X0n-data-3} is equivalent to a subdivision 
\[
\qq \tst S^{\infty, \pi} =  \bigsqcup_{m\mid n}  S^{\infty}_{\pi, (m)},
\]
subject to the requirement that $S^{\infty}_{\pi, (m)} \subset S^{\infty, \pi, d(m)}$ for every $m$. In this notation,
\[
\qq \tst S_{\pi, (m)} = S - \p{\bigcup_{m' \neq m} S^{\infty}_{\pi, (m')}}.
\]

\remi
The subgroup $\cG_{(m)}$ determines an ample cyclic $S_{\pi, (m)}$-subgroup 
\[
\qq G_{(m)} \subset E_{S_{\pi, (m)}}^\sm
\]
of order $\f{n}{m} \cdot d(m)$ such that $(G_{(m)})_{S - S^{\infty, \pi}}$ is a standard cyclic subgroup of $G$. To build $G_{(m)}$, we choose an fppf cover $S'$ of $S_{\pi, (m)}$ for which there is an object $(E' \ra S', \iota')$ of $S_{\pi, (m)} \times_{\EEl_{d(m)}} \EEl_m$, let $G' \subset E'^\sm$ be the pullback of $\cG_{(m)}$, and use \Cref{std-subgp}~\ref{SS-d} to set 
\[
\qq (G_{(m)})_{S'} \ce (\iota')\i(G'_{\f{n}{m} \cdot d(m)}).
\]
\Cref{hop-G}~\ref{hG-iv-pr} shows the agreement of the two pullbacks of $(G_{(m)})_{S'}$ to $S' \times_{S_{\pi, (m)}} S'$, and hence also the effectivity of descent to the sought $G_{(m)}$ over $S_{\pi, (m)}$, as well as the independence of the resulting $G_{(m)}$ on the choice of $S'$ and $(E', \iota')$.

By construction and \Cref{hop-G}~\ref{hG-iv-pr}, $\iota_{\pi, (m)}(G_{(m)})$  is a standard cyclic subgroup of $\cG_{(m)}$.
\erems

The principal reason why the stack $\sX_0(n)$ that we are about to introduce is practical to work with even when $n$ is not squarefree is \Cref{X0-cart}~\ref{X0C-a} below. 

\begin{pp-tweak}[The stack $\sX_0(n)$] \lab{X0n-def}
In order to construct this $\bZ$-stack, we begin by letting $S$ be a variable scheme and by defining the categories $\sX_0(n)(S)$. 

The objects of $\sX_0(n)(S)$ are the tuples
\[
(E \xra{\pi} S,\, G,\, \{ S_{\pi, (m)} \}_{m\mid n}, \, \{\cG_{(m)} \}_{m\mid n})
\]
consisting of a generalized elliptic curve $E \xra{\pi} S$ and a $\Gamma_0(n)$-structure on $E$.

In $\sX_0(n)(S)$, a morphism 
\[
(E_1 \xra{\pi_1} S,\, G_1,\, \{ S_{\pi_1, (m)} \},  \, \{\cG_{(m), 1} \}) \ra (E_2 \xra{\pi_2} S,\, G_2,\, \{ S_{\pi_2, (m)} \},  \, \{\cG_{(m),2 } \}) 
\]
between two tuples such that $S_{\pi_1, (m)} = S_{\pi_2, (m)}$ for every positive divisor $m$ of $n$ consists of 
\begin{enumerate}[label={(\Roman*)}]
\item
an $S$-isomorphism $i_E\colon E_1 \isomto E_2$ of generalized elliptic curves such that 
\[
\q (i_E)_{S - S^{\infty, \pi_1}}(G_1) = G_2;
\]

\item \lab{X0n-def-2}
for each positive divisor $m$ of $n$, an isomorphisms $i_{(m)}$ of stacks over $S_{\pi_1, (m)} = S_{\pi_2, (m)}$ and an isomorphism $i_{\cE_{(m)}}$ of generalized elliptic curves that fit into the commutative diagram
\[
\q\xymatrix{
\cE_{\pi_1, (m)} \ar[r]_-{i_{\cE_{(m)}}}^-{\sim} \ar[d] & \cE_{\pi_2, (m)} \ar[d] \\
S_{\pi_1, (m)} \times_{\EEl_{d(m)}} \EEl_m \ar[r]^{\sim}_{i_{(m)}} & S_{\pi_2, (m)} \times_{\EEl_{d(m)}} \EEl_m}
\]
and such that $i_{\cE_{(m)}}$ induces the isomorphism $(i_E)_{S_{\pi_1, (m)} \times_{\EEl_{d(m)}} \EEl_m}$ between the contractions of $\cE_{\pi_1, (m)}$ and $\cE_{\pi_2, (m)}$ with respect to $\cE_{\pi_1, (m)}^\sm[d(m)]$ and $\cE_{\pi_2, (m)}^\sm[d(m)]$, respectively, and satisfies
\[
\q i_{\cE_{(m)}}(\cG_{(m), 1}) = \cG_{(m), 2}.
\]
\eenum

There are no morphisms between tuples for which $S_{\pi_1, (m)} \neq S_{\pi_2, (m)}$ for some $m$.

In concrete terms, the datum $(i_{(m)}, i_{\cE_{(m)}})$ of \ref{X0n-def-2} amounts to 
\begin{enumerate}[label={(\Roman*$'$)}] \addtocounter{enumi}{1}
\item
an $S_{\pi_1, (m)}$-isomorphism 
\[
\q i_{(m)}\colon S_{\pi_1, (m)} \times_{\EEl_{d(m)}} \EEl_m \isomto S_{\pi_2, (m)} \times_{\EEl_{d(m)}} \EEl_m
\]
together with: for every object $(E'_1 \ra S',\, \iota'_1)$ of $S_{\pi_1, (m)} \times_{\EEl_{d(m)}} \EEl_m$ with $i_{(m)}$-image $(E'_2 \ra S',\, \iota'_2)$, a generalized elliptic curve isomorphism
\[
\q i_{E'_1, E_2'}\colon E'_1 \isomto E_2'
\]
that is compatible with $(i_E)_{S'}$ (via $\iota_1'$ and $\iota_2'$), brings the pullback of $\cG_{(m), 1}$ to the pullback of $\cG_{(m), 2}$, and whose formation commutes with isomorphisms and base change of pairs $(E_1', \iota_1')$.
\eenum

A compatible with $i_E$ pair of isomorphisms $(i_{(m)}, i_{\cE_{(m)}})$ always exists (send $(E_1', \iota_1')$ to $(E_1', \iota_1' \circ (i_E)_{S'}\i)$) and, thanks to $i_{\cE_{(m)}}$, is unique up to a unique isomorphism. However, this unique $(i_{(m)}, i_{\cE_{(m)}})$ may not automatically respect $\cG_{(m), 1}$ and $\cG_{(m), 2}$. In practice, the uniqueness up to a unique isomorphism means that the lack of canonicity in the choice of $(i_{(m)}, i_{\cE_{(m)}})$ does not matter and that the construction of $\sX_0(n)$ stays in the realm of $2$-categories.

The existence of a unique $(i_{(m)}, i_{\cE_{(m)}})$ compatible with $i_E$ ensures that
\begin{itemize}
\item
$\sX_0(n)(S)$ is a groupoid; and 

\item 
the base change functor $\sX_0(n)(S) \ra \sX_0(n)(S')$ along variable scheme morphisms $S' \ra S$ turns $\sX_0(n)$ into a $\bZ$-stack for the fppf topology (see \cite{SP}*{\href{http://stacks.math.columbia.edu/tag/026F}{026F}} for stack axioms).
\end{itemize}

We let 
\[
\sX_0(n)^\infty \subset \sX_0(n) \qq \text{and} \qq \sY_0(n) \subset \sX_0(n)
\]
be the closed substack cut out by the degeneracy loci $S^{\infty, \pi}$ and its open complement (the elliptic curve locus), respectively. By Remark~\ref{same-thing}, there is an identification 
\[
\sY_0(n) = \sY_0(n)^\naive.
\]
By Remark \ref{sqfree-same-thing}, if $n$ is squarefree, then $\sX_0(n)$ is identified with $\sX_0(n)^\naive$.

For a positive divisor $m$ of $n$, we let 
\[
\sX_0(n)_{(m)} \subset \sX_0(n)
\]
be the open substack cut out by the subschemes $S_{\pi, (m)}$. For every tuple classified by $\sX_0(n)_{(m)}$, the degenerate geometric fibers of $E$ are $d(m)$-gons.
\end{pp-tweak}

\begin{pp-tweak}[The contraction $\sX_0(n)^\naive \ra \sX_0(n)$] \lab{forget-naive}
Let $E \xra{\pi} S$ be a generalized elliptic curve equipped with a naive $\Gamma_0(n)$-structure, i.e.,~with an ample cyclic $S$-subgroup $G \subset E^\sm$ of order $n$. To build a $\Gamma_0(n)$-structure on a generalized elliptic curve $\wt{E} \xra{\wt{\pi}} S$ out of $(E, G)$, we first construct $\wt{E}$ by letting $S_{\wt{\pi}, (m)}$, for a positive divisor $m$ of $n$, be the largest open subscheme of $S$ over which the degenerate geometric fibers of $E$ are $m$-gons and by letting $\wt{E}$ be the glueing of the contractions $c_{E^\sm[d(m)]}(E_{S_{\wt{\pi}, (m)}})$ along $E_{S - S^{\infty, \pi}}$. We endow $\wt{E}_{S - S^{\infty, \wt{\pi}}}$ with the cyclic subgroup $G_{S - S^{\infty, \pi}}$ of order $n$. This produces the data \ref{X0n-data-2} and \ref{X0n-data-3}, so it remains to explain how to get \ref{X0n-data-5}.

For a fixed positive divisor $m$ of $n$, each $S_{\wt{\pi}, (m)}$-scheme $S'$, and each generalized elliptic curve $E' \ra S'$ whose degenerate geometric fibers are $m$-gons and that is equipped with an $S'$-isomorphism 
\[
\iota' \colon \wt{E}_{S'} = c_{E^\sm[d(m)]}(E_{S'}) \isomto c_{E'^\sm[d(m)]}(E'),
\]
we endow $E'$ with the unique cyclic $S'$-subgroup $G'$ of order $n$ supplied by \Cref{hop-G}. Due to the uniqueness, the formation of $G'$ commutes with base change and with isomorphisms of pairs $(E', \iota')$. In other words, the subgroups $G'$ give rise to a cyclic subgroup $\cG_{(m)} \subset \cE_{\pi, (m)}^\sm$ of order $n$, which agrees with $G$ on the elliptic curve locus due to \Cref{hop-G} \ref{hG-i}, is ample due to \Cref{hop-G}~\ref{hG-iii}, and is coherent due to \Cref{hop-G}~\ref{hG-iii-pr}. This gives the sought datum \ref{X0n-data-5}.

The construction of $\wt{E}$ and of its $\Gamma_0(n)$-structure respects isomorphisms and base change of pairs $(E, G)$, so we obtain the sought contraction morphism
\[
\sX_0(n)^\naive \ra \sX_0(n),
\]
which for each positive divisor $m$ of $n$ restricts to a morphism
\[
\sX_0(n)^\naive_{(m)} \ra \sX_0(n)_{(m)}.
\]
\end{pp-tweak}

The following lemma together with \Cref{hop-G} is the driving force of our analysis of $\sX_0(n)$.

\begin{lemma-tweak} \lab{X0-cart}
Let $m$ be a positive divisor of $n$.
\benum
\item \lab{X0C-a}
The square
\[
\qq\xymatrix{
\sX_0(n)^\naive_{(m)} \ar[r] \ar[d] & \EEl_m \ar[d] \\
\sX_0(n)_{(m)} \ar[r] & \EEl_{d(m)}
}
\]
is Cartesian.

\item \lab{X0C-b}
The map $\sX_0(n)_{(m)} \ra \EEl_{d(m)}$ is representable by schemes, of finite presentation, separated, quasi-finite, and flat{\upshape;} moreover, it is \'{e}tale over $\bZ[\f{1}{n}]$.
\eenum
\end{lemma-tweak}

\bpf \hfill
\benum
\item
For a generalized elliptic curve $E \xra{\pi} S$, part of the data of a $\Gamma_0(n)$-structure $\gA$ on $E$ with $S_{\pi, (m)} = S$ is the datum of a naive $\Gamma_0(n)$-structure $G'$ on $E'$ for every $(E' \xra{\pi'} S, \iota')$ classified by $S_{\pi, (m)} \times_{\EEl_{d(m)}}\EEl_m$. The assignment of this naive $\Gamma_0(n)$-structure gives the morphism
\[
\qq \sX_0(n)_{(m)} \times_{\EEl_{d(m)}} \EEl_m \ra \sX_0(n)^\naive_{(m)},
\]
which, by construction of the contraction $\sX_0(n)^\naive_{(m)} \ra \sX_0(n)_{(m)}$ in \S\ref{forget-naive}, is a left inverse to the induced morphism 
\[
\qq \sX_0(n)^\naive_{(m)} \ra \sX_0(n)_{(m)} \times_{\EEl_{d(m)}} \EEl_m.
\]
To prove that it is also a right inverse, we need to argue that $\gA$ agrees with the $\Gamma_0(n)$-structure on $E$ determined as in \S\ref{forget-naive} by the naive $\Gamma_0(n)$-structure $G'$ on $E'$. For this, the key point is the coherence requirement \ref{X0n-data-53} on the $\cG_{(m)}$ that is part of $\gA$: thanks to it and to the uniqueness aspect of \Cref{hop-G}, for every $(E'' \xra{\pi''} S, \iota'')$ classified by $S_{\pi, (m)} \times_{\EEl_{d(m)}}\EEl_m$, the naive $\Gamma_0(n)$-structure $G''$ on $E''$ that is part of $\gA$ is also the one determined by $G'$ through \Cref{hop-G}, and likewise over any $S$-scheme $S'$.

\item
We prove the asserted properties with the representability by schemes requirement replaced by representability by algebraic spaces---due to \Cref{rep-schemes}, this loses no generality.

By \Cref{from-bottom}~\ref{FB-a} (applied with $m = 1$ there), $\sX_0(n)^\naive_{(m)} \ra \EEl_m$
enjoys all the properties in question. Moreover, these properties are fppf local on the target (for the representability by algebraic spaces, see \cite{SP}*{\href{http://stacks.math.columbia.edu/tag/04SK}{04SK}} or \cite{LMB00}*{10.4.2}) and, by \Cref{Bnm-input}~\ref{Bnm-a}, $\EEl_m \ra \EEl_{d(m)}$ is surjective, flat, and of finite presentation. With the help of \ref{X0C-a}, we therefore conclude that $\sX_0(n)_{(m)} \ra \EEl_{d(m)}$ inherits the properties in question.
\qedhere
\eenum
\epf

We are ready for the sought identification $\sX_0(n) = \sX_{\Gamma_0(n)}$ and for the regularity of $\sX_{\Gamma_0(n)}$.

\begin{thm-tweak} \lab{X0n-main} \hfill
\benum
\item \lab{X0nM-a}
The stack $\sX_0(n)$ is Deligne--Mumford and regular. The map $\sX_0(n) \ra \sX(1)$ that forgets the $\Gamma_0(n)$-structure and contracts with respect to the identity section induces the identification
\[
\qq \sX_0(n) = \sX_{\Gamma_0(n)};
\]
more precisely, $\sX_0(n)$ and $\sX_{\Gamma_0(n)}$ are the normalizations of $\sX(1)$ in $\sY_0(n)_{\bZ[\f{1}{n}]} \cong \sY_{\Gamma_0(n)}[\f{1}{n}]$.

\item \lab{X0nM-b}
The substack $\sX_0(n)^\infty \subset \sX_0(n)$ is a reduced relative effective Cartier divisor over $\Spec \bZ$ that meets every irreducible component of every geometric fiber of $\sX_0(n) \ra \Spec \bZ$ and is smooth over $\bZ[\f{1}{n}]$.

\eenum
\end{thm-tweak}

\bpf \hfill
\benum
\item
We will use the axiomatic \Cref{axiom}. To apply its part \ref{axiom-a}, and hence to prove the algebraicity of $\sX_0(n)$ and the quasi-compactness and separatedness of $\Delta_{\sX_0(n)/\bZ}$, we use the open cover $\sX_0(n) = \bigcup_{m\mid n} \sX_0(n)_{(m)}$ and appeal to \Cref{X0-cart}~\ref{X0C-b}. To then apply \Cref{axiom}~\ref{axiom-b}, and hence to prove the regularity of $\sX_0(n)$, we let $\sX(n) \ra \sX_0(n)$ be the composition of the contractions 
\[
\qq \sX(n) \ra \sX_0(n)^\naive \qq \text{and} \qq  \sX_0(n)^\naive \ra \sX_0(n)
\]
of \S\ref{naive-X0n-def} and \S\ref{forget-naive} and note that this composition is proper, flat, and surjective due to \S\ref{naive-X0n-def}, \Cref{X0-cart}~\ref{X0C-a}, and \Cref{Bnm-input}~\ref{Bnm-a}. Finally, in order to prove that $\sX_0(n)$ is Deligne--Mumford and $\sX_0(n) = \sX_{\Gamma_0(n)}$, by \Cref{axiom}~\ref{axiom-c}, we need to prove that the map 
\[
\qq \sX_0(n) \ra \sX(1)
\]
is representable by algebraic spaces and that its base change to $\sY(1)_{\bZ[\f{1}{n}]}$ is identified with 
\[
\qq\tst \sY_{\Gamma_0(n)}[\f{1}{n}] \ra \sY(1)_{\bZ[\f{1}{n}]}.
\]
Since $\sY_0(n) = \sY_0(n)^\naive$, the latter identification results from the fact that the image of $\Gamma_0(n)$ in $\GL_2(\bZ/n\bZ)$ is the stabilizer of the subgroup $\bZ/n\bZ \times \{ 0\}$ in $(\bZ/n\bZ)^2$ (compare with the proof of \Cref{X1-grand}~\ref{X1G-c}).

Due to \Cref{rep-crit}~\ref{RC-b}, the representability of $\sX_0(n) \ra \sX(1)$ will follow once we prove that, for every Artinian local algebra $A$ over an algebraically closed field $\ov{k}$ and every $\xi \in \sX_0(n)(\ov{k})$, no nonidentity automorphism of $\xi|_A$ maps to an identity automorphism in $\sX(1)(A)$. More concretely, by \Cref{auto-n-gon}, we need to prove that for every positive divisor $d$ of $n$ and every prime divisor $p$ of $d$, there is no $\Gamma_0(n)$-structure $\gA$ on the standard $d$-gon $E$ over $\ov{k}$ such that some nonidentity automorphism $i \in \mu_p(A) \subset \Aut(E)(A)$ fixes the pullback $\gA_A$ of $\gA$ to $A$. For the sake of contradiction, we fix such $\gA$ and $i$.

We let $m$ be such that $\gA$ has $S_{\pi, (m)} \neq \emptyset$, so, in particular, $d(m) = d$. We let $(\wt{E}, \iota)$ be the standard $m$-gon over $\ov{k}$ equipped with the canonical isomorphism $\iota\colon E \isomto c_{\wt{E}^\sm[d]}(\wt{E})$. Up to unique isomorphism, the pair of isomorphisms $(i_{(m)}, i_{\cE_{(m)}})$ that extends $i$ as in \S\ref{X0n-def} sends $(\wt{E}_A, \iota_A)$ to $(\wt{E}_A, \iota_A \circ i\i)$, so the ample cyclic $A$-subgroups $\wt{G} \subset \wt{E}_A^\sm$ and  $\wt{G}' \subset \wt{E}_A^\sm$ of order $n$ that are the pullbacks of $\cG_{(m)}$ corresponding to $(\wt{E}_A, \iota_A)$ and $(\wt{E}_A, \iota_A \circ i\i)$ must be equal: 
\[
\qqq \wt{G} = \wt{G}' \qq \text{inside}\q \wt{E}_A.
\]
We replace $A$ by an Artinian local fppf cover to assume that the automorphism $\iota_A \circ i \circ \iota_A\i$ of $c_{\wt{E}_A^\sm[d]}(\wt{E}_A)$ is the contraction of an automorphism 
\[
\qq \wt{i} \in \mu_{m}(A) \subset \Aut(\wt{E})(A).
\]
Then $\wt{i}$ gives an isomorphism $(\wt{E}_A, \iota_A \circ i\i)\isomto (\wt{E}_A, \iota_A)$, so must satisfy 
\[
\qqq \wt{i}(\wt{G}') = \wt{G}, \qq \text{i.e.,} \qq \wt{i}(\wt{G}) = \wt{G}.
\]
The latter equality means that $\wt{i}$ also lies in $\wt{G} \cap (\wt{E}^\sm_A)^0 = (\mu_{\f{n}{m}})_A$, that is, 
\[
\qq \wt{i} \in \mu_{\gcd(m, \f{n}{m})}(A).
\]
However, $\mu_{\gcd(m, \f{n}{m})}$ acts trivially on $c_{\wt{E}^\sm[d(m)]}(\wt{E})$ by the definition of $d(m)$ (see \S\ref{notation-cm}), which means that $\iota_A \circ i \circ \iota_A\i = \id$ and contradicts the assumption that $i \neq \id$.

\item
By the proof of \ref{X0nM-a}, $\sX(n) \ra \sX_0(n)$ is surjective, so the claim about the geometric fibers follows from the corresponding claim for $\sX(n)^\infty \subset \sX(n)$ proved in \Cref{Xn-Elln}~\ref{XE-b}. 

For the rest, we may work on $\sX_0(n)_{(m)}$ and may focus on the corresponding claims for 
\[
\qq \sX_0(n)^\infty_{(m)} \ce \sX_0(n)_{(m)} \cap \sX_0(n)^\infty,
\]
so it suffices to observe that $\sX_0(n)_{(m)}^\infty$ is the preimage of $\EEl_{d(m)}^\infty$ under the map 
\[
\qq \sX_0(n)_{(m)} \ra \EEl_{d(m)},
\]
to apply \Cref{elln-props}~\ref{EN-c}--\ref{EN-d} and \Cref{X0-cart}~\ref{X0C-b}, and to use the (R$_0$)$+$(S$_1$) criterion for reducedness. \qedhere
\eenum
\epf


\section{Implications for coarse moduli spaces} \lab{coarse-spaces}

The main goal of this chapter is to take advantage of the moduli interpretation of $\sX_0(n)$ presented in Chapter \ref{Gamma-0-case} to prove that the coarse moduli space $X_0(n)$ is regular at the cusps (and, in fact, regular on a large open subscheme, see \Cref{coarse-reg}). This regularity is not new: \cite{Edi90}*{\S1.2} uses the results of Katz and Mazur to verify via an explicit computation that the completion of $X_0(n)$ along the cusps is regular (such regularity is also a special case of an earlier assertion of Gross and Zagier, see \cite{GZ86}*{Prop.~III.1.4}). In contrast, the proof given below rests on \Cref{X0n-main}~\ref{X0nM-a}, but requires no computation of completions. 

We also exploit \Cref{coarse-miracle} to obtain a base change result for coarse moduli spaces $X_H$ of arbitrary congruence level $H$ (see \Cref{coarse-H}). To prepare for it, we review general properties of $X_H$.

\begin{pp-tweak}[The coarse moduli space of $\sX_H$] \lab{coarse-def}
For an open subgroup $H \subset \GL_2(\wh{\bZ})$, the finite type Deligne--Mumford $\bZ$-stack $\sX_H$ of \S\ref{genl-level} is separated, so it has a coarse moduli space $X_H$ (by \cite{KM97}*{1.3~(1)}, for instance). We let 
\[
Y_H \subset X_H
\]
 be the open that is the coarse moduli space of the ``elliptic curve locus'' 
 \[
\sY_H \subset \sX_H.
 \]
 We write $X(n)$, $Y_0(n)$, etc.~for $X_{\Gamma(n)}$, $Y_{\Gamma_0(n)}$, etc. 

Since $X(1) = \bP^1_\bZ$ (see \Cref{Elln-coarse}) and $X_H$ inherits $\bZ$-properness from $\sX_H$ (see \cite{Ryd13}*{6.12}), the induced map 
\[
X_H \ra X(1)
\]
is finite, so $X_H$ is a projective $\bZ$-scheme. Moreover, $X_H$ inherits normality from $\sX_H$ (see \cite{AV02}*{2.2.3} and compare with the proof of \Cref{coarse-miracle}), so $X_H \ra X(1)$ is even locally free of constant rank by \cite{EGAIV2}*{6.1.5}. In particular, $X_H$ is flat and of relative dimension $1$ over $\Spec \bZ$ at every point.
\end{pp-tweak}

Due to \Cref{Xh-sch} (and the sentence preceding it), $\sX_H = X_H$ whenever $H$ is small enough. The analysis of the case of arbitrary $H$ is facilitated by the following lemma:

\begin{lemma-tweak}[\cite{DR73}*{IV.3.10 (iii)}] \lab{coarse-quot}
For an open subgroup $H \subset \GL_2(\wh{\bZ})$ and an $n \ge 1$, if 
\[
\Gamma(n) \subset H \qq \text{and} \qq \ov{H} \ce \im(H \ra \GL_2(\bZ/n\bZ)),
\]
then $X_H$ is identified with the categorical quotient $X(n)/\ov{H}$. \QED
\end{lemma-tweak}

The coarse moduli spaces $Y_H$ and $X_H$ have been studied extensively in \cite{KM85}, albeit with somewhat different  terminology, notation, and setup. In order to put the results below in the context of the work of \cite{KM85}, we explicate the relationship between the terminology of \emph{op.~cit.}~and that of the approach based on the systematic use of the theory of algebraic stacks.

\begin{prop-tweak} \lab{KM-comp}
Let $H \subset \GL_2(\wh{\bZ})$ be an open subgroup, let $n \in \bZ_{\ge 1}$ be such that $\Gamma(n) \subset H$, and let $\ov{H}$ be the image of $H$ in $\GL_2(\bZ/n\bZ)$. 
\benum
\item \lab{KMC-a}
The ``quotient moduli problem'' $[\Gamma(n)]/\ov{H}$ in the sense of \cite{KM85}*{\S7.1} is identified with $\sY_H$.

\item \lab{KMC-b}
The ``coarse moduli scheme'' $\text{\upshape{M}}([\Gamma(n)]/\ov{H})$ in the sense of \cite{KM85}*{\S8.1} is identified with $Y_H$.

\item \lab{KMC-c}
The ``compactified coarse moduli scheme'' $\ov{\text{\upshape{M}}}([\Gamma(n)]/\ov{H})$ in the sense of \cite{KM85}*{\S8.6} is identified with $X_H$.
\eenum
\end{prop-tweak}

\bpf \hfill
\benum
\item
In the case $H = \Gamma(n)$, the identification $[\Gamma(n)] = \sY(n)$ over $\Ell$ amounts to the definitions given in \cite{KM85}*{\S5.1 and \S3.1} and \S\ref{def-Xn}, so the identification $[\Gamma(n)] = \sY_{\Gamma(n)}$ is part of \Cref{Xn-agree}. Therefore, in general, the desired identification over $\Spec \bZ[\f{1}{n}]$ results by \cite{KM85}*{7.1.3~(2)}, and hence also over all of $\Spec \bZ$ by \cite{KM85}*{7.1.3 (5)--(6)}.

\item
If $\sY_H$ is representable, then the claim follows from \ref{KMC-a} and the definition of \cite{KM85}*{8.1.1}. Therefore, in general, the claim follows from \Cref{coarse-quot}.

\item
Since $\ov{\text{\upshape{M}}}([\Gamma(n)]/\ov{H})$ is defined as the normalization of $\bP^1_{\bZ} = X(1)$ in $\text{\upshape{M}}([\Gamma(n)]/\ov{H})$, due to \ref{KMC-b} it suffices to observe that $X_H$ is the normalization of $X(1)$ in $Y_H$.
\qedhere
\eenum
\epf

Before turning to the case $H = \Gamma_0(n)$, we record the following general result that holds for every $H$.   Its part \ref{CH-a} has been proved in \cite{DR73}*{VI.6.7} by a different method, and the proof given below is in essence due to Katz and Mazur. Its part \ref{CH-b} complements \cite{KM85}*{8.5.3}.

\begin{prop-tweak} \lab{coarse-H}
Let $H \subset \GL_2(\wh{\bZ})$ be an open subgroup, and let $n \in \bZ_{\ge 1}$ be such that $\Gamma(n) \subset H$.
\benum
\item \lab{CH-a}
The coarse moduli space $(X_H)_{\bZ[\f{1}{n}]}$ of $(\sX_H)_{\bZ[\f{1}{n}]}$ is $\bZ[\f{1}{n}]$-smooth.

\item \lab{CH-b}
For any $\bZ[\f{1}{\gcd(6, n)}]$-scheme $S$, the canonical map from the coarse moduli space of $(\sX_H)_S$ to $(X_H)_S$ is an isomorphism.
\eenum
\end{prop-tweak}

\bpf 
Let $\ov{H}$ denote the image of $H$ in $\GL_2(\bZ/n\bZ)$.
\benum
\item
The coarse moduli space $X(n^2)$ may be covered by $\GL_2(\bZ/n^2\bZ)$-invariant open subschemes that are affine over $\bZ$ and are preimages of $\bZ$-affine open subschemes of $X(1)$, so \Cref{coarse-quot} and \cite{KM85}*{Theorem on p.~508 in the section ``Notes on Chapters 8 and 10''} reduce the proof to the case when $H = \Gamma(n^2)$. For this $H$, the $n = 1$ case is clear and if $n \ge 2$, then the geometric points of $\sX(n^2)_{\bZ[\f{1}{n}]}$ have no nontrivial automorphisms by \cite{KM85}*{2.7.2~(1)} and \Cref{auto-n-gon}. Thus, if $n\ge 2$, then \Cref{rep-crit}~\ref{RC-a} ensures that 
\[
\qq X(n^2)_{\bZ[\f{1}{n}]} = \sX(n^2)_{\bZ[\f{1}{n}]}
\]
and \cite{DR73}*{IV.2.5} provides the sought $\bZ[\f{1}{n}]$-smoothness of $X(n^2)_{\bZ[\f{1}{n}]}$.

\item
We work locally on $\bZ[\f{1}{\gcd(6, n)}]$, so we assume that $S$ is either a $\bZ[\f{1}{6}]$-scheme or a $\bZ[\f{1}{n}]$-scheme.

Since $\sX_H \ra \sX(1)$ is representable, the automorphism group of every geometric point of $\sX_H$ is of order dividing $24$. Therefore, by \cite{Ols06}*{2.12}, \'{e}tale locally on its coarse moduli space, $\sX_H$ is the quotient of an affine scheme $\Spec A$ by an action of a finite group $G$ whose order divides $24$. Thus, the case when $S$ is a $\bZ[\f{1}{6}]$-scheme follows from the fact that the formation of the ring of invariants $A^G$ commutes with arbitrary base change if $\#G$ is invertible in $A$.

For the remainder of the proof we assume that $S$ is a $\bZ[\f{1}{n}]$-scheme, so applying \Cref{coarse-miracle} with $\sX = (\sX_H)_{\bZ[\f{1}{n}]}$ reduces the proof to the case when $S = \Spec \bF_p$ with $p \nmid n$. We therefore let $X'$ be the coarse moduli space of $(\sX_H)_{\bF_p}$ and seek to prove that the finite map 
\[
\qq f\colon X' \ra (X_H)_{\bF_p}
\]
is an isomorphism. The source and the target curves of $f$ are $\bF_p$-smooth (equivalently, normal): the target due to \ref{CH-a} and the source due to the $\bF_p$-smoothness of $(\sX_H)_{\bF_p}$ ensured by \cite{DR73}*{IV.6.7}. Therefore, $f$ is locally free by \cite{EGAIV2}*{6.1.5}. To conclude that its rank is $1$, it suffices to exhibit a fiberwise dense open substack $\sU \subset \sY_H[\f{1}{n}]$ whose coarse moduli space is of formation compatible with base change to $\bF_p$.

We choose $\sU$ to be the preimage of the complement of $j = 0$ and $j = 1728$ in $\bA^1_{\bZ[\f{1}{n}]}$, let $\cE \ra \sU$ denote the universal elliptic curve, and let 
\[
\qq \cF \ce \ov{H} \setminus \Isom(\cE[n], (\bZ/n\bZ)^2)
\]
be the finite \'{e}tale $\sU$-stack of level $H$ structures on $\cE$ (compare with \S\ref{genl-level}). The universal level $H$-structure is a section $\gA$ of $\cF \ra \sU$, as is $[-1]_\cE^*(\gA)$. Since $\cF \ra \sU$ is finite \'{e}tale, the substack $\sV \subset \sU$ over which $\gA = [-1]_\cE^*(\gA)$ is both open and closed. By \cite{Del75}*{5.3~(III)}, the automorphism stack of $\cE$ is the constant $\{\pm 1\}_\sU$, so the open complement $\sU \setminus \sV$ is its own coarse moduli space, whereas the coarse moduli space of $\sV$ is the rigidification $\sV \!\!\! \fatslash \{\pm 1\}$ (in the notation of \cite{AOV08}*{Appendix}). Since the formation of $\sV \!\!\! \fatslash \{\pm 1\}$ commutes with arbitrary base change, so does the formation of the coarse moduli space of $\sU$.
\qedhere
\eenum
\epf

\bremstweak 
\remitweak \lab{KM-version} For a version of \Cref{coarse-H}~\ref{CH-a} in residue characteristics dividing $n$ and suitable $H$, see \cite{KM85}*{10.10.3~(5)}.

\remitweak
 In \Cref{coarse-H}~\ref{CH-b}, for some subgroups $H$ one cannot remove the requirement that $\gcd(6, n)$ be invertible on $S$. For instance, by \cite{Ces17}*{Thm.~3.2}, the canonical map from the coarse moduli space of $(\sX_{\Gamma_1(4)})_{\bF_2}$ to $(X_{\Gamma_1(4)})_{\bF_2}$ is not an isomorphism.
\eremstweak

We are ready for the promised regularity of $X_0(n)$ at the cusps. Similar techniques may be used to prove analogous regularity results for $X(n)$ or $X_1(n)$ (or even for $\wt{X}_1(n; n')$, $X_1(n; n')$, or $X_0(n; n')$ with $n$ and $n'$ as in \Cref{X1Nn-finale}), but we do not explicate them because in many cases $X(n) = \sX(n)$ and $X_1(n) = \sX_1(n)$ (see \Cref{Xn-sch} and \Cref{Xh-sch}), and in these cases the entire $X(n)$ or $X_1(n)$ is regular by \Cref{Xn-agree} or \Cref{X1-grand}~\ref{X1G-a}.

\begin{thm-tweak} \lab{coarse-reg}
For an $n \in \bZ_{\ge 1}$, the open subscheme $U \subset X_0(n)$ obtained by removing the closed points corresponding to $j = 0$ or $j = 1728$ in residue characteristics dividing $n$ is regular.
\end{thm-tweak}

\bpf
The regularity of $X_0(n)_{\bZ[\f{1}{n}]}$ follows from \Cref{coarse-H}~\ref{CH-a}, so it suffices to prove the regularity of the coarse moduli space of the preimage 
\[
\sU \subset \sX_0(n)
\]
of the open subscheme of $\bP^1_\bZ$ obtained by removing the sections $j = 0$ and $j = 1728$.

Due to the moduli interpretation of $\sX_0(n)$ given in \S\ref{X0n-def} and \Cref{X0n-main}~\ref{X0nM-a}, the constant group $\{ \pm 1\}_\sU$ is a subgroup of the automorphism group of the universal object of $\sU$. In fact, due to \cite{Del75}*{5.3~(III)} and the representability of $\sU \ra \sX(1)$, this automorphism group equals $\{ \pm 1\}_\sU$. Therefore, the coarse moduli space of $\sU$ is the rigidification $\sU\!\!\! \fatslash \{\pm 1\}$. By \cite{AOV08}*{A.1}, the map 
\[
\sU \surjects \sU\!\!\! \fatslash \{\pm 1\}
\]
is \'{e}tale, and, by \Cref{X0n-main}~\ref{X0nM-a}, the stack $\sU$ is regular, so $\sU\!\!\! \fatslash \{\pm 1\}$ is also regular, as desired.
\epf

\begin{rem-tweak} \lab{rem-more-reg}
One may use the structure of the fibers $X_0(n)_{\bF_p}$ with $p \mid n$ to sharpen \Cref{coarse-reg}. For instance, if $n$ is squarefree, then, due to \Cref{coarse-H} and \cite{KM85}*{13.5.6 and Thm.~on p.~508}, in \Cref{coarse-reg} one may require that the removed points are in addition supersingular (and likewise for general $n$ and those removed points that lie on the reduced components of $X_0(n)_{\bF_p}$). For a more thorough analysis of the coarse space $X_0(n)$, see \cite{Edi90}.
\end{rem-tweak}

We end by proving that $\sX_0(n)^\naive$ yields the same coarse moduli space $X_0(n)$, and hence suffices for many purposes (however, the proof of \Cref{coarse-reg} does rely on the finer $\sX_0(n)$ through the representability of $\sX_0(n) \ra \sX_0(1)$).

\begin{prop-tweak}
For every $n \in \bZ_{\ge 1}$, the contraction morphism
\[
\sX_0(n)^\naive \ra \sX_0(n)
\]
defined in {\upshape\S\ref{forget-naive}} induces an isomorphism on coarse moduli spaces. 
\end{prop-tweak}

\bpf
The coarse moduli space $X_0(n)'$ of $\sX_0(n)^\naive$ exists due to the finiteness of the diagonal of $\sX_0(n)^\naive$ supplied by \Cref{Con07-feast}~\ref{feast-a} (see \cite{Ryd13}*{6.12}). As in \S\ref{coarse-def}, the map 
\[
X_0(n)' \ra \bP^1_{\bZ}
\]
is finite, so, since $\sY_0(n)^\naive = \sY_0(n)$, it suffices to prove that $X_0(n)'$ is normal.

 For the normality, we work Zariski locally on $X_0(n)'$ and note that each open substack 
 \[
 \sU \subset \sX_0(n)^\naive
 \]
 that has an affine coarse moduli space $\Spec A$ satisfies $A = \Gamma(\sU, \sO_\sU)$ by the universal property for maps to $\bA^1_\bZ$. To then see that $\Gamma(\sU, \sO_\sU)$ is integrally closed in its total ring of fractions it suffices to use the normality of $\sU$ supplied by \Cref{Con07-feast}~\ref{feast-a} and the fact that generizations lift along smooth morphisms from algebraic spaces to $\sU$ (see \cite{LMB00}*{5.7.1}).
\epf

\begin{rem-tweak}
The same proof shows that, in the notation of section \ref{Gamma-1-Nn}, for every $n, n' \in \bZ_{\ge 1}$ the coarse moduli spaces of $\sX_1(n; n')$ and $\sX_0(n; n')$ agree with those of $\sX_{\Gamma_1(n; n')}$ and $\sX_{\Gamma_0(n; n')}$.
\end{rem-tweak}


\end{document}